\documentclass[11pt]{amsart}

\usepackage[english]{babel}
\usepackage[utf8x]{inputenc}
\usepackage[T1]{fontenc}

\usepackage[top=3cm,bottom=2cm,left=3cm,right=3cm,marginparwidth=1.75cm]{geometry}
\usepackage{lipsum}
\usepackage{authblk}
\usepackage{blindtext}
\usepackage{amsmath, amsthm, amssymb, esint, mathtools, mathabx}
\usepackage{amsfonts}
\usepackage{graphicx}
\graphicspath{{/texfiles/}}
\usepackage[colorinlistoftodos]{todonotes}
\usepackage[colorlinks=true, allcolors=blue]{hyperref}
\usepackage{multicol}
\usepackage[foot]{amsaddr}

\theoremstyle{thmstyletwo}%
\newtheorem{theorem}{Theorem}[section]

\newtheorem{lemma}[theorem]{Lemma}%
\newtheorem{remark}{Remark}[section]
\newtheorem{proposition}[theorem]{Proposition}
\newtheorem{corollary}[theorem]{Corollary}
\newtheorem{claim}{Claim}[section]%

\numberwithin{equation}{section}

\begin{document}

\author[A. Martina Neuman]{A. Martina Neuman}
\address{Department of Mathematics, University of California, Berkeley}

\title{Functions of nearly maximal Gowers-Host-Kra norms on Euclidean spaces}

\begin{abstract} 
Let $k\geq 2, n\geq 1$ be integers. Let $f: \mathbb{R}^{n} \to \mathbb{C}$. The $k$th Gowers-Host-Kra norm of $f$ is defined recursively by 
\begin{equation*} 
    \| f\|_{U^{k}}^{2^{k}} =\int_{\mathbb{R}^{n}} \| T^{h}f \cdot \bar{f} \|_{U^{k-1}}^{2^{k-1}} \, dh 
\end{equation*}
with $T^{h}f(x) = f(x+h)$ and $\|f\|_{U^1} = | \int_{\mathbb{R}^{n}} f(x)\, dx |$. These norms were introduced by Gowers in his work on Szemer\'edi's theorem, and by Host-Kra in ergodic setting. 
It's shown by Eisner and Tao that for every $k\geq 2$ there exist $A(k,n)< \infty$ and $p_{k} = 2^{k}/(k+1)$ such that $\| f\|_{U^{k}} \leq A(k,n)\|f\|_{p_{k}}$, for all $f \in L^{p_{k}}(\mathbb{R}^{n})$. The optimal constant $A(k,n)$ and the extremizers for this inequality are known. In this exposition, it is shown that if the ratio $\| f \|_{U^{k}}/\|f\|_{p_{k}}$ is nearly maximal, then $f$ is close in $L^{p_{k}}$ norm to an extremizer.
\end{abstract}


\maketitle

\tableofcontents

\section{Introduction} \label{sec:intro}

\noindent Let $k\geq 2, n\geq 1$ be integers and $f:\mathbb{R}^{n}\to\mathbb{C}$ be a measurable function. The $k$th Gowers-Host-Kra norm of $f$ is defined recursively by:
\begin{equation} \label{def:Uk}
    \| f\|_{U^{k}}^{2^{k}} = \int_{\mathbb{R}^{n}} \| T^{h} f \cdot\bar{f} \|_{U^{k-1}}^{2^{k-1}} \, dh.
\end{equation} 
Here, $T^{h}f(x) = f(x+h)$, and $\|f\|_{U^1} = | \int_{\mathbb{R}^{n}} f(x)\, dx|$. These norms were introduced by Gowers \cite{gowers2001new} in his work on Szemer\'edi's theorem, and by Host-Kra \cite{host2005nonconventional} in ergodic setting. They are also discussed in \cite{tao2006additive}. There is an alternative expression to \eqref{def:Uk} \cite{tao2006additive}: 
\begin{equation} \label{def:Ukmult}
    \|f\|_{U^{k}}^{2^{k}} = \int_{\mathbb{R}^{(k+1)n}} \prod_{\alpha\in\{0,1\}^{k}} \mathcal{C}^{\omega_{\alpha}}f (x+\alpha\cdot\vec{h})\, dxd\vec{h}.
\end{equation}
Here $\alpha = (\alpha_1, \text{ ... }, \alpha_{k})\in\{0,1\}^{k}$, $\vec{h} = (h_1, \text{ ... }, h_{k})\in\mathbb{R}^{kn}$ and $\alpha\cdot\vec{h} = \sum_{i=1}^{k} \alpha_{i}h_{i}\in\mathbb{R}^{n}$. $\mathcal{C}f=\bar{f}$ is the conjugation operator and $\omega_{\alpha}$ is the parity of the number of ones in $\alpha\in\{0,1\}^{k}$. For instance, if $k=2$, then \eqref{def:Ukmult} becomes
\begin{equation*} 
    \|f\|_{U^2}^4 = \int_{\mathbb{R}^{3n}} f(x)\overline{f(x+h_1)}\overline{f(x+h_2)}f(x+h_1+h_2)\, dxdh_1dh_2.
\end{equation*}
It follows that, after some changes in variables, 
\begin{align} 
    \nonumber \|f\|_{U^2}^4 &= \int_{\mathbb{R}^{3n}} f(x)\overline{f(y)}\overline{f(z)}f(x+y-z)\, dxdydz = \int_{\mathbb{R}^{2n}} f(x)\overline{f(y)} (f\ast\overline{f})(x+y)\, dxdy \\ 
    \label{U2ex} &= \int_{\mathbb{R}^{2n}} f(x)\overline{f(y-x)}(f\ast\overline{f})(y) \,dxdy = \int_{\mathbb{R}^{n}} (f\ast\overline{f})^2(y)\, dy = \int_{\mathbb{R}^{n}} |\hat{f}|^4(x)\, dx
\end{align}
or $\|f\|_{U^2} = \|\hat{f}\|_4$. Hence the second Gowers-Host-Kra norm is an actual norm. The same has been shown for higher order Gowers-Host-Kra norms \cite{tao2006additive}. It's also shown in \cite{eisner2012large} that for every $k\geq 2, n\geq 1$ there exist $A(k,n) = A(k,1)^{n}$ with $A(k,1) =2^{k/2^{k}}/(k+1)^{(k+1)/2^{k+1}}$ and $p_{k} = 2^{k}/(k+1)$ such that, for every $f \in L^{p_{k}}(\mathbb{R}^{n})$, 
\begin{equation} \label{GHKineq}
    \|f\|_{U^{k}} \leq A(k,n) \|f\|_{p_{k}}.
\end{equation}
We call $f$ an extremizer of \eqref{GHKineq} if equality is attained with this choice of $f$. See Section \ref{sec:basics} below for terminology used in this exposition. It's shown in \cite{eisner2012large} that this constant $A(k,n)$ is best possible and that equality is attained iff $f(x) = C\exp(-(x-c)\cdot M(x-c)) \exp(2\pi iP(x))$, with $C\in\mathbb{C}$, $c\in\mathbb{R}^{n}$, $M$ is a positive-definite $n\times n$ matrix and $P: \mathbb{R}^{n}\to\mathbb{R}$ is a polynomial of degree at most $k-1$. In this exposition, we seek to answer the question, of what happens if the ratio $\|f\|_{U^{k}}/\|f\|_{p_{k}}$ is nearly maximized on $\mathbb{R}^{n}$. This question was not addressed in \cite{eisner2012large}. We set out to show that $f$ must then be close in $L^{p_{k}}$ norm to an extremizer. 

\begin{theorem} \label{thm:main}
Let $k\geq 2, n\geq 1$ be integers. For every $\epsilon > 0$ there exists $\delta > 0$ with the following property. Suppose $f\in L^{p_{k}}(\mathbb{R}^{n})$ and 
\begin{equation} \label{thmstate}
    \| f \|_{U^{k}} \geq (1-\delta)A(k,n)\|f\|_{p_{k}}.
\end{equation}
Then there exists an extremizer $\mathcal{F}$ of \eqref{GHKineq} on $\mathbb{R}^{n}$ such that $\| \mathcal{F}-f\|_{p_{k}} \leq \epsilon\|f\|_{p_{k}}$. 
\end{theorem}

\noindent Theorem \ref{thm:main} can also be equivalently stated in terms of sequences.

\begin{theorem} \label{thm:second}
Let $k\geq 2, n\geq 1$ be integers. Suppose $\{f_{j}\}_{j}$ is a sequence of functions on $\mathbb{R}^{n}$ such that $\|f_{j}\|_{p_{k}} = 1$ and
\begin{equation*} 
    \|f_{j}\|_{U^{k}}\geq (1-\delta_{j})A(k,n)\|f_{j}\|_{p_{k}} = (1-\delta_{j})A(k,n)
\end{equation*}
with $\delta_{j}\to 0$ as $j\to\infty$. Then there exist an extremizer $\mathcal{F}$ of \eqref{GHKineq} on $\mathbb{R}^{n}$, $\lambda_{j}\in\mathbb{R}_{>0}, c_{j}\in\mathbb{R}^{n}$, real-valued polynomial of degree at most $k-1$, $P_{j}$, and positive definite $M_{j}\in M^{n\times n}(\mathbb{R})$ such that 
\begin{equation*}
    \|\mathcal{F} - \lambda_{j}^{n/p_{k}} f_{j}e^{iP_{j}}(\lambda_{j}(M_{j}(\cdot - c_{j})))\|_{p_{k}}\leq o(\delta_{j})\|f_{j}\|_{p_{k}} = o(\delta_{j}).
\end{equation*}
\end{theorem}

\noindent
Theorem \ref{thm:main} and Theorem \ref{thm:second} are equivalent to each other for the following reasons. Note that \eqref{thmstate} is invariant under an orthogonal action on the underlying Euclidean space $\mathbb{R}^{n}$ as well as translation and dilation. Now \eqref{thmstate} is also unchanged if $f$ is replaced by $fe^{iQ}$, with $Q$ being a real-valued function, which will be shown in Section \ref{sec:complexcase} to be "nearly" polynomial, in the sense defined there. Conversely, if $\{f_{j}\}_{j}$ is a sequence as in Theorem \ref{thm:second}, then $\{|f_{j}|\}_{j}$ must first be appropriately translated, dilated and transformed so that all $|f_{j}|$ can be close to the same nonnegative extremizer of \eqref{GHKineq} in $L^{p_{k}}$ norm. If the $f_{j}$ are complex-valued then their complex phases of $f_{j}$ must also be adjusted, by multiplying with appropriate polynomial phases. These modifications are done in order to introduce compactness to such sequence $\{f_{j}\}_{j}$ - and these are the only modifications required, based on the description of extremizers for \eqref{GHKineq}. Hence, we will not distinguish the two theorems in this discussion and simply refer to both as Theorem \ref{thm:main}. Via the relation between the second Gowers-Host-Kra norm of $f$ and the $L^4$ norm of $\hat{f}$ shown in \eqref{U2ex}, this near extremization question has been resolved in \cite{christ2014sharpened}, as $p_2 = 4/3$, and hence $\|f\|_{U^2}\geq (1-\delta)A(2,n)\|f\|_{4/3}$ becomes $\|\hat{f}\|_4\geq (1-\delta)A(2,n)\|f\|_{4/3}$.\\

\noindent
The alternative expression \eqref{def:Ukmult} of the Gowers-Host-Kra norm leads to an inner product inequality. The Gowers inner product \cite{tao2006additive} of degree $2^{k}$, $\mathcal{T}_{k}$, is defined as follows. For each $\alpha\in\{0,1\}^{k}$, let $f_{\alpha}:\mathbb{R}^{n}\to\mathbb{C}$ be a measurable function. Then 
\begin{equation*}   
    \mathcal{T}_{k}(f_{\alpha}:\alpha\in\{0,1\}^{k})= \int_{\mathbb{R}^{(k+1)n}} \prod_{\alpha\in\{0,1\}^{k}} \mathcal{C}^{\omega_{\alpha}}f_{\alpha}(x+\alpha\cdot\vec{h}) \, dxd\vec{h}.
\end{equation*}
The mentioned inner product inequality is the following
\begin{equation} \label{innineq}
    |\mathcal{T}_{k}(f_{\alpha}:\alpha\in\{0,1\}^{k})|\leq A(k,n)^{2^{k}}\prod_{\alpha\in\{0,1\}^{k}}\|f_{\alpha}\|_{p_{k}};
\end{equation}
see Section \ref{sec:basics} below. Observe that \eqref{GHKineq} is a consequence of \eqref{innineq} if $f_{\alpha} = f$ for all $\alpha\in\{0,1\}^{k}$. A more general version of \eqref{thmstate} is then
\begin{equation*} 
    |\mathcal{T}_{k}(f_{\alpha}:\alpha\in\{0,1\}^{k})|\geq (1-\delta) A(k,n)^{2^{k}}\prod_{\alpha\in\{0,1\}^{k}}\|f_{\alpha}\|_{p_{k}}.
\end{equation*}
We will not be investigating the near extremization question in this full generality; however, the arguments presented in the second part of this exposition, which spans from Section \ref{sec:localization} to Section \ref{sec:complexcase}, can be extended for the said case. The conclusion of Theorem \ref{thm:main} is qualitative; the dependence of $\epsilon$ on $\delta$ is not made explicit in this exposition. This non-quantitative dependence $\epsilon = \epsilon(\delta)$ is a result of implicitness arisen in Proposition \ref{prop:lwbrearr} and Subsection \ref{sec:alignment}; other steps could be made explicit without this implicitness. There will be two parts in this exposition. The first part, Section \ref{sec:shortpf}, is a very short proof of Theorem \ref{thm:main} in the case of nonnegative functions, which utilizes a stability result of Young's convolution inequality in \cite{christ2019near}. The second part, spanning from Section \ref{sec:localization} to Section \ref{sec:complexcase}, concludes the theorem for general complex-valued functions and does not rely on the said stability result. The framework used in the second part parallels that laid out in \cite{christ2019near}. Another such inequality that is closely related to the Gowers product inequality is the following:
\begin{equation*}    
    |\mathcal{T}_{k}(f_{\alpha}:\alpha\in\{0,1\}^{k})| \leq C(k,n,\vec{p})\prod_{\alpha\in\{0,1\}^{k}}\|f_{\alpha}\|_{p_{\alpha}}
\end{equation*}
with $\vec{p} = (p_{\alpha}:\alpha\in\{0,1\}^{k})$. The $p_{\alpha}$ satisfy the necessary scaling condition $\sum_{\alpha\in\{0,1\}^{k}}p_{\alpha}^{-1} = k+1$. In Section \ref{sec:localization} below, it's shown that the said inequality holds (in other words, $C(k,n,\vec{p}) <\infty$) if $\vec{p}$ lies in a small neighborhood of the point $\vec{P} = (2^{k}/(k+1),\text{ ... }, 2^{k}/(k+1))$ on the hyperplane $\sum_{\alpha\in\{0,1\}^{k}} p_{\alpha}^{-1} = k+1$. No general sufficient algebraic conditions on $p_{\alpha}$ are known as of now; however, see \cite{bennett2008brascamp} for a geometric description. This discovery at least confirms that the characteristic polytope, as discussed in \cite{bennett2008brascamp}, of the Gowers inner product structure - which is a Brascamp-Lieb structure - is not a degenerate one-point set. The Gowers-Host-Kra norm inequality is an important step toward the study of the near extremization problems of these more general cases. 

\section{Basics} \label{sec:basics}

\noindent
The notation $\mathbb{N}$ denotes the set of nonnegative integers. The notation $x, y$ or $h$ indicates an element in a Euclidean space while the notation $\vec{x},\vec{y}$ or $\vec{h}$ indicates a tuple of such elements. In particular, the notation $(x,\vec{h})$ always means $(x,\vec{h})\in\mathbb{R}^{(k+1)n}$, with $x\in\mathbb{R}^{n}$, $\vec{h} = (h_1,\text{ ... }, h_{k})$ and $h_{i}\in\mathbb{R}^{n}$, $i\in\{1,\text{ ... }, k\}$. The notation $crd(S)$ denotes the cardinality of a finite set $S$. The notation $\mathcal{L}(S)$ denotes the $n$-dimensional Lebesgue measure of a Euclidean set $S\in\mathbb{R}^{n}$. If $S$ is a ball with center $c$ and $r > 0$, then $rS$ denotes the dilated set $\{r(x-c) \, | \, x \in S\}$. The notation $spt(f)$ denotes the essential support of a measurable function $f$ on a Euclidean space. If $f$ is continuous then this support is taken to be the closure of the set of points at which $f$ does not vanish. The symbol $\cdot$ denotes several meanings of multiplication, which will be clear from context. The notation $\langle a, b\rangle_{\mathbb{R}^{n}}$ denotes the usual real inner product between two Euclidean elements $a,b$.\\

\noindent
Let $a,b >0$. The symbol $\asymp$ in $a\asymp b$ means the following: There exist $C >c > 0$ such that $c a \leq b \leq Ca$. In our context, $c, C$ will be integral powers of $2$, say $2^{j}a \leq b \leq 2^{l}a$; in some cases, $j$ and $l$ are chosen so that $2^{j}a$ and $2^{l}a$ are among the closest numbers of this type to $b$.\\
Let $S\subset\mathbb{R}^{n}$. Then 
\begin{equation*} 
    \tilde{S}^{k+1} = \{(x,\vec{h})\in\mathbb{R}^{(k+1)n} \, | \, x, h_{i}, x+\alpha\cdot\vec{h}\in S, \forall i\in\{1,\text{ ... }, k\}, \forall\alpha\in\{0,1\}^{k}\}.
\end{equation*}

\noindent By "$f$ is of unit norm" we mean $\|f\|_{p} = 1$, if the $L^{p}$ norm is understood. By "a Gaussian on $\mathbb{R}^{n}$" we mean the function, $m\exp(-a|x-c|^2)$, with $m, a\in\mathbb{R}_{>0}$ and $c\in\mathbb{R}^{n}$. By "centered Gaussian" we mean $m\exp(-a|x|^2)$. By "the standard centered Gaussian" we mean $m\exp(-|x|^2)$, with $m>0$ chosen so that $\|\mathcal{G}\|_{p} = 1$, for some appropriate $L^{p}$ norm. By "measurable" we mean "Lebesgue measurable." By "a centered ball" we mean "a ball that is centered at the origin." Finally, we employ some short-hand notations - for instance, if $f:\mathbb{R}^{n}\to\mathbb{R}$ then $\{ f > \lambda\}$ means the super-level set $\{ x\in\mathbb{R}^{n}: f(x) > \lambda\}$ - and we make an effort to indicate parametric dependence of our constants whenever it's important for our calculations. When the dimension $n=1$, we simply write $C(\text{other parameters}, 1) = C(\text{other parameters})$.\\
 
\noindent
We call \eqref{GHKineq} "the $k$th Gowers-Host-Kra norm inequality" and \eqref{innineq} "the $k$th Gowers inner product inequality." Let $f_{\alpha}, f:\mathbb{R}^{n}\to\mathbb{C}$ be measurable functions, $\alpha\in\{0,1\}^{k}$, and $\vec{f} = (f_{\alpha}:\alpha\in\{0,1\}^{k})$. We say $f$ is an extremizer of the $k$th Gowers-Host-Kra norm inequality if:
\begin{equation*} 
    \| f\|_{U^{k}} = A(k,n) \|f\|_{p_{k}}. 
\end{equation*}
We say $f$ is a $(1-\delta)$ near extremizer of the $k$th Gowers-Host-Kra norm inequality if:
\begin{equation*} 
    \| f\|_{U^{k}} \geq (1-\delta)A(k,n) \|f\|_{p_{k}}
\end{equation*}
and we say $\{f_{i}\}_{i}$ is an extremizing sequence of the $k$th Gowers-Host-Kra norm inequality if $\|f_{i}\|_{p_{k}} =1$ and,
\begin{equation*} 
    \|f_{i}\|_{U^{k}}\geq (1-\delta_{i})A(k,n)\|f_{i}\|_{p_{k}} = (1-\delta_{i})A(k,n)
\end{equation*}
with $\delta_{i}\to 0$ as $i\to\infty$. We say $\vec{f}$ is an extremizing tuple of the $k$th Gowers inner product inequality if:
\begin{equation*} 
    \mathcal{T}_{k}(\vec{f}) = \mathcal{T}_{k}(f_{\alpha}:\alpha\in\{0,1\}^{k}) = A(k,n)^{2^{k}} \prod_{\alpha\in\{0,1\}^{k}} \|f_{\alpha}\|_{p_{k}}.
\end{equation*}
We say $\vec{f}$ is a $(1-\delta)$ near extremizing tuple of the $k$th Gowers inner product inequality if:
\begin{equation*} 
    \mathcal{T}_{k}(\vec{f}) = \mathcal{T}_{k}(f_{\alpha}:\alpha\in\{0,1\}^{k}) \geq (1-\delta) A(k,n)^{2^{k}} \prod_{\alpha\in\{0,1\}^{k}} \|f_{\alpha}\|_{p_{k}}.
\end{equation*}

\noindent
\textbf{Gowers-Cauchy-Schwarz inequality implies Gowers product inequality:}
\begin{equation*} 
    \mathcal{T}_{k}(f_{\alpha}:\alpha\in\{0,1\}^{k})\leq\prod_{\alpha\in\{0,1\}^{k}}\|f_{\alpha}\|_{U^{k}}\leq A(k,n)^{2^{k}}\prod_{\alpha\in\{0,1\}^{k}}\|f_{\alpha}\|_{p_{k}}.
\end{equation*}
The first inequality is the Gowers-Cauchy-Schwarz inequality \cite{tao2006additive}; the second is the Gowers-Host-Kra norm inequality. We called the resulting inequality, the Gowers product inequality.\\

\noindent
\textbf{Relation of Gowers-Host-Kra norm inequality to sharp Young's inequality:} \\
Specializing the version of the sharp Young's inequality given in \cite{eisner2012large} to the exponents $s= 2^{k-1}$, $r=p_{k-1}$ and $t=q=p_{k}$, we have:
\begin{equation*} 
    \left (\int_{\mathbb{R}^{n}} \| (T^{h}f)\bar{g}\|_{L^{p_{k-1}}(\mathbb{R}^{n})}^{2^{k-1}}\, dh\right )^{1/2^{k-1}}\leq \left (\frac{C_{2k/(k+1)}^2}{C_{k}} \right)^{nk/2^{k-1}} \|f\|_{L^{p_{k}}(\mathbb{R}^{n})}\|g\|_{L^{p_{k}}(\mathbb{R}^{n})}
\end{equation*}
with $C_{t} = \left (\frac{t^{1/t}}{t'^{1/t'}} \right)^{1/2}$, $1/t + 1/t'=1$. It's shown in \cite{eisner2012large} that $A(k)\leq 1 = A(1)$ and that $A(k) = \left ( \frac{C_{2k/(k+1)}^2}{C_{k}} \right )^{k/2^{k}} C_{k-1}^{1/2}$ for all $k\geq 1$. Hence:
\begin{equation*}
    \left ( \int_{\mathbb{R}^{n}} \| (T^{h}f)\bar{g}\|_{L^{p_{k-1}}(\mathbb{R}^{n})}^{2^{k-1}}\, dh \right )^{1/2^{k-1}} \leq \left (\frac{A(k,n)}{A(k-1,n)^{1/2}} \right)^2 \|f\|_{L^{p_{k}}(\mathbb{R}^{n})}\|g\|_{L^{p_{k}}(\mathbb{R}^{n})}.
\end{equation*}
Now letting $f=g$, incorporating this with the ($k-1$)th Gowers-Host-Kra norm inequality and using the definition of the Gowers-Host-Kra norms, we obtain:
\begin{multline*} 
    \|f\|_{U^{k}(\mathbb{R}^{n})}^2 = \left (\int_{\mathbb{R}^{n}}\|(T^{h}f)\bar{f}\|_{U^{k-1}(\mathbb{R}^{n})}^{2^{k-1}}\, dh\right)^{1/2^{k-1}}\\ \leq A(k-1,n)\left (\frac{A(k,n)}{A(k-1,n)^{1/2}}\right)^2\|f\|_{L^{p_{k}}(\mathbb{R}^{n})}^2 = A(k,n)^2\|f\|^2_{L^{p_{k}}(\mathbb{R}^{n})}
\end{multline*}
which is the $k$th Gowers-Host-Kra norm inequality. Due to this relationship between the two inequalities, the optimal constants $A(k,n)$ satisfy: $A(k,m\cdot n) = (A(k,m))^{n}$; in particular, $A(k,m)= (A(k))^{m}$.\\

\noindent
\textbf{Facts about Lorentz semi-norms:} \cite{stein1971introduction} \\
Let $f:\mathbb{R}^{n}\to\mathbb{R}_{\geq 0}$ be a measurable function. There exists a unique decomposition $f = \sum_{j\in\mathbb{Z}} 2^{j} F_{j}$ such that $1_{\mathcal{F}_{j}} \leq F_{j}< 2\cdot 1_{\mathcal{F}_{j}}$, with the measurable sets $\mathcal{F}_{j}$ being pairwise disjoint up to null sets. This decomposition comes from the layer cake representation: $f(x) = \int_0^{\infty} 1_{\{f > t\}}(x)\, dt = \sum_{j\in\mathbb{Z}} \int_{2^{j}}^{2^{j+1}} 1_{\{f > t\}}(x)\, dt =\sum_{j\in\mathbb{Z}} 2^{j} F_{j}(x)$, with $F_{j}(x) = 2^{-j}\int_{2^{j}}^{2^{j+1}} 1_{\{f\geq t\}}(x)\, dt$. It's readily checked that $1_{\mathcal{F}_{j}}\leq F_{j}< 2\cdot 1_{\mathcal{F}_{j}}$, with $\mathcal{F}_{j} = \{ 2^{j}\leq f < 2^{j+1}\}$. Let $f_{*}$ be the nonincreasing rearrangement of $f$ \cite{stein1971introduction}, defined as follows. If $\lambda$ denotes the distribution function of $f$, then $f_{*}(t) = \inf\{s:\lambda(s)\leq t\}$. We claim that:

\begin{lemma} Let $1\leq q<\infty, 1\leq\tilde{q}\leq\infty$. Then:
\begin{equation*} 
    \| f\|_{(q,\tilde{q})} = \left ( \frac{\tilde{q}}{q} \int_0^{\infty} (t^{1/q}f_{*}(t))^{\tilde{q}} \, dt \right )^{1/\tilde{q}} \asymp \left (\sum_{j\in\mathbb{Z}} 2^{j\tilde{q}} (\mathcal{L}(\mathcal{F}_{j}))^{\tilde{q}/q}\right )^{1/\tilde{q}}
\end{equation*}
if $\tilde{q}<\infty$, and, 
\begin{equation*} 
    \|f\|_{(q,\infty)} = \sup_{t > 0} t^{1/q}f_{*}(t) \asymp \sup_{j:\mathcal{F}_{j}\not=\emptyset} 2^{j} (\mathcal{L}(\mathcal{F}_{j}))^{1/q}.
\end{equation*}
\end{lemma}

\begin{proof} First note that if $f=1_{E}$, $E$ being a measurable set of $\mathbb{R}^{n}$, then $\|f\|_{(q,\tilde{q})} = (\mathcal{L}(E))^{1/q}$, $1\leq\tilde{q}\leq\infty$. Since $1_{\mathcal{F}_{j}} \leq F_{j} < 2\cdot 1_{\mathcal{F}_{j}}$, we first consider the case $f = \sum_{j=-l}^{l} 2^{j} 1_{E_{j}}$ with pairwise disjoint measurable sets $E_{j}$ of $\mathbb{R}^{n}$ and $l\in\mathbb{Z}_{>0}$. Then:
\begin{equation*} 
    \|f\|_{(q,\tilde{q})} = \left (\sum_{j=-l}^{l} 2^{j\tilde{q}}(B_{j}^{\tilde{q}/q}-B_{j+1}^{\tilde{q}/q}) \right)^{1/\tilde{q}} \asymp \left (\sum_{j=-l}^{l} 2^{j\tilde{q}} (\mathcal{L}(E_{j}))^{\tilde{q}/q} \right )^{1/\tilde{q}}.
\end{equation*} 
Here $B_{j} = \sum_{i=j}^{l} (\mathcal{L}(E_{i}))$ and $B_{l+1} = 0$. The first equality is easily deduced from the case of a single set indicator function - see also \cite{stein1971introduction}. The approximation that follows is an application of the algebraic fact: $(a^{t}+b^{t})\leq (a+b)^{t}\leq 2^{t} (a^{t}+b^{t})$, for $t\geq 1$ and $a,b>0$. That means, the constants in the approximation depend only on $q,\tilde{q}$ and not on $l$. Then for the general case:
\begin{equation*} 
    \| f\|_{(q,\tilde{q})} = \| \sum_{j\in\mathbb{Z}} 2^{j}F_{j}\|_{(q,\tilde{q})} =\lim_{l\to\infty} \| \sum_{j=-l}^{l} 2^{j} F_{j}\|_{(q,\tilde{q})}\asymp \lim_{l\to\infty} \left ( \sum_{j=-l}^{l} 2^{j\tilde{q}} (\mathcal{L}(\mathcal{F}_{j}))^{\tilde{q}/q}\right)^{1/\tilde{q}}
\end{equation*}
for all sufficiently positively large $l$ depending on $f$. Hence this ultimately allows us to write:
\begin{equation*} 
    \| f\|_{(q,\tilde{q})}\asymp \lim_{l\to\infty} \left ( \sum_{j=-l}^{l} 2^{j\tilde{q}} (\mathcal{L}(\mathcal{F}_{j}))^{\tilde{q}/q}\right)^{1/\tilde{q}}=\left ( \sum_{j\in\mathbb{Z}} 2^{j\tilde{q}} (\mathcal{L}(\mathcal{F}_{j}))^{\tilde{q}/q}\right)^{1/\tilde{q}}.
\end{equation*}
Similarly, if $\tilde{q}=\infty$ and if $f=\sum_{j=-l}^{l} 2^{j} 1_{E_{j}}$, then $\|f\|_{(q,\infty)} = \sup_{-l\leq j\leq l} 2^{j}B_{j}^{1/q}\asymp \sup_{-l\leq j\leq l} 2^{j}(\mathcal{L}(E_{j}))^{1/q}$. Then for the general case:
\begin{align*} 
    \| f\|_{(q,\infty)} = \| \sum_{j\in\mathbb{Z}} 2^{j}F_{j}\|_{(q,\infty)} &=\lim_{l\to\infty} \| \sum_{j=-l}^{l} 2^{j} F_{j}\|_{(q,\infty)}\\ &\asymp\lim_{l\to\infty}\left ( \sup_{-l\leq j\leq l} 2^{j} (\mathcal{L}(\mathcal{F}_{j}))^{1/q}\right ) \asymp \sup_{j:\mathcal{F}_{j}\not=\emptyset} 2^{j}(\mathcal{L}(\mathcal{F}_{j}))^{1/q}.
\end{align*}
The last approximation follows because the first approximation holds for all sufficiently positively large $l$. 
\end{proof}

\section{A short proof of Theorem \ref{thm:main} for nonnegative functions} \label{sec:shortpf}

\noindent
\cite{personalcommunication} Let $\beta\in\{0,1\}^{k+1}$, then $\beta = (\alpha,0)$ or $\beta=(\alpha,1)$ for some $\alpha\in\{0,1\}^{k}$. Let $\vec{f} = (f_{\beta}:\beta\in\{0,1\}^{k+1})$, and suppose $f_{\beta}\geq 0$ for all $\beta\in\{0,1\}^{k+1}$. Then
\begin{align}
    \nonumber \mathcal{T}_{k+1}(f_{\beta}:\beta\in\{0,1\}^{k+1}) &= \int_{\mathbb{R}^{n}} \mathcal{T}_{k}(T^{h}f_{(\alpha,0)}\cdot f_{(\alpha,1)}:\alpha\in\{0,1\}^{k}) \, dh \\
    \nonumber &\leq \int_{\mathbb{R}^{n}} \prod_{\alpha\in\{0,1\}^{k}} \| T^{h} f_{(\alpha,0)}\cdot f_{(\alpha,1)}\|_{U^{k}}\, dh\\ \nonumber &\leq A(k,n)^{2^{k}}\int_{\mathbb{R}^{n}} \prod_{\alpha\in\{0,1\}^{k}} \| T^{h}f_{(\alpha,0)}\cdot f_{(\alpha,1)}\|_{p_{k}} \, dh\\
    \label{shortpf1} &\leq A(k,n)^{2^{k}} \prod_{\alpha\in\{0,1\}^{k}} \left ( \int_{\mathbb{R}^{n}} \|T^{h}f_{(\alpha,0)}\cdot f_{(\alpha,1)}\|_{p_{k}}^{2^{k}}\, dh \right )^{1/2^{k}}.
\end{align}
The first inequality is due to the Gowers-Cauchy-Schwarz inequality, the second to the Gowers-Host-Kra norm inequality and the last to H\"older's inequality. Substituting $p_{k} = 2^{k}/(k+1)$ in \eqref{shortpf1}:
\begin{align} 
    \nonumber \int_{\mathbb{R}^{n}} \| T^{h}f_{(\alpha,0)}\cdot f_{(\alpha,1)}\|_{p_{k}}^{2^{k}}\, dh &= \int_{\mathbb{R}^{n}} \left ( \int_{\mathbb{R}^{n}} f_{(\alpha,0)}^{p_{k}}(y+h)f_{(\alpha,1)}^{p_{k}}(y) \, dy \right )^{2^{k}/p_{k}}\, dh\\
    \nonumber & = \int_{\mathbb{R}^{n}} \left ( \int_{\mathbb{R}^{n}} f_{(\alpha,0)}^{p_{k}}(y+h)f_{(\alpha,1)}^{p_{k}}(y) \, dy \right )^{k+1}\, dh =\| \tilde{f}^{p_{k}}_{(\alpha,0)}\ast f_{(\alpha,1)}^{p_{k}}\|_{k+1}^{k+1}\\
    \label{shortpf2} & \leq Y(k+1,n)^{k+1}\|\tilde{f}_{(\alpha,0)}^{p_{k}}\|_{q}^{k+1}\|f^{p_{k}}_{(\alpha,1)}\|_{q}^{k+1}
    =Y(k+1, n)^{k+1}\|f_{(\alpha,0)}\|_{p_{k}q}^{2^{k}}\|f^{p_{k}}_{(\alpha,1)}\|_{p_{k}q}^{2^{k}}.
\end{align}
The sole inequality in the display \eqref{shortpf2} above follows from Young's convolution inequality. Here $\tilde{f}(y) = f(-y)$, $2q^{-1} = (k+1)^{-1}+1$, or equivalently, $q = p_{k+1}/p_{k}$ and $Y(k+1,n)$ is the optimal constant in Young's convolution inequality for the involved exponents. Combining \eqref{shortpf1} and \eqref{shortpf2}, we obtain the following majorization:
\begin{equation*}   
    \mathcal{T}_{k+1}(f_{\beta}:\beta\in\{0,1\}^{k+1}) \leq A(k,n)^{2^{k}} Y(k+1,n)^{k+1}\prod_{\beta\in\{0,1\}^{k+1}} \|f_{\beta}\|_{p_{k+1}} = A(k+1,n)^{2^{k+1}}\prod_{\beta\in\{0,1\}^{k+1}} \|f_{\beta}\|_{p_{k+1}}
\end{equation*}
as it happens that $A(k,n)^{2^{k}}Y(k+1,n)^{k+1} = A(k+1,n)^{2^{k+1}}$ - see \cite{eisner2012large} and Section \ref{sec:basics}. Now suppose that there exists $\delta >0$ such that
\begin{equation*} 
    \mathcal{T}_{k+1}(f_{\beta}:\beta\in\{0,1\}^{k+1}) \geq (1-\delta) A(k+1,n)^{2^{k+1}}\prod_{\beta\in\{0,1\}^{k+1}}\|f_{\beta}\|_{p_{k+1}}.
\end{equation*}
Then each inequality in \eqref{shortpf1} and \eqref{shortpf2} must hold in reverse, up to a factor of $1-c(k)\delta$, for some small $c(k)>0$. In particular, for each $\alpha\in\{0,1\}^{k}$,
\begin{equation} \label{alpha0}
    \| \tilde{f}^{p_{k}}_{(\alpha,0)}\ast f^{p_{k}}_{(\alpha,1)}\|_{k+1}^{k+1}\geq (1-c(k)\delta)Y(k+1,n)^{k+1}\|f^{p_{k}}_{(\alpha,0)}\|_{q}^{k+1}\|f_{(\alpha,1)}^{p_{k}}\|_{q}^{k+1}.
\end{equation}
Since $(\alpha,0),(\alpha,1)\in\{0,1\}^{k+1}$ if $\alpha\in\{0,1\}^{k}$, a stability result for Young's convolution inequality in \cite{christ2019near} implies from \eqref{alpha0} that, for each $\beta\in\{0,1\}^{k+1}$, there exists a Gaussian function $G_{\beta}$ such that:
\begin{equation} \label{Gb}
    \|G_{\beta}-f_{\beta}^{p_{k}}\|_{q}\leq o(\delta)\|f_{\beta}^{p_{k}}\|_{q}.
\end{equation}
Since $\|G_{\beta}^{1/p_{k}}-f_{\beta}\|_{p_{k+1}}\leq C(k)\|G_{\beta}-f_{\beta}^{p_{k}}\|_{q}^{1/p_{k}}$, by a simple calculation, \eqref{Gb} leads to $\|G_{\beta}^{1/p_{k}}-f_{\beta}\|_{p_{k+1}}\leq o(\delta)\| f_{\beta}\|_{p_{k+1}}$. \\

\noindent This discussion so far applies to $k+1\geq 3$, since $k\geq 2$. For the characterization of a near extremizing tuple $\vec{f} = (f_1, f_2, f_3, f_4)$ of the second Gowers inner product inequality, note that
\begin{equation} \label{neq4}
    \mathcal{T}_2(f_1,f_2,f_3, f_4) = \langle \tilde{f}_1\ast f_2, \tilde{f}_3\ast f_4 \rangle\leq\|\tilde{f}_1\ast f_2\|_2\|\tilde{f}_3\ast f_4\|_2\leq A(2,n)^4\prod_{i=1}^4 \|f_{i}\|_{4/3}.
\end{equation}
The first inequality is an application of the Cauchy-Schwarz inequality while the second follows from Young's convolution inequality. The appearance of the constant in \eqref{neq4} is due to the facts that $A(j,n)^{2^{j}}B(j+1,n)^{j+1} = A(j+1,n)^{2^{j+1}}$ actually holds for all $j\geq 1$ and that $A(1,n) = A(1)^{n} = 1$. If for some $\delta >0$,
\begin{equation*} 
    \mathcal{T}_2(f_1,f_2,f_3, f_4) \geq (1-\delta) A(2,n)^4\prod_{i=1}^4 \|f_{i}\|_{4/3}
\end{equation*}
then with a similar argument as above, we trade each inequality in \eqref{neq4} for a reverse inequality with a factor of $1-c\delta$ and use the said stability result for Young's convolution inequality to conclude that, for each $i\in\{1, \text{ ... }, 4\}$, there exists a Gaussian $G_{i}$ such that, $\|G_{i}-f_{i}\|_{4/3}\leq o(\delta)\|f_{i}\|_{4/3}$. As Gaussians are the only nonnegative maximizers of the Gowers-Host-Kra norm inequality; hence this concludes Theorem \ref{thm:main} for the case of nonnegative near extremizers. $\qed$

\section{Localization around a single scale} \label{sec:localization}

\subsection{Normalization}

Let $\Theta:\mathbb{R}_{> 0}\to\mathbb{R}_{>0}$ be such that $\Theta(t)\to 0$ if $t\to\infty$. Let $f:\mathbb{R}^{n}\to\mathbb{C}$ and $f\in L^{q}$, $q\in [1,\infty)$. We say that $f$ is normalized (with respect to $\Theta$) if the following two conditions occur:
\begin{align*}
    \int_{\{|f|>\rho\}} |f|^{q} \, dx &\leq \Theta(\rho)\\
    \int_{\{|f|<\rho\}} |f|^{q} \, dx &\leq \Theta(\rho^{-1}).
\end{align*}
We often say "$f$ is normalized" for short as the presence of such a function $\Theta$ will always be implied, and the selection of $\Theta$ will precede the selection of the involved near extremizers, so as to create no ambiguity. We will also allow the following relaxed version of normalization. Let $\nu > 0$ be a small number and $\Theta$ be as above. Then $f$ is said to be $\nu$-normalized, with respect to $\Theta$, if $f = g + h$, $g$ is normalized with respect to $\Theta$ and $\|h\|_{q}\leq\nu$. In our application, $q = p_{k} = \frac{2^{k}}{k+1}$. From now on, whenever we say a "growth function" we mean a function $\Theta: \mathbb{R}_{>0} \to \mathbb{R}_{>0}$ satisfying $\Theta(\rho) \to 0$ as $\rho\to\infty$. Our claim is that, there exist a growth function $\Theta$ and a $\delta_0 >0$, so that the following property occurs. If $0 < \delta \leq \delta_0$ and $f$ is a $(1-\delta)$ near extremizer of the Gowers-Host-Kra norm inequality, then there exists $\lambda >0$, such that, $\tilde{f}(x) = \lambda^{n/p_{k}}f(\lambda x)$ is normalized, with respect to $\Theta$. toward this end, we will prove an equivalent result below. \\

\noindent
Since the normalization condition concerns the absolute value of $f$ and since $|f|$ is a $(1-\delta)$ near extremizer if $f$ is, as $\|f\|_{U^{k}}\leq \||f|\|_{U^{k}}$, we consider only nonnegative functions. Let $f:\mathbb{R}^{n}\to\mathbb{R}_{\geq 0}$ be a measurable function and $f = \sum_{j\in\mathbb{Z}} 2^{j} F_{j}$ be the unique decomposition of $f$ mentioned in Section \ref{sec:basics}. Recall here that $1_{\mathcal{F}_{j}} \leq F_{j}< 2\cdot 1_{\mathcal{F}_{j}}$, and the measurable sets $\mathcal{F}_{j}$ are pairwise disjoint up to null sets. This decomposition will set up the stage for the use of Lorentz semi-norms which will be needed. In particular, we want to show:

\begin{proposition} \label{prop:norm}
There exist positive functions $\phi, \Phi$ satisfying $\lim_{t\to\infty} \phi(t) = 0$ and $\lim_{t\to 0} \Phi(t) = 0$, and positive constants $\delta_0, c_0$ such that the following holds.\\
Let $0\leq f \in L^{p_{k}}(\mathbb{R}^{n})$. Let $f = \sum_{j\in\mathbb{Z}} 2^{j} F_{j}$ with the $F_{j}$ related to $f$ as indicated above. Suppose that $\|f\|_{p_{k}}=1$, and for some $0<\delta\leq\delta_0$, we have
\begin{equation*}
    \| f \|_{U^{k}} \geq (1-\delta)A(k,n)\|f\|_{p_{k}} = (1-\delta)A(k,n).
\end{equation*} 
Then there exists $l \in \mathbb{Z}$ such that:
\begin{align}
    \label{lwbd} 2^{l}(\mathcal{L}(\mathcal{F}_{l}))^{1/p_{k}} &\geq c_0\|f\|_{p_{k}} = c_0 \\
    \label{upbd} \sum_{|j-l| \geq m} 2^{jp} \mathcal{L}(\mathcal{F}_{j}) &\leq (\phi(m) + \Phi(\delta))\cdot\|f\|_{p_{k}}^{p_{k}} = \phi(m)+\Phi(\delta).
\end{align}
\end{proposition}

\begin{remark} 
Statement \eqref{lwbd} ensures that there exists a scale $2^{l}$ so that the contribution of the quantity $2^{l}(\mathcal{L}(\mathcal{F}_{l}))^{1/p_{k}}$ to the total norm $\| f \|_{p_{k}}$ is non-negligible. Choose $\lambda$ to be an integral power of $2$ to dilate the set $\mathcal{F}_{l}$, so that its dilated version has measure as close to one as possible. Statement \eqref{upbd} is equivalent to saying, this dilated version of $f$ is then $\nu$-normalized with respect to an appropriate growth function $\Theta$ and a positive $\nu=\nu(\delta)$ that tends to zero as $\delta\to 0$. 
\end{remark}

\noindent
Let $1\leq q, \tilde{q} < \infty$ and $f = \sum_{j\in\mathbb{Z}} 2^{j} F_{j}$ be as above. Recall from Section \ref{sec:basics} that $\| f \|_{(q, \tilde{q})} \asymp (\sum_{j\in\mathbb{Z}} 2^{j\tilde{q}}(\mathcal{L}(\mathcal{F}_{j}))^{\tilde{q}/q})^{1/\tilde{q}}$ if $\tilde{q}<\infty$, and $\|f\|_{(q,\infty)}\asymp\sup_{j:\mathcal{F}_{j}\not=\emptyset} 2^{j}(\mathcal{L}(\mathcal{F}_{j}))^{1/q}$.

\subsection{A digression}

Let $\vec{f}$ denote $(f_{\alpha} : \alpha\in\{0,1\}^{k})$. We claim that, there exists $C(k,n) <\infty$ such that:
\begin{equation} \label{multlin}
    |\mathcal{T}_{k} (\vec{f})| = \bigg | \int_{\mathbb{R}^{(k+1)n}} \prod_{\alpha\in\{0,1\}^{k}} f_{\alpha}(x+\alpha\cdot\vec{h}) \, dx d\vec{h} \, \bigg | \leq C(k,n) \prod_{\alpha\in\{0,1\}^{k}} \| f_{\alpha} \|_{(p_{k},2^{k})}. 
\end{equation}
This is a consequence of Lorentz space interpolation. Indeed:\\
Let $F(x) = \int_{\mathbb{R}^{nk}} \prod_{\alpha\in\{0,1\}^{k};\alpha\not=\vec{0}} f_{\alpha}(x+\alpha\cdot\vec{h}) \, d\vec{h}$ and $(q')^{-1}+q^{-1}=1$. By H\"older's inequality for $1\leq r,\tilde{r} \leq\infty$ for the Lorentz semi-norms \cite{stein1971introduction},
\begin{equation} \label{Tkdominated}
    |\mathcal{T}_{k} (\vec{f})| \leq C\| f_{\vec{0}}\|_{(r,\tilde{r})} \| F \|_{(r',\tilde{r}')}. 
\end{equation}
Let $Q_{r}(\vec{P})$ denote a closed cube centered at $\vec{P} = (p_{k},\text{ ... }, p_{k})\in\mathbb{R}^{2^{k}}$ with radius $r$ and $\Omega = \{ \vec{q} = (q_{\alpha}: \alpha\in\{0,1\}^{k}) \in Q_{r}(\vec{P}): \sum_{\alpha\in\{0,1\}^{k}} (q_{\alpha})^{-1} = k+1\}$. Then if $r=r(k)$ is sufficiently small, there exists a constant $C(k,n)<\infty$ such that for every $\vec{q}\in\Omega$, 
\begin{equation} \label{multlinver}
    |\mathcal{T}_{k}(\vec{f})| \leq C(k,n) \prod_{\alpha\in\{0,1\}^{k}} \| f_{\alpha}\|_{q_{\alpha}}.
\end{equation}
A proof of \eqref{multlinver} will be given at the end of this section. Assuming \eqref{multlinver}, then if $\vec{q}\in\Omega$,
\begin{equation} \label{remdominated1}
    \| F \|_{q_{\vec{0}}'} \leq C(k,n) \prod_{\alpha\in\{0,1\}^{k};\alpha\not=\vec{0}} \| f_{\alpha} \|_{q_{\alpha}}. 
\end{equation}
If we vary one $q_{\alpha}, \alpha\not= \vec{0}$, while keeping other values $q_{\beta}, \beta\not=\vec{0},\alpha$, fixed and apply interpolation for the Lorentz norms, we obtain from \eqref{remdominated1}:
\begin{equation*} 
    \| F \|_{(r_{\alpha}',\infty)} \leq C(k,n) \|f_{\alpha} \|_{(p_{k},\infty)}\prod_{\beta\in\{0,1\}^{k};\beta\not=\vec{0},\alpha} \| f_{\beta} \|_{q_{\beta}}
\end{equation*}
whenever $p_{k}^{-1}+r_{\alpha}^{-1} + \sum q_{\beta}^{-1} = k+1$ and $r_{\alpha}, q_{\beta} \in [p_{k}-2^{-k},p_{k}+2^{-k}]$. We continue applying interpolation on other exponents $q_{\beta}$, one by one. Suppose that $q_{\gamma}$ is the last exponent to be interpolated, then:
\begin{equation*} 
    \| F\|_{(r_{\gamma}',\infty)} \leq C(k,n) \|f_{\gamma}\|_{q_{\gamma}} \prod_{\beta\in\{0,1\}^{k};\beta\not=\vec{0},\gamma}\|f_{\beta}\|_{(p_{k},\infty)}
\end{equation*}
whenever $r_{\gamma}^{-1}+q_{\gamma}^{-1}+(2^{k}-2)p_{k}^{-1} = k+1$ and $r_{\gamma}, q_{\gamma} \in [p_{k}-2^{-k},p_{k}+2^{-k}]$. Applying interpolation one last time, we have:
\begin{align*}
    \|F\|_{(p_{k}',1)} &\leq C(k,n) \|f_{\gamma}\|_{(p_{k},1)}\prod_{\alpha\in\{0,1\}^{k};\alpha\not=\vec{0},\gamma} \| f_{\alpha}\|_{(p_{k},\infty)}\\
    \|F\|_{(p_{k}',\infty)} &\leq C(k,n) \|f_{\gamma}\|_{(p_{k},\infty)} \prod_{\alpha\in\{0,1\}^{k};\alpha\not=\vec{0},\gamma} \| f_{\alpha}\|_{(p_{k},\infty)}.
\end{align*}
Since $\gamma$ can be any of the values $\alpha\in\{0,1\}^{k}$, these calculations yield, 
\begin{align}
    \label{remdom1} \| F \|_{(p_{k}',1)} &\leq C(k,n) \prod_{\alpha\in\{0,1\}^{k};\alpha\not=\vec{0}} \|f_{\alpha}\|_{(p_{k},r_{\alpha})}\\
    \label{remdom2} \| F\|_{(p'_{k},\infty)} &\leq C(k,n) \prod_{\alpha\in\{0,1\}^{k};\alpha\not=\vec{0}} \|f_{\alpha}\|_{(p_{k},\infty)}
\end{align}
whenever $\sum r_{\alpha}^{-1} = 1$ and $r_{\alpha} \in \{1,\infty\}$. Then \eqref{Tkdominated}, \eqref{remdom1} and \eqref{remdom2} yield:
\begin{equation} \label{remdominated2}
    |\mathcal{T}_{k}(\vec{f})| \leq C\|f_{\vec{0}}\|_{(p_{k},r_{\vec{0}})}\| F\|_{(p_{k}',r_{\vec{0}}')} \leq C(k,n) \prod_{\alpha\in\{0,1\}^{k}} \|f_{\alpha}\|_{(p_{k},r_{\alpha})} 
\end{equation}
whenever $\sum r_{\alpha}^{-1} = 1$ and $r_{\alpha}\in\{1,\infty\}$. To complete the proof of \eqref{multlin}, we need the following short lemma:

\begin{lemma} \label{lem:Tklast} \cite{personalcommunication} Let $\mathcal{T}$ be a scalar-valued $m$-linear form on $\mathbb{R}^{n}$ and $q\in [1,\infty)$. Suppose that 
\begin{equation} \label{Tkri}
    |\mathcal{T}(\vec{f})|\leq C(m) \prod_{1\leq i\leq m} \|f_{i}\|_{(q,r_{i})}
\end{equation}
whenever $r_{i}\in\{1,\infty\}$ and $\sum r_{i}^{-1} = 1$. Suppose furthermore that $\vec{q} = (q_{i})_{i=1}^{m}$ with $\sum q_{i}^{-1} = 1$ and $q_{i} \in (1,\infty)$. Then 
\begin{equation} \label{Tkqi}
    |\mathcal{T}(\vec{f})|\leq C(m,\vec{q})\prod_{1\leq i \leq m} \|f_{i}\|_{(q,q_{i})}.
\end{equation}   
\end{lemma}

\begin{proof} We use the described decomposition $f_{i} = \sum_{j\in\mathbb{Z}} 2^{j}F_{i,j}$, with $1_{\mathcal{F}_{i,j}}\leq F_{i,j} < 2\cdot 1_{\mathcal{F}_{i,j}}$ and the measurable sets $\mathcal{F}_{i,j}$ pairwise disjoint. Let $S_{i,\lambda} = \{j: 2^{\lambda} \leq 2^{j}(\mathcal{L}(\mathcal{F}_{i,j}))^{1/p_{k}} < 2^{\lambda+1} \}$ for each $\lambda\in\mathbb{Z}$. Define $f_{i,\lambda} = \sum_{j\in S_{i,\lambda}} 2^{j}F_{i,j}$. Note that $f_{i} = \sum_{\lambda\in\mathbb{Z}} 2^{\lambda}f_{i,\lambda}$ and for $s\in[1,\infty)$:
\begin{equation} \label{Lorentzs}
    \| f_{i} \|_{(q,s)} \asymp (\sum_{\lambda\in\mathbb{Z}} 2^{\lambda s} crd(S_{i,\lambda}))^{1/s}.
\end{equation}
If $s=\infty$:
\begin{equation} \label{Lorentzinf}
    \|f_{i}\|_{(q,\infty)} \asymp \lim_{t\to\infty}  (\sum_{\lambda\in\mathbb{Z}} 2^{\lambda t} crd(S_{i,\lambda}))^{1/t} = \sup_{\lambda: S_{i,\lambda}\not=\emptyset} 2^{\lambda}.
\end{equation}
Let $\vec{\lambda} = (\lambda_{i})_{i=1}^{m} \in\mathbb{Z}^{m}$. Then from \eqref{Tkri}, \eqref{Lorentzs} and \eqref{Lorentzinf}:
\begin{equation} \label{TLorentz}
    |\mathcal{T}(\vec{f})| \leq C(m) \sum_{\vec{\lambda}} 2^{\sum_{l} \lambda_{l}} \min_{1\leq i\leq m} (crd(S_{i,\lambda_{i}})) \leq C(m) \sum_{\vec{\lambda}} \theta((crd(S_{i,\lambda_{i}}))_{i=1}^{m}) \prod_{1\leq i\leq m} 2^{\lambda_{i}}crd(S_{i,\lambda_{i}})^{1/q_{i}} 
\end{equation}
with $\sum q_{i}^{-1} = 1$ and $\theta((t_{i})_{i=1}^{m}) = (\min_{i} t_{i}) \cdot \prod_{1\leq i\leq m} t_{i}^{-1/q_{i}}$. Consider all the vectors $\vec{\lambda}$ for which:
\begin{align}
    \label{cardmin} crd(S_{1,\lambda_1}) &= \min_{i} crd(S_{i,\lambda_{i}})\\
    \label{cardmax} crd(S_{m,\lambda_{m}}) &= \max_{i} crd(S_{i,\lambda_{i}}). 
\end{align}
The following analysis applies the same way for other vectors $\vec{\lambda}$ with only minor changes in notation. For each $i\in\{1,\text{ ... }, m\}$, define a nonnegative integer $n_{i}$ so that 
\begin{equation} \label{cardcomp}
    2^{-n_{i}-1}crd(S_{1,\lambda_1}) \leq crd(S_{i,\lambda_{i}}) \leq 2^{-n_{i}}crd(S_{1,\lambda_1}).
\end{equation}
Fix a tuple $\vec{n}=(n_{i})_{i=1}^{m} \in\mathbb{N}^{m}$; we sum the right hand side of \eqref{TLorentz} over the $\vec{\lambda}$ for which all three \eqref{cardmin}, \eqref{cardmax} and \eqref{cardcomp} are satisfied with this fixed $\vec{n}$. Define $I_1= 1$ and $I_{i}$, $i\not= 0$, is the largest index $l$ such that $2^{-n_{i}-1}crd(S_{1,\lambda_1}) \leq crd(S_{i,l}) \leq 2^{-n_{i}}crd(S_{1,\lambda_1})$. With these choices of $I_{i}$, it's clear that $I_1 \mapsto I_{i}$ are bounded-to-one functions, uniformly in $\vec{n}$. Then for this fixed $\vec{n}$, the sum on the right hand side of \eqref{TLorentz} will be majorized by:
\begin{equation} \label{tailcon1}
    C(m) \sum_{I_{m}} 2^{-n_{m}} \prod_{1\leq i < m} 2^{n_{i}/q_{i}} \prod_{1\leq t\leq m} 2^{I_{t}} crd(S_{t, I_{t}})^{1/q_{t}}\leq C(m) 2^{-r\sum u_{i}} \sum_{I_{m}} \prod_{1\leq t\leq m} 2^{I_{t}} crd(S_{t,I_{t}})^{1/q_{t}}
\end{equation}
with $u_{m} = n_{m}$ and $u_{i} = n_{m} -n_{i} \geq 0$ and $0< r = \min_{i} 1/q_{i}$. Once again, with these choices of $I_{i}$, the constant $C(m)$ in \eqref{tailcon1} doesn't depend on $\vec{n}$. Applying H\"older's inequality to the right hand side sum in \eqref{tailcon1} to have it further majorized by:
\begin{equation} \label{tailcon2}
    C(m) 2^{-r\sum u_{i}} \prod_{1\leq i\leq m} (\sum_{I_{m}} 2^{I_{i} \cdot q_{i}} crd(S_{i, I_{i}}))^{1/q_{i}} \leq C(m) 2^{-r\sum u_{i}} \prod_{1\leq i\leq m} \| f_{i}\|_{(q,q_{i})}.
\end{equation}
Summing the sum in the right hand side of \eqref{TLorentz} over $\vec{n}\in\mathbb{N}^{m}$ and utilizing the bound in \eqref{tailcon2} and the convergence of geometric series, we conclude that \eqref{Tkri} indeed implies \eqref{Tkqi}.    
\end{proof}

\noindent
Applying this Lemma \ref{lem:Tklast} to \eqref{remdominated2}, with $q=p_{k}$ and $m = 2^{k}$ with indices $\alpha\in\{0,1\}^{k}$, we obtain \eqref{multlin} with the desired exponents $q_{i} = 2^{k}$.

\begin{remark}
The interpolation result of this discussion is a very restricted result that is suited to the task at hand. A more general interpolation result can be found in \cite{janson2006interpolation}.    
\end{remark}

\begin{proof}[Proof of \eqref{multlinver}]
\cite{personalcommunication} We prove the following by induction. For every $k\geq 2$, there exists $r_{k}$ such that, if $\vec{p} = (p_{\alpha}:\alpha\in\{0,1\}^{k})$ is such that $\sum_{\alpha} p_{\alpha}^{-1} = k+1$ and $|p_{\alpha}-2^{k}/(k+1)|\leq r_{k}$, for every $\alpha\in\{0,1\}^{k}$, then
\begin{equation*} 
    \mathcal{T}_{k}(\vec{f})\leq \prod_{\alpha\in\{0,1\}^{k}}\|f_{\alpha}\|_{p_{\alpha}}
\end{equation*}
for every $\vec{f} = (f_{\alpha}:\alpha\in\{0,1\}^{k})$ with each $f_{\alpha}$ nonnegative and measurable.\\
If $k=2$ then the conclusion is a consequence of Young's and H\"older's inequalities:
\begin{equation*} 
    \mathcal{T}_{k}(\vec{f}) = \langle f_1\ast\tilde{f}_2,f_3\ast\tilde{f}_4\rangle_\mathbb{R}\leq \|f_1\ast\tilde{f}_2\|_{p_{1,2}}\|f_3\ast\tilde{f}_4\|_{p_{3,4}}\leq\prod_{i=1}^4 \|f_{i}\|_{p_{i}}.
\end{equation*}
Here, $\tilde{f}(x) = f(-x)$. The condition on $p_{i}, i\in\{1,\text{ ... },4\}$ is simply, $p_{1,2}^{-1}+p_{3,4}^{-1}=1$, or equivalently, $\sum_{i=1}^4 p_{i}^{-1} = 2+1$. \\
Assume the conclusion holds for the case $k-1$. Let $\vec{f} = (f_{\alpha}\geq 0:\alpha\in\{0,1\}^{k}\}$. For each $\alpha\in\{0,1\}^{k}$, there is a unique $\beta\in\{0,1\}^{k-1}$ such that $\alpha = (\beta,0)$ or $\alpha = (\beta,1)$. Then for each $\beta\in\{0,1\}^{k}$, let $g_{\beta}^{t}(x) = f_{(\beta,0)}(x)f_{(\beta,1)}(x+t)$. Write $\vec{g}^{t} = (g^{t}_{\beta}:\beta\in\{0,1\}^{k-1})$, then
\begin{align*} 
    \mathcal{T}_{k}(\vec{f}) = \int_{\mathbb{R}^{n}} \mathcal{T}_{k-1}(\vec{g}^{t})\, dt &\leq \int_{\mathbb{R}^{n}} \prod_{\beta\in\{0,1\}^{k-1}} \|g_{\beta}^{t}\|_{p_{\beta}}\, dt \\ 
    &\leq \prod_{\beta\in\{0,1\}^{k-1}} \left (\int_{\mathbb{R}^{n}} \|g^{t}_{\beta}\|_{p_{\beta}}^{2^{k-1}}\, dt \right )^{1/2^{k-1}}\leq\prod_{\beta\in\{0,1\}^{k-1}} \|f_{(\beta,0)}\|_{p_{(\beta,0)}}\|f_{(\beta,1)}\|_{p_{(\beta,1)}}.
\end{align*}
The first inequality follows from the induction hypothesis for some $r_{k-1} >0$ such that, $\sum_{\beta\in\{0,1\}^{k-1}} p_{\beta}^{-1} = k$ and $|p_{\beta} - 2^{k-1}/k|\leq r_{k-1}$, for every $\beta\in\{0,1\}^{k-1}$. The second inequality follows from Minkowski's inequality, and the last from Young's inequality, with the following relations between the exponents: $1/2^{k-1} + 1/p_{\beta} =1/p_{(\beta,0)} +1/p_{(\beta,1)}$. Now it's easily observed that if $p_{(\beta,0)} = 2^{k}/(k+1) = p_{(\beta,1)}$ then $p_{\beta} = 2^{k-1}/k$, if $\sum_{\beta\in\{0,1\}^{k-1}} p_{\beta}^{-1} = k$ then $\sum_{\beta\in\{0,1\}^{k}} p_{(\beta,0)}^{-1} + p_{(\beta,1)}^{-1} = k+1$, and for each $\beta\in\{0,1\}^{k-1}$, $p_{\beta}$ is a continuous function of $p_{(\beta,0)}$ and $p_{(\beta,1)}$. Hence if for each $\beta\in\{0,1\}^{k-1}$, $p_{(\beta,0)}$ and $p_{(\beta,1)}$ are taken sufficiently close to $2^{k}/(k+1)$, so that $p_{\beta}$ is within an $r_{k-1}$ distance to $2^{k-1}/k$, then the conclusion for the case $k$ follows from the induction hypothesis for the case $k-1$. The proof of \eqref{multlinver} is now complete.    
\end{proof}

\subsection{Proof of Proposition \ref{prop:norm}}

Let $\eta = \eta(\delta) > 0$ be a small number to be chosen later. Define $S = S(\eta) = \{ j \in \mathbb{Z}: 2^{j} (\mathcal{L}(\mathcal{F}_{j}))^{1/p_{k}} > \eta \}$. Statement \eqref{lwbd} will follow if this set $S$ is nonempty, which is indeed true if $\eta$ is sufficiently small. To see this, consider $f_{S} = \sum_{j\in S} 2^{j} F_{j}$ and $f_{S^{c}} = \sum_{j\not\in S} 2^{j} F_{j}$. Then:
\begin{align*} 
    \| f_{S^{c}} \|_{(p_{k}, 2^{k})}^{2^{k}} = \| \sum_{j\not\in S} 2^{j} F_{j} \|_{(p_{k},2^{k})}^{2^{k}} &\asymp \sum_{j\not\in S} 2^{j2^{k}} (\mathcal{L}(\mathcal{F}_{j}))^{2^{k}/p_{k}}\\ 
    &\leq \max_{j\not\in S} (2^{j} (\mathcal{L}(\mathcal{F}_{j}))^{1/p_{k}})^{2^{k}-p_{k}} \sum_{j\not\in S} (2^{j} (\mathcal{L}(\mathcal{F}_{j}))^{1/p_{k}})^{p_{k}} \leq \eta^{2^{k}-p_{k}} \| f_{S^{c}} \|_{p_{k}}^{p_{k}} 
\end{align*} 
or 
\begin{equation} \label{compdom}
    \| f_{S^{c}} \|_{(p_{k},2^{k})} \leq C\eta^{1-p_{k}/2^{k}} \| f_{S^{c}} \|_{p_{k}}^{p_{k}/2^{k}}. 
\end{equation}

\noindent As above, for a set $A\subset\mathbb{Z}$, we denote $f_{A} = \sum_{j\in A} 2^{j} F_{j}$ and define a $2^{k}$ tuple-valued function on $\mathbb{R}^{n}$, $\vec{f}_{A} = (f_{A}, f, \text{ ... }, f)$. With this notation:
\begin{multline} \label{compdomnorm}
    |\mathcal{T}_{k} (\vec{f}_{S^{c}})| = \big | \int_{\mathbb{R}^{(k+1)n}} f_{S^{c}}(x) \prod_{\alpha\in\{0,1\}^{k};\alpha\not= 0} f(x+\alpha\cdot\vec{h}) \, dx d\vec{h} \, \big | \\ \leq C(k,n) \| f_{S^{c}} \|_{(p_{k},2^{k})} \| f\|_{p_{k}}^{2^{k}-1} \leq C(k,n) \eta^{1-p_{k}/2^{k}} \|f_{S^{c}}\|_{p_{k}}^{p_{k}/2^{k}} \|f\|_{p_{k}}^{2^{k}-1}.
\end{multline}
The first inequality follows from \eqref{multlin} and the fact that $\| f\|_{(p_{k},2^{k})}\leq \|f\|_{p_{k}}$ and the last inequality from \eqref{compdom}. Recall that $\| f\|_{p_{k}} = 1$. Hence $S=\emptyset$ in the context of \eqref{compdomnorm} would imply that:
\begin{equation} \label{Uknormdom}
    \| f\|_{U^{k}}^{2^{k}} \leq C(k,n)\eta^{1-p_{k}/2^{k}} \| f\|_{p_{k}}^{2^{k}} = C(k,n)\eta^{1-p_{k}/2^{k}}.
\end{equation}
Since $1-p_{k}/2^{k} > 0$, the right hand side of \eqref{Uknormdom} is small if $\eta$ is small. From the near extremizing hypothesis, $\| f\|_{U^{k}} \geq A(k,n) (1-\delta)$. Hence, if $\eta$ is sufficiently small then $S = S(\eta)$ must be nonempty; otherwise \eqref{Uknormdom} implies $\|f\|_{U^{k}} = 0$, a contradiction. Note that the argument also shows, for given $\eta$, if $\delta$ is sufficiently small then $S = S(\eta) \not=\emptyset$. In addition, from the definition of $S$ and Chebyshev's inequality, 
\begin{equation} \label{cardom}
    crd(S) \leq C\eta^{-p_{k}}\|f\|_{p_{k}}^{p_{k}} = C\eta^{-p_{k}}.
\end{equation}

\noindent
We now prove statement \eqref{upbd}. Let $\eta > 0$ be a small number, and consequently let $\delta > 0$ be small enough if needed, so that $S=S(\eta) \not=\emptyset$. Let $M = \max_{i,i' \in S} |i-i'|$. In order to prove \eqref{upbd}, it suffices to show that $M$ is bounded by a finite upper bound that does not depend on a particular $(1-\delta)$ near extremizer $f$. Assume that $S$ has more than one element. Note that since $S$ is a set of integers, $M \geq crd(S)$. Let $N > 0$ be a large integer to be chosen below, which will depend only on $n, k$ and $\eta$ - and consequently on $\delta$. Suppose that with this choice of $N$ we have $M \leq 10 N crd(S)$, then we obtain the desired bound, since from \eqref{cardom}:
\begin{equation*} 
    M\leq CNcrd(S) \leq C(k,n,\delta)\eta(\delta)^{-p_{k}}.
\end{equation*}
If $M>10Ncrd(S)$, then there must exist integers $J > I$ such that: $S\cap (-\infty,I] \not= \emptyset$, $S\cap [J,\infty) \not=\emptyset$, $S\cap (I,J) = \emptyset$ and $J-I \geq M/(N crd(S))$. Indeed, choose an integer $L$ such that $M/(2N crd(S)) \leq L \leq M/(N crd(S))$. Consider the intervals $I_{j} = [jL, (j+1)L) \subset (\inf(S), \sup(S))$ and let $K$ be the set of all these indices $j$. Then $crd(K) \asymp M/L \asymp N crd(S)$. Hence $I_{j}\cap S =\emptyset$ for at least one of these intervals. Since 
\begin{equation*} 
    \| f_{S^{c}}\|_{(p_{k},2^{k})}^{2^{k}} = \sum_{l\not\in S} \| 2^{l}F_{l}\|_{(p_{k},2^{k})}^{2^{k}} \geq \sum_{j \in K: I_{j}\cap S = \emptyset} \| \sum_{l\in I_{j}} 2^{l}F_{l}\|_{(p_{k},2^{k})}^{2^{k}}
\end{equation*}
and $crd(K)-crd(S) \geq c(N-1)crd(S) \geq cN-1$, there exists one such interval $I_{j}$ such that, 
\begin{equation*} 
    \| f_{S^{c}}\|_{(p_{k},2^{k})}^{2^{k}} \geq CN^{-1}\| \sum_{l\in I_{j}} 2^{l}F_{l}\|_{(p_{k},2^{k})}^{2^{k}}.
\end{equation*}
Take such an interval, then $I = jL$ and $J = (j+1)L$ are our desired integers. Moreover, with this pair of $(I,J)$:
\begin{equation} \label{IJ}
    \| \sum_{I < l < J} 2^{l} F_{l}\|_{(p_{k},2^{k})}^{2^{k}} \leq CN^{-1} \| f_{S^{c}} \|_{(p_{k},2^{k})}^{2^{k}} \leq CN^{-1}\eta^{2^{k}-p_{k}}.
\end{equation}
The second inequality follows from \eqref{compdom}. Define $f_{u} = \sum_{i\geq J} 2^{i}F_{i}$, $f_{d} = \sum_{i\leq I} 2^{i}F_{i}$ and $f_{b} = f -f_{u}-f_{d}$. As \eqref{IJ} states, $\| f_{b} \|_{p_{k}} \leq CN^{-1}\eta^{c}$ for some $c>0$. Now in terms of these $f_{u}, f_{d}, f_{b}$:
\begin{equation*} 
    \|f\|_{U^{k}}^{2^{k}} = | \mathcal{T}_{k}(f_{u}+f_{d}+f_{b}, \text{ ... }, f_{u}+f_{d}+f_{b})|.
\end{equation*}
Expanding the above sum, we will have $C(k)$ terms in which at least one of the components is $f_{b}$, $C(k)$ mixed terms in which the components consist of only $f_{u}$ and $f_{d}$ but no $f_{b}$, and two pure terms in which the components are either $f_{u}$ or $f_{d}$, which we denote $\mathcal{T}_{k}(\vec{f}_{u})$ and $\mathcal{T}_{k}(\vec{f}_{d})$, respectively. We claim that the contribution of all the mixed terms is majorized by $C(k,n) 2^{-\eta^{p_{k}}M/N}$. We will first need the following fact about the $2^{k}$ linear form $\mathcal{T}_{k}$ on indicators of sets. Recall that $\Omega = \{ \vec{q} = (q_{\alpha}: \alpha\in\{0,1\}^{k}) \in Q_{r}(\vec{P}): \sum_{\alpha\in\{0,1\}^{k}} (q_{\alpha})^{-1} = k+1\}$ and $Q_{r}(\vec{P})$ denotes a closed cube of sufficiently small size $r$ centered at $\vec{P} = (p_{k}, \text{ ... }, p_{k})\in\mathbb{R}^{2^{k}}$. Then:

\begin{lemma} \label{lem:Tksets}
Suppose $\vec{p} = (p_{\alpha}: \alpha\in\{0,1\}^{k})\in\Omega$. Then there exist $\tau > 0$ and $C$ such that, if $E_{\alpha} \subset \mathbb{R}^{n}$ are sets with finite measures, $\alpha\in\{0,1\}^{k}$, then:
\begin{equation} \label{lemTksets}
    \mathcal{T}_{k}(1_{E_{\alpha}}:\alpha\in\{0,1\}^{k}) \leq C(k,n) \rho (\mathcal{L}(E_{\alpha}): \alpha\in\{0,1\}^{k})^{\tau}\prod_{\alpha\in\{0,1\}^{k}} (\mathcal{L}(E_{\alpha}))^{1/p_{\alpha}}
\end{equation}
with $\rho(c_{\alpha}:\alpha\in\{0,1\}^{k}) = \min_{\beta\not=\gamma\in\{0,1\}^{k}} \frac{c_{\beta}}{c_{\gamma}}$.
\end{lemma}

\begin{proof}
Select $\vec{q} = (q_{\alpha}: \alpha\in\{0,1\}^{k})\in\Omega$ that differs from $\vec{p} = (p_{\alpha}: \alpha\in\{0,1\}^{k})$ only in two components, for convenience, say $\alpha = \vec{0}$ and $\alpha = (1,\text{ ... }, 1) = \vec{1}$, so that $q_{\vec{0}}^{-1} = p_{\vec{0}}^{-1} + \tau$ and $q_{\vec{1}}^{-1} = p_{\vec{1}}^{-1} - \tau$ for some $\tau > 0$ small enough. Then from \eqref{multlinver}, for this pair $(\vec{0},\vec{1})$ we obtain:
\begin{equation*} 
    \mathcal{T}_{k}(1_{E_{\alpha}}:\alpha\in\{0,1\}^{k}) \leq C(k,n) \prod_{\alpha\in\{0,1\}^{k}} |E_{\alpha}|^{1/q_{\alpha}} \leq C(k,n) \left ( \frac{\mathcal{L}(E_{\vec{0}})}{\mathcal{L}(E_{\vec{1}})} \right )^{\tau} \prod_{\alpha\in\{0,1\}^{k}} (\mathcal{L}(E_{\alpha}))^{1/p_{\alpha}}.
\end{equation*}
Repeating this process to the other pairs $(\beta, \gamma)$, $\beta,\gamma\in\{0,1\}^{k}$, we obtain the conclusion \eqref{lemTksets}.
\end{proof}

\begin{lemma} \label{lem:mixed}
Using the terminology above, then the contribution of all the quantities $\mathcal{T}_{k}$ with mixed-term inputs is majorized by $C(k,n) 2^{-\eta^{p_{k}}M/N}$ in absolute value.    
\end{lemma}

\begin{remark}
This lemma confirms the intuition that if the components are not "compatible" then their multilinear product doesn't contribute much.    
\end{remark}

\begin{proof}[Proof of Lemma \ref{lem:mixed}] Take a typical mixed term $\mathcal{T}_{k} (f_{u}, \text{ ... }, f_{d})$ in the expansion of $\mathcal{T}_{k}(f_{u}+f_{d}+f_{b}, \text{ ... }, f_{u}+f_{d}+f_{b})$, whose at least two components are assumed to be different. The notation used here is only for convenience; the argument will be independent of which two components are different. Let $\mathbb{Z}^{2^{k}} \supset R = \{ \vec{i} = (i_{\alpha}: \alpha\in\{0,1\}^{k}) : i_{\vec{0}} \geq J, i_{\vec{1}} \leq I\}$. Let $\epsilon > 0$ be a small number and let $R_{\epsilon} = \{ \vec{i} \in R: 2^{i_{\alpha}}\mathcal{L}(\mathcal{F}_{i_{\alpha}}))^{1/p_{k}} \geq \epsilon, \forall\alpha\in\{0,1\}^{k}\}$. We take $\epsilon$ sufficiently small so that $R_{\epsilon}\not=\emptyset$. By Chebyshev's inequality, 
\begin{equation} \label{cardRe}
    crd(R_{\epsilon}) \leq C\epsilon^{-2^{k}p_{k}}.
\end{equation}
By the definition of $R_{\epsilon}$ and the Gowers product inequality:
\begin{equation} \label{TkRe}
    \sum_{\vec{i} \in R_{\epsilon}^{c}} 2^{\sum_{\alpha\in\{0,1\}^{k}} i_{\alpha}} \mathcal{T}_{k}(1_{\mathcal{F}_{i_{\alpha}}}:\alpha\in\{0,1\}^{k}) \leq C\epsilon.
\end{equation}
Next we need to find an upper bound for the corresponding sum over $\vec{i}\in R_{\epsilon}$. Note that if $\vec{i} = (i_{\alpha}:\alpha\in\{0,1\}^{k}) \in R_{\epsilon}$, then $2^{i_{\alpha}}(\mathcal{L}(\mathcal{F}_{i_{\alpha}}))^{1/p_{k}} \leq C$, since $\|f\|_{p_{k}} = 1$, and that $2^{i_{\vec{1}}}(\mathcal{L}(\mathcal{F}_{i_{\vec{1}}}))^{1/p_{k}} \geq \epsilon$, which then implies
\begin{equation} \label{extremeL1}
    \mathcal{L}(\mathcal{F}_{i_{\vec{1}}}) \geq \epsilon^{p_{k}} 2^{-i_{\vec{1}}p_{k}} \geq \epsilon^{p_{k}}2^{-Ip_{k}}.
\end{equation} 
Furthermore, since $i_{\vec{0}} \geq J \geq I + M/(Ncrd(S)) \geq I + c\eta^{p_{k}}M/N$,
\begin{equation} \label{extremeL0}
    \mathcal{L}(\mathcal{F}_{i_{\vec{0}}}) \leq C2^{-i_{\vec{0}}p_{k}} \leq C2^{-Ip_{k}-c\eta^{p_{k}}M/N}.
\end{equation}
\eqref{extremeL1} and \eqref{extremeL0} conclude:
\begin{equation} \label{ratL0L1}
    \frac{\mathcal{L}(\mathcal{F}_{i_{\vec{0}}})}{\mathcal{L}(\mathcal{F}_{i_{\vec{1}}})} \leq C\epsilon^{-p_{k}}2^{-c\eta^{p_{k}}M/N}.
\end{equation}
Applying Lemma \ref{lem:Tksets} to the sets $\mathcal{L}(\mathcal{F}_{i_{\alpha}})$ and utilizing \eqref{cardRe} and \eqref{ratL0L1}, we obtain a bound on the sum over $\vec{i}\in R_{\epsilon}$:
\begin{align} 
    \nonumber \sum_{\vec{i} \in R_{\epsilon}} 2^{\sum_{\alpha\in\{0,1\}^{k}} i_{\alpha}} \mathcal{T}_{k}((1_{\mathcal{F}_{i_{\alpha}}}:\alpha\in\{0,1\}^{k})) &\leq C(k,n)\epsilon^{-p_{k}}2^{-c\eta^{p_{k}}M/N} crd(R_{\epsilon}) \\
    \label{Tkalpha} &\leq C(k,n)\epsilon^{-C(k)}2^{-c\eta^{p_{k}}M/N}. 
\end{align}
\eqref{TkRe} and \eqref{Tkalpha} give:
\begin{equation} \label{Tkmixed}
    |\mathcal{T}_{k} (f_{u}, \text{ ... }, f_{d})| \leq C\epsilon + C(k,n)\epsilon^{-C(k)}2^{-c\eta^{p_{k}}M/N}.
\end{equation}
Choose $\epsilon = e^{-\tau \eta^{p_{k}} M/N}$ for a sufficiently small $\tau = \tau(k) > 0$; we have from \eqref{Tkmixed}:
\begin{equation*} 
    |\mathcal{T}_{k} (f_{u}, \text{ ... }, f_{d})| \leq C(k,n)2^{-c\eta^{p_{k}}M/N}
\end{equation*}
which is the conclusion of Lemma \ref{lem:mixed}.
\end{proof}

\noindent
Recall that \eqref{IJ} implies $\| f_{b} \|_{p_{k}} \leq CN^{-1}\eta^{c}$. Then \eqref{Tkmixed}, the definition of the Gowers-Host-Kra norms and the Gowers-Host-Kra norm inequality give
\begin{multline} \label{Ukdu}
    \|f\|_{U^{k}}^{2^{k}} \leq |\mathcal{T}_{k}(\vec{f}_{u})| + |\mathcal{T}_{k}(\vec{f}_{d})| + C(k,n) 2^{-\eta^{p_{k}}M/N} + CN^{-1}\eta^{c} \\
    \leq A(k,n)^{2^{k}}( \|f_{u}\|_{p_{k}}^{2^{k}} + \|f_{d}\|_{p_{k}}^{2^{k}} )+ C(k,n) 2^{-\eta^{p_{k}}M/N} + CN^{-1}\eta^{c}.
\end{multline}
From the definition of $f_{u}$ and $f_{d}$, we have $\min(\|f_{u}\|_{p_{k}}, \|f_{d}\|_{p_{k}}) \geq \eta$. Since $(\|f_{u}\|_{p_{k}}^{p_{k}} + \|f_{d}\|_{p_{k}}^{p_{k}})^{1/p_{k}} \leq \|f\|_{p_{k}} = 1$, $\max(\|f_{u}\|_{p_{k}}, \|f_{d}\|_{p_{k}}) \leq (1-c\eta^{p_{k}})\|f\|_{p_{k}} = 1-c\eta^{p_{k}}$. Hence: 
\begin{equation} \label{du}
    \|f_{u}\|_{p_{k}}^{2^{k}}+\|f_{d}\|_{p_{k}}^{2^{k}} \leq \max (\|f_{u}\|_{p_{k}}, \|f_{d}\|_{p_{k}})^{2^{k}-p_{k}} (\|f_{u}\|_{p_{k}}^{p_{k}}+\|f_{d}\|_{p_{k}}^{p_{k}}) \leq 1-c\eta^{p_{k}}.
\end{equation}
Inserting the bound in \eqref{du} into \eqref{Ukdu} and utilizing the fact that $f$ is a $(1-\delta)$ near extremizer:
\begin{equation*} 
    A(k,n)^{2^{k}}(1-\delta) \leq A(k,n)^{2^{k}}(1-c\eta^{p_{k}}) + C(k,n) 2^{-\eta^{p_{k}}M/N} + CN^{-1}\eta^{c}
\end{equation*}
or,
\begin{equation} \label{pkcon1}
    2^{-c\eta^{p_{k}}M/N}\geq c(k,n)\eta^{p_{k}} - c(k,n)N^{-1}-C(k,n)\delta.
\end{equation}
Now choose $N$ large enough to be the nearest strictly positive integer $\asymp c\eta^{-p_{k}}$ for some small $c$, so that from \eqref{pkcon1}:
\begin{equation} \label{pkcon2}
    2^{-c\eta^{2p_{k}}M} \geq c(k,n)\eta^{p_{k}} - C(k,n)\delta.
\end{equation}
Choose $\delta_0$ sufficiently small so that $c(k,n)\eta^{p_{k}} - C(k,n)\delta_0 > 0$, which yields, $\eta \geq C_0(k,n)\delta_0^{1/p_{k}}$ for some $C_0(k,n) >0$. We replace the first term in the right hand side of \eqref{pkcon2} with this new lower bound of $\eta$; the inequality then yields:
\begin{equation*} 
    M \leq C\eta^{-2p_{k}}
\end{equation*}
with $C = C(k,n,\delta)$, which is a desired bound on $M$. This completes the proof of statement \eqref{upbd} and consequently, the proof of Proposition \ref{prop:norm}. $\qed$

\section{Control of distribution functions}

\noindent
We wish to obtain two conclusions in this chapter; one is that, for a near extremizer, the contribution in $L^{p_{k}}$ norm of the superposition of its super-level sets associated with values outside a large compact range is negligible, and the other is, the measure of a super-level set with value in the said compact range closely approximates the measure of a super-level set of an extremizer of the same value. 

\subsection{Precompactness of symmetric rearrangements}

Note that if $f\geq 0$ is a $(1-\delta)$ near extremizer then so is its symmetric rearrangement, $f^{*}$ \cite{lieb2001analysis}:
\begin{align*} 
    \| f^{*}\|_{U^{k}}^{2^{k}} &= \int_{\mathbb{R}^{(k+1)n}}\prod_{\alpha\in\{0,1\}^{k}} f^{*}(x+\alpha\cdot\vec{h}) \, dxd\vec{h}\geq\int_{\mathbb{R}^{(k+1)n}}\prod_{\alpha\in\{0,1\}^{k}} f(x+\alpha\cdot\vec{h}) \, dxd\vec{h} \\ 
    &= \|f\|_{U^{k}}^{2^{k}} \geq (1-\delta)A(k,n)^{n2^{k}}\| f\|_{p_{k}}^{2^{k}} = (1-\delta)A(k,n)^{n2^{k}}\|f^{*}\|_{p_{k}}^{2^{k}}.
\end{align*}
The first inequality is due to the general rearrangement inequality and the second to the $L^{p_{k}}$-norm preservation property of symmetric rearrangements \cite{lieb2001analysis}. By definition, the normalization condition preserves measure. That means if $\Theta$ is a growth function and if $f\geq 0$ is normalized with respect to $\Theta$ then so is $f^{*}$:
\begin{align*}
    \int_{\{f^{*} > \rho\}} f^{*}(x) \, dx &= \int_{\{|f| > \rho\}} |f| (x) \, dx \leq \Theta(\rho)\\
    \int_{\{f^{*} < \rho\}} f^{*}(x) \, dx &= \int_{\{|f| < \rho\}} |f| (x) \, dx \leq \Theta(\rho^{-1}).
\end{align*}
Let $\{f_{i}\}_{i}$ be a sequence of measurable functions on $\mathbb{R}^{n}$. We say $\{f_{i}\}_{i}$ is a normalized extremizing sequence if $\|f_{i}\|_{p_{k}} =1$, $f_{i}$ is normalized with respect to $\Theta$ for all $i$ and if there exists a sequence $\{\delta_{i}\}_{i}$ satisfying $\delta_{i}\to 0$ as $i\to\infty$, such that:
\begin{equation*} 
    \| f_{i}\|_{U^{k}} \geq (1-\delta_{i})A(k,n)\|f_{i}\|_{p_{k}}.
\end{equation*}
Suppose furthermore that $f_{i}\geq 0$ for all $i$ and let $\{f^{*}_{i}\}_{i}$ be the corresponding symmetric rearrangement sequence of $\{f_{i}\}_{i}$. Then as noted, $\{ f^{*}_{i}\}_{i}$ is also a normalized extremizing sequence. 

\begin{lemma} \label{lem:rearrangcont}
Let $\{f^{*}_{i}\}_{i}$ be as above. Let $\epsilon > 0$. Then there exist $r> 0, R >0$ such that, for all sufficiently large $i$:
\begin{align}
    \label{lemrear1} \int_{\{|x| \leq r\}} (f^{*}_{i})^{p_{k}} \, dx &\leq \epsilon\\
    \label{lemrear2} \int_{\{|x| \geq R\}} (f^{*}_{i})^{p_{k}} \, dx &\leq \epsilon.
\end{align}
\end{lemma}

\begin{proof}
Indeed, by normalization, 
\begin{align} 
    \nonumber \int_{\{|x| \leq r\}} (f^{*}_{i})^{p_{k}} \, dx &\leq \int_{\{|x| \leq r\}\cap\{f^{*}_{i}\leq \alpha\}} (f_{i}^{*})^{p_{k}} \, dx+ \int_{\{f^{*}_{i} \geq \alpha\}} (f_{i}^{*})^{p_{k}} \, dx \\ 
    \label{fstari} &\leq c(n)\alpha^{p_{k}}r^{n} + \int_{\{f^{*}_{i}\geq \alpha\}} (f_{i}^{*})^{p_{k}} \, dx \leq c(n)\alpha^{p_{k}}r^{n} + \Theta (\alpha).
\end{align}
The third inequality follows from the fact that $f^{*}_{i}$ is normalized. Choose $\alpha$ sufficiently large so that $\Theta(\alpha) \leq (1/2)\epsilon$ and with this $\alpha$, $r$ small enough so that $c(n)\alpha^{p_{k}}r^{n} \leq (1/2)\epsilon$. Putting all these upper bounds in \eqref{fstari}, we obtain \eqref{lemrear1}. The fact that $f^{*}_{i}$ is radially symmetric and nonincreasing gives us,
\begin{equation*} 
    c(n)|s|^{n}(f_{i}^{*})^{p_{k}}(s) \leq \int_{\{|x| \leq s\}} (f_{i}^{*})^{p_{k}} \, dx \leq 1.
\end{equation*}
Hence, if $ R \leq |s| \leq R_0$:
\begin{equation} \label{fstaripk}
    (f_{i}^{*})^{p_{k}}(s) \leq C(n)|s|^{-n} \leq C(n) R^{-n}.
\end{equation}
In other words, $f_{i}^{*}$ is bounded above by $C(n)R^{-n/p_{k}}$ in an annulus $ \{ R\leq |x| \leq R_0\}$, with $R_0 > R$, uniformly for all $i$, and hence, $\{|x|\geq R\}\subset\{(f^{*}_{i})^{p_{k}}\leq C(n)R^{-n}\}$. Thus: 
\begin{equation} \label{fstaripktail}
    \int_{\{|x| \geq R\}} (f^{*}_{i})^{p_{k}} \, dx \leq \int_{\{(f^{*}_{i})^{p_{k}} \leq C(n)R^{-n}\}} (f^{*}_{i})^{p_{k}} \, dx \leq \Theta ((c(n)R^{n})^{1/p_{k}}).
\end{equation} 
The last inequality in \eqref{fstaripktail} follows from normalization. Now choose $R$ sufficiently large so that $\Theta ((c(n)R^{n})^{1/p_{k}}) \leq\epsilon$. Then \eqref{fstaripktail} implies \eqref{lemrear2}. This completes the proof of the lemma. 
\end{proof}

\noindent Let $r> 0, R> 0$ and $\mathcal{C}_{r,R} = \{ r \leq |x|\leq R\}$. The proof of Lemma \ref{lem:rearrangcont}, particularly \eqref{fstaripk}, shows that the sequence $\{ f^{*}_{i}\}_{i}$ restricted to a ray cutting the annulus $\mathcal{C}_{r,R}$, produces a uniformly bounded, nonincreasing sequence of functions on an interval $[r,R]$. It follows from Helly's selection principle \cite{rudin1976principles} that there exists a subsequence of $\{f_{i}^{*}\}_{i}$ that converges pointwise on every such ray. Since $f^{*}_{i}$ is radial for all $i$, such a subsequence also converges pointwise on the annulus $\mathcal{C}_{r,R}$. We denote such a subsequence by $\{ f^{*}_{i;r,R}\}_{i}$ and the corresponding pointwise limit by $\mathcal{F}_{r,R}$. It's clear that $\mathcal{F}_{r,R}$ is also radial. We select $r$ to be an increasing function of $\epsilon$ and $R$ a decreasing function of $\epsilon$. Choose a sequence $\{\epsilon_{i}\}_{i}$ that decreases to zero and denote $r(\epsilon_{i}) = r_{i}, R(\epsilon_{i}) = R_{i}$. Then by a Cantor diagonal argument and passing to a subsequence of $\{f^{*}_{i;r_{i},R_{i}}\}_{i}$, we obtain a sequence of nonincreasing sequence of radial functions $\{\mathcal{F}_{r_{i},R_{i}}\}_{i}$ such that $\mathcal{F}_{r_{i},R_{i}}$ is supported on the annulus $\mathcal{C}_{r_{i},R_{i}}$, for all $i$, and $\mathcal{F}_{r_{i}, R_{i}} = \mathcal{F}_{r_{n},R_{n}}$ on $\mathcal{C}_{r_{n},R_{n}}$ if $i\geq n$. Let $\mathcal{F}$ be the pointwise limit of $\{\mathcal{F}_{r_{i},R_{i}}\}_{i}$. Note that $\mathcal{F}_{r_{n}, R_{n}} \leq \mathcal{F}_{r_{i},R_{i}}\leq \mathcal{F}$, if $i \geq n$. For convenience, we denote the subsequence of $\{f^{*}_{i;r_{i},R_{i}}\}_{i}$ in the construction of $\mathcal{F}_{r_{i},R_{i}}$ at every index $i$ as simply $\{ f^{*}_{i}\}_{i}$. We argue that $\mathcal{F}$ must be a centered Gaussian function on $\mathbb{R}^{n}$. Indeed, by Lemma \ref{lem:rearrangcont}, for all sufficiently large $i$,
\begin{equation} \label{doubletailcon}
    \int_{\{|x|\leq r_{i}\}} (f^{*}_{i})^{p_{k}} \, dx + \int_{\{|x| \geq R_{i}\}} (f^{*}_{i})^{p_{k}} \, dx \leq 2\epsilon_{i}.
\end{equation}
Take $\eta > 0$. Then take $I$ sufficiently large so that $2\epsilon_{i} < \eta$, if $i \geq I$. By letting $i\to\infty$, it follows from \eqref{doubletailcon} and the fact that $\mathcal{F}_{r_{i}, R_{i}} = \mathcal{F}_{r_{n},R_{n}}$ on $\mathcal{C}_{r_{n},R_{n}}$ if $i\geq n\geq I$,
\begin{align*} 
    \int_{\mathbb{R}^{n}} |\mathcal{F}_{r_{n},R_{n}} - \mathcal{F}_{r_{i},R_{i}}|^{p_{k}} \, dx &\leq \int_{\{|x| \leq r_{I}\}} |\mathcal{F}_{r_{n},R_{n}} - \mathcal{F}_{r_{i},R_{i}}|^{p_{k}} \\ 
    & + \int_{\{|x| \geq R_{I}\}}|\mathcal{F}_{r_{n},R_{n}} - \mathcal{F}_{r_{i},R_{i}}|^{p_{k}} + \int_{\{r_{I}\leq |x| \leq R_{I}\}} |\mathcal{F}_{r_{n},R_{n}} - \mathcal{F}_{r_{i},R_{i}}|^{p_{k}} \, dx \leq \eta. 
\end{align*}
Hence $\{\mathcal{F}_{r_{i},R_{i}}\}_{i}$ is a Cauchy sequence in $L^{p_{k}}$ and thus $\mathcal{F}$ is also its $L^{p_{k}}$-limit. On the other hand, by the Dominated Convergence Theorem, $\mathcal{F}_{r_{i},R_{i}}$ is the pointwise limit and thus the $L^{p_{k}}$-limit of $\{f^{*}_{i}\}_{i}$ on $\mathcal{C}_{r_{i},R_{i}}$. Hence $\mathcal{F}$ is the $L^{p_{k}}$-limit of $\{f^{*}_{i}\}_{i}$ on $\mathbb{R}^{n}$. By the Gowers-Host-Kra norm inequality, $\| \mathcal{F}-f^{*}_{i}\|_{U^{k}} \leq A(k,n) \| \mathcal{F}-f^{*}_{i}\|_{p_{k}}$, and hence, as $i\to\infty$, $\|f^{*}_{i}\|_{U^{k}} \to\|\mathcal{F}\|_{U^{k}}$. But since $\| f^{*}_{i} \|_{U^{k}} \geq (1-\delta_{i})A(k,n)\|f^{*}_{i}\|_{p_{k}}$, these imply that $\|\mathcal{F} \|_{U^{k}}= A(k,n)\|\mathcal{F}\|_{p_{k}}$. By the characterization of extremizers of the Gowers-Host-Kra norm inequality, $\mathcal{F}$ must be a Gaussian, which then is a centered Gaussian, as it is the pointwise limit of radial functions.

\subsection{Control of distribution functions}

Recall that if $f:\mathbb{R}^{n}\to\mathbb{R}$ is a $(1-\delta)$ near extremizer then so is $|f|$. Hence we only consider nonnegative near extremizers for now. Suppose $\| f \|_{U^{k}} \geq (1-\delta)A(k,n)\|f\|_{p_{k}}$ and $ \| \mathcal{F}- f^{*} \|_{p_{k}} \leq \delta$ with $\mathcal{F}$ being a centered Gaussian on $\mathbb{R}^{n}$. Denote, for $s>0$, $F_{s} = \{ f > s\}$, $F^{*}_{s} = \{ f^{*} > s\}$ and $\mathcal{F}_{s} = \{ \mathcal{F} > s\}$.

\begin{lemma} \label{lem:fandrear}
Let $f, f^{*}, \mathcal{F}, F_{s}, F^{*}_{s}, \mathcal{F}_{s}$ be as above. Given $\delta > 0$. There exists $\eta = \eta(\delta) >0$ satisfying $\eta(\delta)\to 0$ as $\delta\to 0$, such that if $s \in [\eta,\| \mathcal{F} \|_{\infty}-\eta]$ then
\begin{align}
    \label{lemfandrear1} \mathcal{L}(\mathcal{F}_{s}\Delta F^{*}_{s}) &= o(\delta)\\
    \label{lemfandrear2} |\mathcal{L}(\mathcal{F}_{s}) - \mathcal{L}(F_{s})| &= o(\delta).
\end{align}
\end{lemma}

\begin{proof}
Let $\eta \in (0,1)$ be a small number to be chosen below. Consider $ F^{*}_{s} \setminus \mathcal{F}_{s - \eta}$. If $x \in F^{*}_{s} \setminus \mathcal{F}_{s-\eta}$ then $|(\mathcal{F}-f^{*})(x)| \geq \eta$. Chebyshev's inequality gives, $\mathcal{L}(F^{*}_{s} \setminus \mathcal{F}_{s-\eta} )\leq \eta^{-p_{k}} \| \mathcal{F}- f^{*} \|_{p_{k}}^{p_{k}} \leq \eta^{-p_{k}} \delta^{p_{k}}$, and, $\mathcal{L}(\mathcal{F}_{s+\eta} \setminus F^{*}_{s}) \leq \eta^{-p_{k}}\delta^{p_{k}}$. Hence:
\begin{equation} \label{diff1}
    \mathcal{L}(\mathcal{F}_{s+\eta} \Delta \mathcal{F}_{s-\eta}) \leq c\eta^{-p_{k}} \delta^{p_{k}}.
\end{equation}
On the other hand, since $\mathcal{F}$ is a Gaussian and all the super-level sets of $\mathcal{F}$ are nested, centered ellipsoids, we also have
\begin{equation} \label{diff2}
    \mathcal{L}(\mathcal{F}_{s+\eta} \Delta \mathcal{F}_{s-\eta}) \leq C\eta^{1/2}.
\end{equation}
Optimizing \eqref{diff1} and \eqref{diff2} to have:
\begin{equation} \label{diff3}
    \mathcal{L}(\mathcal{F}_{s+\eta} \Delta \mathcal{F}_{s-\eta}) \leq C\delta^{p_{k}/(2p_{k}+1)}.
\end{equation}
Since the super-level sets of $f^{*}$ are also centered ellipsoids, if $\eta = \delta^{2p_{k}/(2p_{k}+1)}$ then
\begin{align}
    \label{diff4} \mathcal{L}(\mathcal{F}_{s-\eta}\Delta F^{*}_{s}) &\leq c\delta^{p_{k}/(2p_{k}+1)}\\
    \label{diff5} \mathcal{L}(\mathcal{F}_{s+\eta}\Delta F^{*}_{s}) &\leq c\delta^{p_{k}/(2p_{k}+1)}.
\end{align}
\eqref{diff3}, \eqref{diff4}, \eqref{diff5} then imply $\mathcal{L}(\mathcal{F}_{s} \Delta F^{*}_{s}) \leq C \delta^{p_{k}/(2p_{k}+1)}$, and consequently, $| \mathcal{L}(\mathcal{F}_{s}) - \mathcal{L}(F_{s}) | \leq C \delta^{p_{k}/(2p_{k}+1)}$, which are the promises \eqref{lemfandrear1} and \eqref{lemfandrear2}, respectively.
\end{proof}

\begin{remark} \label{rem:lem5.2}
We can summarize the findings above in the following language. Given a $\delta$ small, there exist $\eta(\delta)$ and $\epsilon(\eta,\delta)$ with the following properties. If $\delta\to 0$ then $\eta(\delta)\to 0$, and if $\eta$ is fixed, $\delta\to 0$ then $\epsilon(\eta,\delta)\to 0$. If $\mathcal{F}_{s}, F_{s}$ are as above, then $\sup_{s\in [\eta(\delta), \| \mathcal{F}_{s}\|_{\infty} - \eta(\delta)]} |\mathcal{L}(\mathcal{F}_{s}) - \mathcal{L}(F_{s})| \leq \epsilon (\eta, \delta)$. Moreover, by the description of $\epsilon(\eta, \delta)$, we have that if $\xi(\eta)$ is a positive continuous function satisfying $\xi(\eta)\to 0$ as $\eta\to 0$, then there exists a function $\eta\mapsto\delta_0(\eta) > 0$ such that $\sup_{0<\delta \leq\delta_0(\eta)} \epsilon(\eta,\delta) \leq \xi(\eta)$. In particular, $\lim_{\eta\to 0} \sup_{0<\delta\leq \delta_0(\eta)} \epsilon(\eta,\delta) = 0$, for some positive function $\delta_0$. These properties of $\epsilon(\eta,\delta)$ will be needed below in Section \ref{sec:superlevel}.    
\end{remark}

\begin{lemma} \label{lem:levelsetstack}
Let $f^{*}, \mathcal{F}, F_{s}$ and $\eta$ be as above. There exists $C > 0$ so that,
\begin{equation*} 
    \bigg \| \int_0^{\eta} 1_{F_{s}} \, ds \bigg \|_{p_{k}} \leq C \| \mathcal{F} - f^{*} \|_{p_{k}} + C\eta (\log (1/\eta))^{C}
\end{equation*}
and,
\begin{equation*} 
    \bigg \| \int_{\| \mathcal{F}\|_{\infty}- \eta}^{\infty} 1_{F_{s}} \, ds  \bigg \|_{p_{k}} \leq C\| \mathcal{F} - f^{*} \|_{p_{k}} + O(\eta).
\end{equation*}    
\end{lemma}

\begin{proof}
Indeed, 
\begin{align*} 
    \bigg \| \int_0^{\eta} 1_{F_{s}} \, ds \bigg \|_{p_{k}} = \bigg \| \int_0^{\eta} 1_{F^{*}_{s}} \, ds \bigg \|_{p_{k}}= \| \min(f, \eta)\|_{p_{k}} &\leq \|\mathcal{F} - f^{*}\|_{p_{k}} + \| \min (\mathcal{F}, \eta)\|_{p_{k}} \\
    &\leq \| \mathcal{F} - f^{*}\|_{p_{k}} + \eta (\log (1/\eta))^{C}.
\end{align*}
Likewise: 
\begin{align*} 
    \bigg \| \int_{\| \mathcal{F}\|_{\infty} - \eta}^{\infty} 1_{F_{s}} \, ds \bigg \|_{p_{k}} &= \bigg \| \int_{\| \mathcal{F} \|_{\infty} - \eta}^{\infty} 1_{F^{*}_{s}} \, ds \bigg \|_{p_{k}} = \| \max (0, f^{*} - (\| \mathcal{F}\|_{\infty} - \eta)\|_{p_{k}} \\ 
    &\leq \|\mathcal{F}-f^{*}\|_{p_{k}} + \| \max (0,\mathcal{F} - (\| \mathcal{F}\|_{\infty} - \eta)) \|_{p_{k}}\leq \| \mathcal{F}-f^{*}\|_{p_{k}} + O(\eta).
\end{align*}
The proof is complete.     
\end{proof}

\begin{remark}
One can appreciate the usefulness of the normalization condition in obtaining the decay properties \eqref{lemrear1} and \eqref{lemrear2} and consequently the precompactness of $\{f^{*}_{i}\}_{i}$ in $L^{p_{k}}$, which then allows us to extract an extremizer $\mathcal{F}$. Moreover, one observes that the only needed hypotheses for Lemma \ref{lem:fandrear} and Lemma \ref{lem:levelsetstack} are that $f^{*}$ is radially symmetric and decreasing, $\mathcal{F}$ is a centered Gaussian, and $\|\mathcal{F}-f\|_{p_{k}}\leq\delta, \|\mathcal{F}-f^{*}\|_{p_{k}}\leq\delta$; the specific fact that $f$ is a near extremizer for the Gowers-Host-Kra norm inequality did not enter. This observation was first used by Christ to obtain similar conclusions about near extremizing triple $(f,g,h)$ of Young's inequality on Euclidean spaces \cite{christ2019near}.   
\end{remark}

\section{Control of super-level sets} \label{sec:superlevel}

\noindent
Let $n=1$. Let $f$ be a $(1-\delta)$ near extremizer of the Gowers-Host-Kra norm inequality and $\mathcal{F}$ a centered Gaussian such that $\|f\|_{p_{k}} = 1$ and $\| \mathcal{F}-f^{*}\|_{p_{k}} \leq \delta$. In this chapter, we only consider $f\geq 0$. Since $f$ is of unit norm, so are $f^{*}$ and $\mathcal{F}$. For $t > 0$, let $\Omega(t)= [t, \| \mathcal{F}\|_{\infty} - t]$. Recall that for $s> 0$, $F_{s} = \{ f > s\}, F^{*}_{s} = \{ f^{*} > s\}, \mathcal{F}_{s} = \{ \mathcal{F} > s\}$.

\begin{proposition} \label{prop:lwbrearr}
For every $\epsilon > 0$ there exists $\delta > 0$ such that for every $(1-\delta)$ near extremizer $f$ and every $\nu_{\alpha}\in\Omega(\epsilon)$, $\alpha\in\{0,1\}^{k}$,
\begin{equation} \label{lwbrearr1}
    \mathcal{T}_{k}(1_{F_{\nu_{\alpha}}}: \alpha\in\{0,1\}^{k}) \geq (1-\epsilon) \mathcal{T}_{k}(1_{F^{*}_{\nu_{\alpha}}}:\alpha\in\{0,1\}^{k}).
\end{equation}    
\end{proposition}

\begin{remark}
We will prove the statement in the following form:\\
For every $\epsilon > 0$ there exists $\zeta > 0$ and $\delta >0$ such that if $f$ is a $(1-\delta)$ near extremizer of the Gowers-Host-Kra norm inequality and if $\nu_{\alpha} \in \Omega(\zeta)$, $\alpha\in\{0,1\}^{k}$ then,
\begin{equation} \label{lwbrearr2}
    \mathcal{T}_{k}(1_{F_{\nu_{\alpha}}}:\alpha\in\{0,1\}^{k}) \geq (1-\epsilon)\mathcal{T}_{k}(1_{F^{*}_{\nu_{\alpha}}}:\alpha\in\{0,1\}^{k}).
\end{equation}
\eqref{lwbrearr2} is equivalent to \eqref{lwbrearr1} if $\zeta = o(\epsilon)$ and $\epsilon = o(\zeta)$, which will be the case. 
\end{remark}

\noindent
Assume Proposition \ref{prop:lwbrearr} for a moment. Assume $\nu_{\alpha} = \nu$ for all $\alpha\in\{0,1\}^{k}$. Then by definition, $\mathcal{T}_{k}(1_{F_{\nu_{\alpha}}}:\alpha\in\{0,1\}^{k}) = \mathcal{T}_{k}(1_{F_{\nu}}:\alpha\in\{0,1\}^{k}) = \| 1_{F_{\nu}}\|_{U^{k}}^{2^{k}}$ and $\mathcal{T}_{k}(1_{F^{*}_{\nu}}:\alpha\in\{0,1\}^{k}) = \|1_{F^{*}_{\nu}}\|_{U^{k}}^{2^{k}}$. Then \eqref{lwbrearr2} becomes $\| 1_{F_{\nu}}\|_{U^{k}}\geq (1-\epsilon)\|1_{F^{*}_{\nu}}\|_{U^{k}}$, which entails that the sets $F_{\nu_{\alpha}}$ will then be nearly intervals, by the following result:

\begin{proposition} \label{prop:christ} \cite{christ2019subsets} Let $n\geq 1$. For every $\epsilon > 0$ there exists $\delta > 0$ such that if $E\subset \mathbb{R}^{n}$ is a measurable set with finite measure and $\| 1_{E}\|_{U^{k}} \geq (1-\delta) \|1_{E^{*}}\|_{U^{k}}$ then there exists an ellipsoid $\mathcal{E}\subset\mathbb{R}^{n}$ such that $\mathcal{L}(E \Delta\mathcal{E}) < \epsilon\mathcal{L}(E)$.     
\end{proposition}

\noindent
In what follows, $\vec{\nu} = (\nu_{\alpha}:\alpha\in\{0,1\}^{k})$.

\subsection{Proof of Proposition \ref{prop:lwbrearr}}

\subsubsection{Set-up:}

\noindent
Let $t\in (0,1)$. By the Gowers product inequality
\begin{equation} \label{Gprineq}
    \int_{\mathbb{R}^{\{0,1\}^{k}}\setminus\Omega(t)^{\{0,1\}^{k}}} \mathcal{T}_{k}(1_{F_{\nu_{\alpha}}}:\alpha\in\{0,1\}^{k}) \, d\vec{\nu} \leq A(k)^{2^{k}} \int_{\mathbb{R}^{\{0,1\}^{k}}\setminus\Omega(t)^{\{0,1\}^{k}}} \prod_{\alpha\in\{0,1\}^{k}} \| 1_{F_{\nu_{\alpha}}}\|_{p_{k}} \, d\vec{\nu}.
\end{equation}
Moreover $\mathbb{R}^{\{0,1\}^{k}}\setminus\Omega(t)^{\{0,1\}^{k}}\subset \cup_{\alpha\in\{0,1\}^{k}} (\mathbb{R}\setminus\Omega(t)) \times \mathbb{R}^{\{0,1\}^{k}-1}$. Hence an immediate consequence of Lemma \ref{lem:levelsetstack}, Fubini's theorem and \eqref{Gprineq} is that there exists a universal $C(k) >0$ such that for all small $t > 0$:
\begin{multline} \label{Cklog}
    \int_{\mathbb{R}^{\{0,1\}^{k}}\setminus\Omega(t)^{\{0,1\}^{k}}} \mathcal{T}_{k}(1_{F_{\nu_{\alpha}}}:\alpha\in\{0,1\}^{k}) \, d\vec{\nu} \\ \leq A(k)^{2^{k}} \prod_{\alpha\in\{0,1\}^{k}} \bigg \| \int_{\mathbb{R}\setminus\Omega(t)} 1_{F_{\nu_{\alpha}}} \, d\nu_{\alpha} \bigg \|_{p_{k}} \leq C(k)\delta + C(k)t(\log (1/t))^{C}. 
\end{multline}
Define:
\begin{align*}
    H(\nu_{\vec{0}}) &= \mathcal{T}_{k}((1_{F_{\nu_{\vec{0}}}}, 1_{F_{\nu_{\alpha}}}): \alpha\in\{0,1\}^{k};\alpha\not=\vec{0})\\
    \mathcal{H}(\nu_{\vec{0}}) &= \mathcal{T}_{k}((1_{\mathcal{F}_{\nu_{\vec{0}}}}, 1_{\mathcal{F}_{\nu_{\alpha}}}): \alpha\in\{0,1\}^{k};\alpha\not=\vec{0}).
\end{align*}
In other words, we fix all the values $\nu_{\alpha}, \alpha\not=\vec{0}$ and consider $H, \mathcal{H}$ as nonnegative functions of only $\nu_{\vec{0}}$. Since level sets are nested, both $H$ and $\mathcal{H}$ are non-increasing. Moreover, since $\mathcal{F}$ is a centered Gaussian of unit norm, if $\nu_{\alpha}\in\Omega(t)$, for all $\alpha\in\{0,1\}^{k}$, then $\mathcal{H}$ is bounded below by a strictly positive function of $t$ as long as $0 <t< \| \mathcal{F}\|_{\infty}$. Similarly, $\mathcal{T}_{k}(1_{\mathcal{F}_{\nu_{\alpha}}}:\alpha\in\{0,1\}^{k}) \geq \phi(t)$, with $\phi$ being a strictly positive continuous function satisfying $\phi\to 0$ as $t\to 0$; this function $\phi$ will be needed below. Moreover, $\mathcal{H}$ is also Lipschitz continuous with a Lipschitz constant $L(t)$ that is independent of a specific $\vec{\nu}$, as long as $\nu_{\alpha}\in\Omega(t)$, for all $\alpha\in\{0,1\}^{k}$. Now:
\begin{align}
    \nonumber \int_{\Omega(t)^{\{0,1\}^{k}}} \mathcal{T}_{k}(1_{F_{\nu_{\alpha}}}: &\alpha\in\{0,1\}^{k}) \, d\vec{\nu} \\ 
    \nonumber = &\| f\|_{U^{k}}^{2^{k}} - \int_{\mathbb{R}^{\{0,1\}^{k}}\setminus\Omega(t)^{\{0,1\}^{k}}} \mathcal{T}_{k}(1_{F_{\nu_{\alpha}}}:\alpha\in\{0,1\}^{k}) \, d\vec{\nu} \\ 
    \nonumber \geq & (1-\delta) \| \mathcal{F}\|_{U^{k}} - C(k)\delta - C(k)t (\log(1/t))^{C} \\ 
    \label{Cklogpf} \geq &(1-\delta) \int_{\Omega(t)^{\{0,1\}^{k}}} \mathcal{T}_{k}(1_{\mathcal{F}_{\nu_{\alpha}}}:\alpha\in\{0,1\}^{k}) \, d\vec{\nu} - C(k)\delta - C(k)t(\log(1/t))^{C}.    
\end{align}
The first inequality follows from \eqref{Cklog} and the near extremizing hypothesis. The second inequality follows from the definition of the Gowers-Host-Kra norms. Recall from Remark \ref{rem:lem5.2}, that if $s\in \Omega(t)$, then:
\begin{equation} \label{compsets}
    \sup_{s\in\Omega(t)} |\mathcal{L}(\mathcal{F}_{s}) - \mathcal{L}(F_{s})| \leq \epsilon(t,\delta)
\end{equation}
and $\epsilon(t,\delta)$ has the following properties:
\begin{enumerate}
    \item $\epsilon(t,\delta) \to 0$ if $\delta\to 0$ and $t$ is fixed.
    \item For every positive continuous function $\xi(t)$ satisfying $\xi(t)\to 0$ as $t\to 0$, there exists $\delta_0(t) > 0$ such that $\epsilon(t,\delta) \leq \xi(t)$ for all $0 < \delta \leq \delta_0(t)$; hence $\lim_{t\to 0} \sup_{0<\delta \leq\delta_0(t)} \epsilon(t,\delta)= 0$.
\end{enumerate}
It follows from \eqref{compsets} and the Gowers product inequality, that if $\vec{\nu}\in \Omega(t)^{\{0,1\}^{k}}$ then:
\begin{multline} \label{compsetsTk}
    | \mathcal{T}_{k}(1_{\mathcal{F}_{\nu_{\alpha}}}:\alpha\in\{0,1\}^{k}) - \mathcal{T}_{k}(1_{F^{*}_{\nu_{\alpha}}}: \alpha\in\{0,1\}^{k})| \\ \leq C(k) \sup_{s\in\Omega(t)} | (\mathcal{L}(\mathcal{F}_{s}))^{1/p_{k}} - (\mathcal{L}(F_{s}))^{1/p_{k}}| \cdot ((\mathcal{L}(\mathcal{F}_{t}))^{1/p_{k}} + \epsilon(t,\delta))^{2^{k}-1} \leq C(k,t)\epsilon(t,\delta). 
\end{multline}
The constant $C(k,t)$ depends on $\mathcal{L}(\mathcal{F}_{t})$ and satisfies $C(k,t)\to\infty$ if $t\to 0$. Combining \eqref{compsetsTk} and the general rearrangement inequality, we obtain:
\begin{equation} \label{compsetsup}
    \mathcal{T}_{k}(1_{F_{\nu_{\alpha}}}:\alpha\in\{0,1\}^{k}) \leq \mathcal{T}_{k}(1_{F^{*}_{\nu_{\alpha}}}: \alpha\in\{0,1\}^{k}) \leq \mathcal{T}_{k}(1_{\mathcal{F}_{\nu_{\alpha}}}:\alpha\in\{0,1\}^{k}) +C(k,t)\epsilon(t,\delta).
\end{equation}
Take two numbers $\rho,\eta\in (0,1)$ satisfying $\rho\leq\eta$, then $\Omega(\eta)\subset\Omega(\rho)$. Integrating $\vec{\nu}$ in \eqref{compsetsup} over $\Omega(\rho)^{\{0,1\}^{k}}\setminus\Omega(\eta)^{\{0,1\}^{k}}$ gives:
\begin{multline} \label{compsetsup1}
    \int_{\Omega(\rho)^{\{0,1\}^{k}}\setminus\Omega(\eta)^{\{0,1\}^{k}}} \mathcal{T}_{k}(1_{F_{\nu_{\alpha}}}:\alpha\in\{0,1\}^{k}) \, d\vec{\nu} \\ \leq \int_{\Omega(\rho)^{\{0,1\}^{k}}\setminus\Omega(\eta)^{\{0,1\}^{k}}} \mathcal{T}_{k}(1_{\mathcal{F}_{\nu_{\alpha}}}: \alpha\in\{0,1\}^{k}) \, d\vec{\nu} + C(k,\rho)\epsilon(\rho,\delta).
\end{multline}
Here again $C(k,\rho)\to\infty$ if $\rho\to 0$. Substituting $t=\rho$ in \eqref{Cklogpf} to have:
\begin{multline} \label{compsetslw}
    \int_{\Omega(\rho)^{\{0,1\}^{k}}} \mathcal{T}_{k}(1_{F_{\nu_{\alpha}}}: \alpha\in\{0,1\}^{k}) \, d\vec{\nu} \geq (1-\delta) \int_{\Omega(\rho)^{\{0,1\}^{k}}} \mathcal{T}_{k}(1_{\mathcal{F}_{\nu_{\alpha}}}:\alpha\in\{0,1\}^{k}) \, d\vec{\nu} \\ - C(k)\delta - C(k)\rho(\log(1/\rho))^{C}.
\end{multline}
Subtracting \eqref{compsetsup1} from \eqref{compsetslw}:
\begin{multline} \label{compsetsdiff}
    \int_{\Omega(\eta)^{\{0,1\}^{k}}} \mathcal{T}_{k}(1_{F_{\nu_{\alpha}}}:\alpha\in\{0,1\}^{k}) \, d\vec{\nu} \geq \int_{\Omega(\eta)^{\{0,1\}^{k}}} \mathcal{T}_{k}(1_{\mathcal{F}_{\nu_{\alpha}}}:\alpha\in\{0,1\}^{k}) \, d\vec{\nu} \\ - C(k) \delta\phi(\eta)- C(k)\delta -C(k)\rho(\log (1/\rho))^{C} - C(k,\rho)\epsilon(\rho,\delta).
\end{multline}
Let $\xi(t)$ be a positive continuous function satisfying $\xi(t)\to 0$ as $t\to 0$. We now require in our selection of $\rho,\eta$, that $C(k)\rho(\log(1/\rho))^{C} + C(k)\phi(\eta) \leq \xi(\eta)$. By the second property of $\epsilon(t,\delta)$, there exists a positive function $\delta_0(t)$ such that for all $0 < \delta \leq \delta_0(t)$, $t\in(0,1)$, the two following conditions are satisfied:
\begin{align*}
    \epsilon(t,\delta) &\leq \xi(t)\\
    \epsilon(t,\delta) &\leq C(k,t)^{-1} \xi(t)
\end{align*}
with $C(k,\cdot)$ having the same meaning as in \eqref{compsetsdiff}. Since $\xi(t)\to 0$ as $t\to 0$, $c_1\xi(\rho)\leq\xi(\eta)\leq c_2\xi(\rho)$, if $\rho\leq\eta$ are both sufficiently small. These observations together with the conditions allow us to rewrite \eqref{compsetsdiff} as:
\begin{equation} \label{compsetsTkint}
    \int_{\Omega(\eta)^{\{0,1\}^{k}}} (\mathcal{T}_{k}(1_{\mathcal{F}_{\nu_{\alpha}}}: \alpha\in\{0,1\}^{k}) - \mathcal{T}_{k}(1_{F_{\nu_{\alpha}}}:\alpha\in\{0,1\}^{k})) \, d\vec{\nu} \leq C(k)\delta + C(k)\xi(\eta).
\end{equation}
Note that if we re-define $\epsilon(s,\delta) = C(k)\delta + C(k)\xi(s)$, then this new $\epsilon(s,\delta)$ still satisfies the two properties of the error function $\epsilon$ mentioned above. We rewrite \eqref{compsetsup} and \eqref{compsetsTkint} respectively in this new language, with $\nu_{\alpha}\in\Omega(t)$, $\alpha\in\{0,1\}^{k}$:
\begin{align}
    \label{compsets1side} \mathcal{T}_{k}(1_{F_{\nu_{\alpha}}}:\alpha\in\{0,1\}^{k}) \leq \mathcal{T}_{k}(1_{F^{*}_{\nu_{\alpha}}}:\alpha\in\{0,1\}^{k}) \leq \mathcal{T}_{k}(1_{\mathcal{F}_{\nu_{\alpha}}}:\alpha\in\{0,1\}^{k}) +C(k)\epsilon(t,\delta) &\\
    \label{compsetsinteps} \int_{\Omega(t)^{\{0,1\}^{k}}} (\mathcal{T}_{k}(1_{\mathcal{F}_{\nu_{\alpha}}}: \alpha\in\{0,1\}^{k}) - \mathcal{T}_{k}(1_{F_{\nu_{\alpha}}}:\alpha\in\{0,1\}^{k})) \, d\vec{\nu} \leq \epsilon(t,\delta) &.
\end{align}
By Fubini's theorem, \eqref{compsetsinteps} is simply,
\begin{multline} \label{Hupb}
    \int_{\Omega(t)^{\{0,1\}^{k}}} (\mathcal{T}_{k}(1_{\mathcal{F}_{\nu_{\alpha}}}:\alpha\in\{0,1\}^{k}) - \mathcal{T}_{k}(1_{F_{\nu_{\alpha}}}:\alpha\in\{0,1\}^{k})) \, d\vec{\nu} \\ = \int_{\Omega(t)^{\{0,1\}^{k}}} (\mathcal{H}(\nu_{\vec{0}}) - H(\nu_{\vec{0}})) \, d\vec{\nu} \leq \epsilon(t,\delta). 
\end{multline}
\eqref{compsets1side} implies, if $\nu_{\alpha}\in\Omega(t)$, for all $\alpha\in\{0,1\}^{k}$:
\begin{equation} \label{Hlwb}
    \mathcal{H}(\nu_{\vec{0}}) - H(\nu_{\vec{0}}) \geq - C(k)\epsilon(t,\delta).
\end{equation}
Integrating \eqref{Hlwb} over $\Omega(t)$, we have uniformly for $(\nu_{\alpha}: \alpha\in\{0,1\}^{k};\alpha\not=\vec{0}) \in\Omega(t)^{\{0,1\}^{k}\setminus\{\vec{0}\}}$:
\begin{equation} \label{Hlwbint}
    \int_{t}^{\|\mathcal{F}\|_{\infty}-t} (\mathcal{H}(\nu_{\vec{0}}) - H(\nu_{\vec{0}})) \, d\nu_{\vec{0}} \geq - C(k)\epsilon(t,\delta). 
\end{equation}

\subsubsection{A process:}

\noindent
Recall that $\mathcal{H}$ is Lipschitz continuous with a Lipschitz constant majorized by a quantity $L(t)$ that is independent of $\vec{\nu}$, as long as $\nu_{\alpha}\in\Omega(t)$, for all $\alpha\in\{0,1\}^{k}$. Let $K= K(t) = \max (1,L(t))$. Note that $K(t)\to\infty$ as $t\to 0$, since $\mathcal{L}(\mathcal{F}_{\|\mathcal{F}\|_{\infty}-t})\to 0$ and $\mathcal{L}(\mathcal{F}_{t}) \to \infty$ as $t\to 0$. So far we have two parameters $\rho\leq\eta$. The parameter $\rho$ is an auxiliary one whose ultimate use was to define the new error function $\epsilon$; our main parameter is the parameter $\eta$. We now will define new parameters in terms of $\eta$. Let $r=r(\eta)$ be a small quantity to be chosen below, and suppose that there exists $\nu_{\vec{0}}' \in [\eta,\|\mathcal{F}\|_{\infty} - \eta -r]$ such that:
\begin{equation} \label{HtildeH}
    H(\nu_{\vec{0}}') \leq \mathcal{H}(\nu_{\vec{0}}') -r.
\end{equation}
Then, we claim, a similar property will hold for a sub-range of $\nu_{\vec{0}} \in [\nu_{\vec{0}}', \nu_{\vec{0}}' + cK^{-1}r]$:
\begin{equation} \label{tildeHcr}
    H(\nu_{\vec{0}})\leq H(\nu_{\vec{0}}') \leq \mathcal{H}(\nu_{\vec{0}}') -r\leq \mathcal{H}(\nu_{\vec{0}}) - cr.
\end{equation}
Indeed, the non-increasing property of $H$ gives the first inequality in \eqref{tildeHcr}; the second is just \eqref{HtildeH}, and the Lipschitz continuity of $\mathcal{H}$ over the selected interval gives the last. Furthermore, if \eqref{HtildeH} happens, we can increase the lower bound in \eqref{Hlwbint}:
\begin{align} 
    \nonumber \int_{\eta}^{\|\mathcal{F}\|_{\infty}-\eta} (\mathcal{H}(\nu_{\vec{0}})  - H(\nu_{\vec{0}})) \, d\nu_{\vec{0}} &= \int_{[\nu_{\vec{0}}', \nu_{\vec{0}}' + cK^{-1}r]} + \int_{[\eta, \|\mathcal{F}\|_{\infty}-\eta]\setminus [\nu_{\vec{0}}', \nu_{\vec{0}}' + cK^{-1}r]}(\mathcal{H}(\nu_{\vec{0}})  - H(\nu_{\vec{0}})) \, d\nu_{\vec{0}}\\ 
    \label{tildeHtail} &\geq cr^2K^{-1} - C(k)\epsilon(\eta,\delta). 
\end{align}
The inequality in \eqref{tildeHtail} follows from applying the respective lower bounds as in \eqref{Hlwbint}, with $t=\eta$, and \eqref{tildeHcr}. It will be in our desire that $r(\eta)\to 0$ as $\eta\to 0$ but at a rate $o(\eta)$. To this end, take $q\in (0,1)$. We consider only sufficiently small values of $\delta$ so that $qcK^{-1}r^2 \geq C\epsilon(\eta,\delta)$. This is possible by the first property of $\epsilon$ described above. This allows us to rewrite \eqref{tildeHtail} as:
\begin{equation} \label{tildeHtailfin}
    \int_{\eta}^{\|\mathcal{F}\|_{\infty} -\eta} (\mathcal{H}(\nu_{\vec{0}}) - H(\nu_{\vec{0}})) \, d\nu_{\vec{0}} \geq c(k)K^{-1}r^2. 
\end{equation}
Let $\mathcal{S}$ be the set of $(\nu_{\alpha}: \alpha\in\{0,1\}^{k}; \alpha\not=\vec{0})\in\Omega(\eta)^{\{0,1\}^{k}\setminus\{\vec{0}\}}$ for which there exists at least one $\nu_{\vec{0}}'$ such that \eqref{HtildeH} is satisfied. By \eqref{Hupb}, with $t=\eta$, \eqref{tildeHtailfin} and Markov's inequality:
\begin{equation} \label{markov}
    \mathcal{L}(\mathcal{S}) \leq C(k)Kr^{-2}\epsilon(\eta,\delta).
\end{equation}
Select the parameter $\delta = \delta(\eta)$ so that $\delta\to 0$ as $\eta\to 0$, with a rate sufficiently rapid so that $K(\eta)\epsilon(\eta,\delta(\eta))\to 0$ as $\eta\to 0$; this is possible due to the first property of $\epsilon(\eta,\delta(\eta))$ mentioned above. Then choose $r=r(\eta)\to 0$ satisfying $r\to 0$ as $\eta\to 0$, with a rate sufficiently slow so that $K^2(\eta)(r^{-2}(\eta)\epsilon(\eta,\delta(\eta)))^{1/2^{k}}\to 0$ as $\eta\to 0$. Since $K(\eta)\to\infty$ as $\eta\to 0$ and $K^2(\eta)(r^{-2}(\eta)\epsilon(\eta,\delta(\eta)))^{1/2^{k}} \to 0$, \eqref{markov} implies $\mathcal{L}(\mathcal{S})\to 0$ as $\eta\to 0$. Introduce another parameter $\zeta = \zeta(\eta)$ satisfying $\zeta\to 0$ as $\eta\to 0$ and $\zeta\geq\eta$, so that $\Omega(\zeta)\subset\Omega(\eta)$. Let $S=S(\eta)$ to be the set of all $\vec{\nu}=(\nu_{\vec{0}}, (\nu_{\alpha}:\alpha\in\{0,1\}^{k};\alpha\not=\vec{0})) \in\Omega(\zeta)^{\{0,1\}^{k}}$ so that $(\nu_{\alpha}: \alpha\in\{0,1\}^{k};\alpha\not=\vec{0}) \in\mathcal{S}$. It's clear from definition and \eqref{markov} that $\mathcal{L}(S)\to 0$ as $\eta\to 0$. By the definition of $S$, if $\vec{\nu}\in \Omega(\zeta)^{\{0,1\}^{k}}\setminus S$ then,
\begin{equation} \label{lwbTktildeF}
    \mathcal{T}_{k}(1_{F_{\nu_{\alpha}}}:\alpha\in\{0,1\}^{k}) \geq\mathcal{T}_{k}(1_{\mathcal{F}_{\nu_{\alpha}}}:\alpha\in\{0,1\}^{k}) - r.
\end{equation}
Recall that $\mathcal{T}_{k}(1_{\mathcal{F}_{\nu_{\alpha}}}:\alpha\in\{0,1\}^{k}) \geq \phi(\zeta)$ if $\nu_{\alpha}\in\Omega(\zeta)$ and $\phi$ is a positive continuous function. Then \eqref{lwbTktildeF} and the general rearrangement inequality imply:
\begin{equation} \label{TkFstar}
    \mathcal{T}_{k}(1_{F^{*}_{\nu_{\alpha}}}:\alpha\in\{0,1\}^{k}) \geq \phi(\zeta) -r.
\end{equation}
\eqref{compsets1side}, with $t=\zeta$, and \eqref{lwbTktildeF} furthermore imply:
\begin{multline} \label{TTstar}
    \mathcal{T}_{k}(1_{F_{\nu_{\alpha}}}:\alpha\in\{0,1\}^{k}) \geq\mathcal{T}_{k}(1_{\mathcal{F}_{\nu_{\alpha}}}:\alpha\in\{0,1\}^{k}) - r \geq \mathcal{T}_{k}(1_{F^{*}_{\nu_{\alpha}}}: \alpha\in\{0,1\}^{k}) - r - C(k)\epsilon(\zeta,\delta). 
\end{multline}
We now further require $\zeta = \zeta(\eta)$ tending to zero as $\eta\to 0$ with a rate sufficiently slow so that $\phi(\zeta)$ also tends to zero slowly and, 
\begin{align}
    \label{cond1} \frac{r+C(k)\epsilon(\eta,\delta)}{\phi(\zeta)-r} &\to 0\\
    \label{cond2} \frac{K(\eta)^2r^{-2}\epsilon(\eta,\delta)}{\phi(\zeta)-r} &\to 0.
\end{align}
\eqref{TkFstar}, \eqref{TTstar} and \eqref{cond1} imply:
\begin{align}
    \nonumber \mathcal{T}_{k}(1_{F_{\nu_{\alpha}}}:\alpha\in\{0,1\}^{k}) &\geq \mathcal{T}_{k}(1_{F^{*}_{\nu_{\alpha}}}:\alpha\in\{0,1\}^{k}) - \frac{r+ C(k)\epsilon(\eta,\delta)}{\phi(\zeta)-r}\cdot (\phi(\zeta)-r) \\
    \label{FFstarcomp} &\geq (1-o(\eta)) \mathcal{T}_{k}(1_{F^{*}_{\nu_{\alpha}}}:\alpha\in\{0,1\}^{k}).
\end{align}
Note that \eqref{FFstarcomp} is precisely \eqref{lwbrearr2} for the case $\vec{\nu}\in\Omega(\zeta)^{\{0,1\}^{k}}\setminus S$. We now investigate the exceptional set $S$. Hence, in addition to the above requirements on $\zeta$, we also enforce $\zeta(\eta)\to 0$ as $\eta\to 0$ with a rate sufficiently slow so that $(\mathcal{L}(S))^{1/2^{k}}/\zeta(\eta) \to 0$ as $\eta\to 0$. This allows us to find, if $\eta$ is sufficiently small, for every $\vec{\nu}\in S\cap \Omega(2\zeta)^{\{0,1\}^{k}}$ two vectors $\vec{\nu}'\not= \vec{\nu}''\in \Omega(\zeta)^{\{0,1\}^{k}} \setminus S$ such that $\nu_{\alpha} - 2(\mathcal{L}(S))^{1/2^{k}} \leq \nu_{\alpha}''\leq \nu_{\alpha}\leq \nu_{\alpha}' \leq \nu_{\alpha} +2(\mathcal{L}(S))^{1/2^{k}}$, for all $\alpha\in\{0,1\}^{k}$. Then for these vectors $\vec{\nu}\in \Omega(2\zeta)^{\{0,1\}^{k}}\subset\Omega(\eta)^{\{0,1\}^{k}}$:
\begin{align} 
    \nonumber \mathcal{T}_{k}(1_{F_{\nu_{\alpha}}}:\alpha\in\{0,1\}^{k}) &\geq \mathcal{T}_{k}(1_{F_{\nu_{\alpha}'}}:\alpha\in\{0,1\}^{k})\\ 
    \nonumber &\geq \mathcal{T}_{k}(1_{\mathcal{F}_{\nu_{\alpha}'}}:\alpha\in\{0,1\}^{k}) - r \geq \mathcal{T}_{k}(1_{\mathcal{F}_{\nu_{\alpha}''}}:\alpha\in\{0,1\}^{k}) - C(\zeta)(\mathcal{L}(S))^{1/2^{k}} - r\\ 
    \label{twovcomp} &\geq \mathcal{T}_{k}(1_{F^{*}_{\nu_{\alpha}''}}:\alpha\in\{0,1\}^{k}) - C(\zeta)(\mathcal{L}(S))^{1/2^{k}}-r-C(k)\epsilon(\eta,\delta).
\end{align}
Since level sets are nested, $\mathcal{T}_{k}(1_{F_{\nu_{\alpha}}}:\alpha\in\{0,1\}^{k})$ is non-increasing in terms of each variable $\nu_{\alpha}$; that explains the first inequality in \eqref{twovcomp}. The second inequality comes from the definition of the set $S$. Since $\Omega(2\zeta)$ is a compact set and $\mathcal{F}$ is a centered Gaussian, $\mathcal{T}_{k}(1_{\mathcal{F}_{\nu_{\alpha}}}:\alpha\in\{0,1\}^{k})$ is Lipschitz continuous with a Lipschitz constant $C(\zeta)$ in each variable $\nu_{\alpha}\in \Omega(2\zeta)$; hence the third inequality follows. Finally the last comes from \eqref{compsets1side}. It remains to show:
\begin{equation} \label{vdprime}
    C(\zeta)(\mathcal{L}(S))^{1/2^{k}}+ r +C(k)\epsilon(\eta,\delta) \leq o(\eta)\mathcal{T}_{k}(1_{F^{*}_{\nu_{\alpha}''}}:\alpha\in\{0,1\}^{k}).
\end{equation}
Indeed, from the definition of $C(\zeta)$ we have $C(\zeta) \leq CK(\zeta)\leq CK(\eta)$, and from the definition of $S$, we have, 
\begin{equation*} C(\zeta)(\mathcal{L}(S))^{1/2^{k}} \leq C(k)K(\eta)(K(\eta)r^{-2}\epsilon(\eta,\delta))^{1/2^{k}}\leq C(k)K^2(\eta)(r^{-2}\epsilon(\eta,\delta))^{1/2^{k}}.\end{equation*}
Hence \eqref{TkFstar}, \eqref{cond1}, \eqref{cond2} and \eqref{twovcomp} then imply \eqref{vdprime}:
\begin{align*} 
    C(\zeta)(\mathcal{L}(S))^{1/2^{k}}+ r +C(k)\epsilon(\eta,\delta) &\leq \frac{C(k)K^2(\eta)(r^{-2}\epsilon(\eta,\delta))^{1/2^{k}}}{\phi(\zeta)-r}\cdot (\phi(\zeta)-r) \\
    & +\frac{r+C(k)\epsilon(\eta,\delta)}{\phi(\zeta)-r}\cdot (\phi(\zeta)-r) \leq o(\eta)\mathcal{T}_{k}(1_{F^{*}_{\nu_{\alpha}''}}:\alpha\in\{0,1\}^{k}).
\end{align*}
Note that $\mathcal{T}_{k}(1_{F^{*}_{\nu_{\alpha}}}:\alpha\in\{0,1\}^{k})$ has a non-increasing property in terms of each variable $\nu_{\alpha}$, similarly to that of $\mathcal{T}_{k}(1_{F_{\nu_{\alpha}}}:\alpha\in\{0,1\}^{k})$. This fact, \eqref{twovcomp} and \eqref{vdprime} together give us the desired conclusion for $\vec{\nu}\in S$:
\begin{equation} \label{conclwbrearr}
    \mathcal{T}_{k}(1_{F_{\nu_{\alpha}}}:\alpha\in\{0,1\}^{k}) \geq (1-o(\eta))\mathcal{T}_{k}(1_{F^{*}_{\nu_{\alpha}}}:\alpha\in\{0,1\}^{k}).
\end{equation}
Finally \eqref{FFstarcomp} and \eqref{conclwbrearr} together give the desired conclusion of \eqref{lwbrearr2}. $\qed$ 

\begin{remark}
It's important for us to analyze the subset $S$, even when its measure is at most $o(\eta)$, in order to prepare for our next discussions. As noted above, the conclusion of Proposition \ref{prop:christ} is only applicable when we have the conclusion of Proposition \ref{prop:lwbrearr} for the diagonal case $\nu_{\alpha} = \nu$, for all $\alpha\in\{0,1\}^{k}$. Hence we can't afford to bypass even a subset of measure zero.  
\end{remark}

\section{Conclusion for one dimension}

\subsection{Preparation} \label{sec:prep}

Let $n=1$. Recall that in the proof of Proposition \ref{prop:lwbrearr}, the parameter $\delta = \delta(\eta)$ satisfies $\delta\to 0$ as $\eta\to 0$. We can choose $\delta(\eta)$ to be a one-to-one function, in which case it allows us to rephrase the findings of the previous two sections as follows:\\

\noindent
Let $\mathcal{F}$ be a centered Gaussian. For every $\delta$ there exists $\eta = \eta(\delta)$ satisfying $\eta\to 0$ as $\delta\to 0$ such that the following occurs. Suppose $f\in L^{p_{k}}(\mathbb{R})$ is a nonnegative $(1-\delta)$ near extremizer with $\|\mathcal{F}-f^{*}\|_{p_{k}} \leq \delta\|f\|_{p_{k}}$. Then for every $s\in\Omega(\eta) = [\eta, \| \mathcal{F}\|_{\infty} - \eta]$, there exists an interval $I_{s}$ such that $\mathcal{L}(I_{s} \Delta F_{s}) = o(\delta) \mathcal{L}(F_{s})$, $\mathcal{L}(\mathcal{F}_{s}\Delta F^{*}_{s}) = o(\delta)$ and consequently $|\mathcal{L}(\mathcal{F}_{s}) - \mathcal{L}(F_{s})| = o(\delta)$. \\

\noindent
As hinted, to further analyze the distribution of $f$, we replace the measurable set $F_{s}$ with one such corresponding interval $I_{s}$. One conclusion of this chapter is that, if $I_{s}$ and $I_{s'}$ are two intervals such that $s, s' \in\Omega(\eta)$ and $c_{s}, c_{s'}$ are centers of $I_{s}, I_{s'}$ respectively, then $c_{s}$ must be close to $c_{s'}$ in an appropriate sense that will be made clear. To this end, we first argue that a selected tuple $(1_{I_{\nu_{\alpha}}}:\alpha\in\{0,1\}^{k})$ nearly achieves equality in the rearrangement inequality:

\begin{lemma}  \label{lem:rearlow}
Let $\delta, \eta$ be as above. Let $\vec{\nu} = (\nu_{\alpha}:\alpha\in\{0,1\}^{k})$ with $\nu_{\alpha}\in \Omega(\eta)$. There exists $\delta_0 > 0$ such that if $\delta \leq \delta_0$ then:
\begin{equation*} 
    \mathcal{T}_{k} (1_{I_{\nu_{\alpha}}}:\alpha\in\{0,1\}^{k}) \geq (1-o(\delta)) \mathcal{T}_{k} (1_{I^{*}_{\nu_{\alpha}}}:\alpha\in\{0,1\}^{k}).
\end{equation*}
\end{lemma}

\begin{proof}
Since $\mathcal{L}(I^{*}_{s} \Delta F^{*}_{s})\leq \mathcal{L}(I_{s} \Delta F_{s})$, $\mathcal{L}(I^{*}_{s} \Delta F^{*}_{s}) = o(\delta)\mathcal{L}(F_{s})$ if $s\in\Omega(\eta)$. Proposition \ref{prop:lwbrearr} and the Gowers product inequality then imply the following three inequalities:
\begin{align*}
    |\mathcal{T}_{k}(1_{F_{\nu_{\alpha}}}:\alpha\in\{0,1\}^{k}) - \mathcal{T}_{k}(1_{F^{*}_{\nu_{\alpha}}}:\alpha\in\{0,1\}^{k}) | &=o(\delta)\prod_{\alpha\in\{0,1\}^{k}}(\mathcal{L}(F_{\nu_{\alpha}}))^{1/p_{k}}\\
    |\mathcal{T}_{k}(1_{F_{\nu_{\alpha}}}:\alpha\in\{0,1\}^{k}) - \mathcal{T}_{k}(1_{I_{\nu_{\alpha}}}:\alpha\in\{0,1\}^{k}) | 
    &= o(\delta)\prod_{\alpha\in\{0,1\}^{k}}(\mathcal{L}(F_{\nu_{\alpha}}))^{1/p_{k}}\\
    |\mathcal{T}_{k}(1_{F^{*}_{\nu_{\alpha}}}:\alpha\in\{0,1\}^{k}) - \mathcal{T}_{k}(1_{I^{*}_{\nu_{\alpha}}}:\alpha\in\{0,1\}^{k}) | &=o(\delta)\prod_{\alpha\in\{0,1\}^{k}}(\mathcal{L}(F_{\nu_{\alpha}}))^{1/p_{k}}
\end{align*}
which then give us:
\begin{equation} \label{rearlow1}
    \mathcal{T}_{k}(1_{I_{\nu_{\alpha}}}:\alpha\in\{0,1\}^{k}) \geq \mathcal{T}_{k}(1_{I^{*}_{\nu_{\alpha}}}:\alpha\in\{0,1\}^{k}) - o(\delta)\prod_{\alpha\in\{0,1\}^{k}}(\mathcal{L}(F_{\nu_{\alpha}}))^{1/p_{k}}.
\end{equation}
Since $\mathcal{L}(I_{\nu_{\alpha}})\geq(1-o(\delta))\mathcal{L}(\mathcal{F}_{\nu_{\alpha}}) \geq (1-o(\delta))\mathcal{L}(\mathcal{F}_{\|\mathcal{F}\|_{\infty}-\eta})\geq\mathcal{L}(I_{\|\mathcal{F}\|_{\infty}-\eta})$, this means $I^{*}_{\|\mathcal{F}\|_{\infty}-\eta}\subset I^{*}_{\nu_{\alpha}}$, for $\nu_{\alpha}\in\Omega(\eta)$ and $\alpha\in\{0,1\}^{k}$. Define:
\begin{equation*} 
    \psi(\eta) = \mathcal{T}_{k}(1_{I^{*}_{\|\mathcal{F}\|_{\infty}-\eta}}, \text{ ... }, 1_{I^{*}_{\|\mathcal{F}\|_{\infty}-\eta}}).
\end{equation*} 
Then the said set inclusion implies:
\begin{equation} \label{reartail}
    \mathcal{T}_{k}(1_{I^{*}_{\nu_{\alpha}}}:\alpha\in\{0,1\}^{k}) \geq\mathcal{T}_{k}(1_{I^{*}_{\|\mathcal{F}\|_{\infty}-\eta}}, \text{ ... }, 1_{I^{*}_{\|\mathcal{F}\|_{\infty}-\eta}}) = \psi(\eta).
\end{equation}
Note that $\mathcal{L}(\mathcal{F}_{\eta}) = C(\log (1/\eta))^{C}$ and that $\psi(\eta)$ stays strictly bounded below as long as $\eta > 0$. We select $\eta = \eta(\delta)\to 0$ sufficiently slowly as $\delta\to 0$ and $\delta_0 > 0$ so that if $\delta \leq \delta_0$ then $o(\delta) (\log (1/\eta(\delta)))^{C} = o(\delta)$. Then \eqref{reartail} gives us:
\begin{align} 
    \nonumber -o(\delta)\prod_{\alpha\in\{0,1\}^{k}}( \mathcal{L}(F_{\nu_{\alpha}}))^{1/p_{k}} &\geq -o(\delta)(\mathcal{L}(\mathcal{F}_{\eta}))^{k+1} \geq -o(\delta)(\log (1/\eta))^{C} \geq -o(\delta)\psi(\eta) \\ 
    \label{rearlow2} &= -o(\delta)\mathcal{T}_{k}(1_{I^{*}_{\|\mathcal{F}\|_{\infty}-\eta}}, \text{ ... }, 1_{I^{*}_{\|\mathcal{F}\|_{\infty}-\eta}}) \geq -o(\delta)\mathcal{T}_{k}(1_{I^{*}_{\nu_{\alpha}}}: \alpha\in\{0,1\}^{k}).
\end{align}
A combination of \eqref{rearlow1} and \eqref{rearlow2} then give the desired conclusion.    
\end{proof}

\noindent
We will show in Subsection \ref{sec:alignment} below that Lemma \ref{lem:rearlow} implies the centers $c_{\nu_{\alpha}}, c_{\nu_{\beta}}$ are close to each other, for $\alpha,\beta\in\{0,1\}^{k}$ and $\nu_{\alpha}, \nu_{\beta}\in \Omega(\eta)$. We now replace a near extremizer with a superposition of interval approximations of its super-level sets. Suppose there exist a sequence $\delta_{i}\to 0$ as $i\to\infty$, a nonnegative sequence of functions $\{f_{i}\}_{i}$ and a centered Gaussian extremizer $\mathcal{F}$, such that $\| f_{i} \|_{U^{k}} \geq (1-\delta_{i})A(k)\|f_{i}\|_{p_{k}}$ and $\| \mathcal{F}-f^{*}_{i}\|_{p_{k}} \leq \delta_{i}\|f_{i}\|_{p_{k}}$. Assume $\|f_{i}\|_{p_{k}} = 1$ for all $i$. Let $s\in\mathbb{R}_{>0}$. Denote $F_{i, s} = \{ f_{i} > s\}$ and $F^{*}_{i,s} = \{ f^{*}_{i} > s\}$. Consider another sequence $\eta_{i}\to 0$ as $i\to\infty$ and suppose further that, if $s\in [\eta_{i}, \| \mathcal{F}\|_{\infty} - \eta_{i}]$, then there exists an interval $I_{i,s}$ such that $\mathcal{L}(I_{i,s} \Delta F_{i, s}) \leq \delta_{i} \mathcal{L}(F_{i,s})$. Let $h_{i}(x) =\int_{\eta_{i}}^{\| \mathcal{F}\|_{\infty} - \eta_{i}} 1_{I_{i,s}}(x) \, ds$.

\begin{lemma} \label{lem:2funcomp}
Let $f_{i}$ and $h_{i}$ be as above. Then $\| f_{i} - h_{i}\|_{p_{k}} \to 0$ as $i\to\infty$.     
\end{lemma}

\begin{proof}
By Lemma \ref{lem:levelsetstack} and Minkowski's integral inequality,
\begin{align*}
    \| f_{i} - h_{i} \|_{p_{k}} &\leq \bigg \| \int_{\eta_{i}}^{\| \mathcal{F}\|_{\infty} - \eta_{i}} 1_{F_{i,s}} \, ds - \int_{\eta_{i}}^{\| \mathcal{F}\|_{\infty} - \eta_{i}} 1_{I_{i,s}} \, ds \bigg \|_{p_{k}} + C\eta_{i}(\log (1/\eta_{i}))^{C} \\
    &\leq \int_{\eta_{i}}^{\|\mathcal{F}\|_{\infty}-\eta_{i}} \| 1_{F_{i,s}} - 1_{I_{i,s}}\|_{p_{k}} \, ds + C\eta_{i}(\log (1/\eta_{i}))^{C}\\
    &= \int_{\eta_{i}}^{\|\mathcal{F}\|_{\infty}-\eta_{i}} (\mathcal{L}(I_{i,s} \Delta F_{i,s})^{1/p_{k}} \, ds + C\eta_{i}(\log (1/\eta_{i}))^{C} \\
    &\leq \delta_{i} \int_{\eta_{i}}^{\|\mathcal{F}\|_{\infty}-\eta_{i}} (\mathcal{L}(F_{i,s}))^{1/p_{k}} \, ds + C\eta_{i}(\log (1/\eta_{i}))^{C} \leq \delta_{i} + C\eta_{i}(\log (1/\eta_{i}))^{C}.
\end{align*}
The last inequality is due to $\|f_{i}\|_{p_{k}} = 1$. Let $i\to\infty$, we get the conclusion.     
\end{proof}

\noindent
By Lemma \ref{lem:2funcomp}, if we can establish that $\{ h_{i} \}_{i}$ is precompact in $L^{p_{k}}$, then we have the same result for $\{f_{i}\}_{i}$. We show here a related compactness result which will be needed. 

\begin{lemma} \label{lem:precompact}
Let $a < b \in \mathbb{R}$ and $1\leq q < \infty$. Let $[-B,B]$ be a closed interval. For each $i\in\mathbb{Z}_{>0}$ and $s\in\mathbb{R}$ let $I_{i,s} \subset [-B,B]$ be an interval. Suppose further that the function $(x,s) \mapsto 1_{I_{i,s}}(x)$ is measurable. Then $\{ \int_{a}^{b} 1_{I_{i,s}}(x) \, ds \}_{i}$ is precompact in $L^{q}$. 
\end{lemma}

\begin{proof}
Let $g_{i}(x) = \int_{a}^{b} 1_{I_{i,s}}(x) \, ds$. It's clear that, there exists $C>0$ such that $\|g_{i}\|_{q}\leq C$. Moreover, $\int_{\{|x|> r\}} |g_{i}(x)|^{q} \, dx = 0$, for every $r > B$, and $\lim_{h\to 0} \| T^{h}g_{i} - g_{i} \|_{q} = 0$ uniformly in $i$. Hence by the Fr\'echet-Kolmogorov theorem \cite{brezis2011functional}, $\{ g_{i}\}_{i}$ is precompact in $L^{q}$.     
\end{proof}

\subsection{A monotonicity result} \label{sec:mono}

Let $I = [-1,1]$ and $J = [-\eta -1, 1+\eta]$ for some $\eta\in [0,\frac{2}{k-1}]$. Define:
\begin{equation*} 
    \phi(t) = \int_{\mathbb{R}^{k+1}} 1_{J+t}(x) \prod_{\alpha\in\{0,1\}^{k};\alpha\not=\vec{0}} 1_{I}(x+\alpha\cdot\vec{h}) \, d\vec{h}dx.
\end{equation*}
$\phi(t)$ is a continuous, nonnegative even function of $t$ and has a compact support. Furthermore, let:
\begin{equation*} 
    H(x) = \int_{\mathbb{R}^{k}} \prod_{\alpha \in \{0,1\}^{k};\alpha \not= \vec{0}} 1_{I}(x+ \alpha \cdot \vec{h}) \, d\vec{h}.
\end{equation*}
Then $\phi(t) = \int_{\mathbb{R}} 1_{J+t}\cdot H(x)\, dx$. $H$ is also a continuous, nonnegative even function whose support is the interval $[-\frac{k+1}{k-1},\frac{k+1}{k-1}]$. Indeed, the interval $[-1,1]$ is clearly contained in the support of $H$. Suppose $x > 1$ and $x\in spt(H)$. Then there exists $\vec{h}=(h_{i})_{i}\in\mathbb{R}^{k}$ such that $|x+\alpha\cdot\vec{h}|\leq 1$ for all $\alpha\in\{0,1\}^{k}$ (there exists, in fact, a set of such $\vec{h}$ of positive measure). Let $j\in\{1, \text{ ... },k\}$ and define $\beta\in\{0,1\}^{k}$ by $\beta_{j}=1$ and $\beta_{i}=0$ if $i\not=j$. Then $|x+\beta\cdot\vec{h}| = |x+h_{j}| \leq 1$ implies $h_{j}\leq 1-x$. Since this holds for every $j\in\{1,\text{ ... }, k\}$, $k^{-1}\sum_{i=1}^{k} h_{i}\leq 1-x$. Now consider $\alpha = \vec{1} = (1,\text{ ... }, 1)$. Then $|x+\vec{1}\cdot\vec{h}| = |x+\sum_{i=1}^{k} h_{i}| \leq 1$ implies $k^{-1}\sum_{i=1}^{k}h_{i}\geq -k^{-1}(1+x)$. By transitivity, $-k^{-1}(1+x) \leq 1-x$, or equivalently $x\leq\frac{k+1}{k-1}$. Similarly, if $x < -1$ and $x\in spt(H)$ then $x \geq -\frac{k+1}{k-1}$. \\
Now suppose $x\in [-\frac{k+1}{k-1},\frac{k+1}{k-1}]$. We define $\vec{h} = (h_{i})_{i}\in\mathbb{R}^{k}$ such that $h_{i} \leq \frac{1-x}{k}$ if $1 < x \leq \frac{k+1}{k-1}$ and $h_{i} \leq \frac{x-1}{k}$ if $-\frac{k+1}{k-1} \leq x < -1$. For the former case, $-1 \leq x(1-k) + k \leq x+\alpha\cdot\vec{h} \leq x + 1-x = 1$, for all $\alpha\in\{0,1\}^{k}$; hence $x\in spt(H)$. We also obtain the same conclusion for the latter case. \\

\noindent
We claim that $\phi(t)$ is a strictly decreasing function of $t\geq 0$ in its support. We make a quick remark that this claim and the fact that $spt(H) = [-\frac{k+1}{k-1},\frac{k+1}{k-1}]$ are why we take $\eta\in [0,\frac{2}{k-1}]$. Since if $\eta > \frac{2}{k-1}$ then there exist $t_1 > t_0$ sufficiently small, so that they are both in the support of $\phi$, say, $t_1 = (\eta-\frac{2}{k-1})/10 > t_0 = (\eta - \frac{2}{k-1})/20$, such that, 
\begin{multline*} 
    \phi(t_1) = \int_{\mathbb{R}^{k+1}} 1_{J+t_1}(x) \prod_{\alpha\in\{0,1\}^{k};\alpha\not= 0} 1_{I}(x+\alpha\cdot\vec{h}) \, d\vec{h} dx \\ = \int_{\mathbb{R}^{k+1}} 1_{J+t_0}(x) \prod_{\alpha\in\{0,1\}^{k};\alpha\not= 0} 1_{I}(x+\alpha\cdot\vec{h}) \, d\vec{h} dx = \phi(t_0).
\end{multline*}
With our choice of $\eta$, $spt(\phi) \subset [-\frac{2(k+1)}{k-1},\frac{2(k+1)}{k-1}]$. Note that:
\begin{equation*} 
    (d/dt) \phi (t) = \int_{\mathbb{R}} (\delta_{t+1+ \eta} - \delta_{t-1-\eta})\cdot H (x)\, dx = H(t+1+\eta) - H(t-1-\eta).
\end{equation*} 
Furthermore, with our choice of $\eta$, $t-1-\eta$ always lies in the support of $H$. Hence our claim on the monotonicity of $\phi$ will follow if we can show $H(x)$ is strictly decreasing for $x \geq 0$ in its support. Let $\mathcal{C} = \{(x,\vec{h})\in\mathbb{R}^{k+1}: x+\alpha\cdot\vec{h}\in I, \forall \alpha\in\{0,1\}^{k}\}$. Then $\mathcal{C}$ is a compact, convex, balanced subset of $\mathbb{R}^{k+1}$. Let $\mathcal{C}_{x} = \{ \vec{h}\in\mathbb{R}^{k}: (x,\vec{h})\in\mathcal{C}\}$. We observe that $H(x) = \mathcal{L}(\mathcal{C}_{x})$. 

\begin{lemma} \label{lem:H}
Let $H$ be as above. Then $H(0) > H(x)$, for all $x > 0$.  
\end{lemma}

\noindent
Lemma \ref{lem:H}, if assumed true, will imply that $H(x) > H(y)$ if $0\leq x < y$ and $x, y$ both lie in the support of $K$. Indeed, let $0 < x < y$ be in the interior of the support of $H$, then $H(x) = \mathcal{L}(C_{x}) > 0$ and $H(y) = \mathcal{L}(C_{y}) > 0$. Let $t\in (0,1)$ be such that $x = (1-t)\cdot 0 + ty$. By the convexity of $\mathcal{C}$, $C_{x} \supset (1-t) C_0 + tC_{y}$. Then by the Brunn-Minkowski inequality, $\mathcal{L}(C_{x}) \geq (\mathcal{L}(C_0))^{1-t}(\mathcal{L}(C_{y}))^{t}$, and by Lemma \ref{lem:H}:
\begin{equation*} 
    H(x) \geq H(0)^{1-t}H(y)^{t} > H(y)^{1-t}H(y)^{t} = H(y).
\end{equation*}

\subsection{Proof of Lemma \ref{lem:H}} \label{sec:lemHpf}

We will deduce Lemma \ref{lem:H} from the following more general result. The setup is as follows:\\

\noindent
Let $M, N\in\mathbb{Z}_{>0}$. For $i\in\{1,\text{ ... }, M\}$, let $J_{i}$ be a closed interval centered at $0$, $\mathcal{L}(J_{i}) = l_{i} \in\mathbb{R}_{>0}$, and $L_{i}:\mathbb{R}^{N}\to\mathbb{R}$ be surjective linear mappings. For $\vec{t} = (t_{i})_{i} \in \mathbb{R}^{M}$, we let:
\begin{equation*} 
    \Psi(\vec{t}) = \int_{\mathbb{R}^{N}} \prod_{i = 1}^{M} 1_{J_{i} + t_{i}} (L_{i}(\vec{y})) \, d\vec{y}.
\end{equation*}
We make two observations. First, if $\vec{v} \in\mathbb{R}^{N}$:
\begin{equation*} 
    \int_{\mathbb{R}^{N}}\prod_{i=1}^{M} 1_{J_{i}}(L_{i}(\vec{y}+\vec{v})) \, d\vec{y} = \int_{\mathbb{R}^{N}}\prod_{i=1}^{M} 1_{J_{i}}(L_{i}(\vec{y})) \, d\vec{y}.
\end{equation*} 
Secondly, by the general rearrangement inequality, for $\vec{t} = (t_{i})_{i}\in\mathbb{R}^{M}$:
\begin{equation*} 
    \int_{\mathbb{R}^{N}}\prod_{i=1}^{M} 1_{J_{i}+t_{i}}(L_{i}(\vec{y})) \, d\vec{y} \leq \int_{\mathbb{R}^{N}}\prod_{i=1}^{M} 1_{J_{i}}(L_{i}(\vec{y})) \, d\vec{y}
\end{equation*} 
or equivalently, $\Psi(\vec{t}) \leq \Psi(\vec{0})$. It then follows that, for $\vec{t}\in\mathbb{R}^{M}$, $\Psi(\vec{t}) = \Psi(\vec{0})$ if there exists $\vec{v} \in\mathbb{R}^{N}$ such that $L_{i}(\vec{v}) = t_{i}$, $i\in\{1, \text{ ... }, M\}$. It will be shown that this is also a necessary condition, provided that $(L_{i},l_{i})_{i=1}^{M}$ is \textbf{an admissible tuple}:\\

\noindent
\emph{Admissibility:} For $i\in\{1,\text{ ... }, M\}$, let $J_{i}$ be a closed interval centered at $0$, $\mathcal{L}(J_{i}) = l_{i} \in\mathbb{R}_{>0}$ and $L_{i}:\mathbb{R}^{N}\to\mathbb{R}$ be surjective linear mappings. Let $\vec{l} = (l_{i})_{i}\in(\mathbb{R}_{>0})^{M}$. Define $\mathcal{K}_{\vec{l}} = \{ \vec{x} \in\mathbb{R}^{N}: |L_{i}(\vec{x})| \leq l_{i}, \forall i \in\{1, \text{ ... }, M\} \}$. We said that $(L_{i},l_{i})_{i}$ is \textbf{an admissible tuple} if for every $m\in\{1, \text{ ... }, M\}$ there exists $\vec{x}_{m}\in\mathcal{K}_{\vec{l}}$ such that $|L_{m}(\vec{x}_{m})|= l_{m}$.

\begin{lemma} \label{lem:personal} \cite{personalcommunication}
Let $L_{i}, J_{i}, l_{i}$ and $\Psi$ be as above. Suppose $(L_{i},l_{i})_{i=1}^{M}$ is an admissible tuple. Let $\vec{t}\in\mathbb{R}^{M}$. Then $\Psi(\vec{0})\geq\Psi(\vec{t})$ and equality holds iff there exists $\vec{v}\in\mathbb{R}^{N}$ such that $L_{i}(\vec{v}) = t_{i}$ for all $i\in\{1, \text{ ... }, M\}$.  
\end{lemma}

\begin{proof}
Let $K = \{ (\vec{x},\vec{t})\in\mathbb{R}^{N}\times\mathbb{R}^{M}: L_{i}(\vec{x}) \in J_{i}+t_{i}, \forall i \in \{1,\text{ ... }, M\}\}$. For each $\vec{t}\in\mathbb{R}^{M}$, let $K(\vec{t}) = \{ \vec{x}\in\mathbb{R}^{N}: (\vec{x},\vec{t}) \in K\}$. It's clear that $K$ is convex, and if $(\vec{x},\vec{t})\in K$ then $|L_{i}(\vec{x}) - t_{i}| \leq l_{i}$. Since $\Psi(\vec{t})$ represents the $N$-dimensional volume of $K(\vec{t})$ and $J_{i}$ are centered at $0$, $\Psi(\vec{t}) = \Psi(-\vec{t})$ or $\mathcal{L}(K(\vec{t})) = \mathcal{L}(K(-\vec{t}))$. Moreover, 
\begin{equation} \label{setin}
    K(\vec{0}) \supset (1/2)K(\vec{t})+(1/2)K(-\vec{t}). 
\end{equation}
Suppose that $\Psi(\vec{t}) = \Psi(\vec{0})$ for some $\vec{t}\in\mathbb{R}^{M}$, which implies $\Psi(\vec{0}) = \mathcal{L}(K(\vec{0})) = \Psi(\vec{t}) = (\mathcal{L}(K(\vec{t})))^{1/2}(\mathcal{L}(K(-\vec{t})))^{1/2}$. Then \eqref{setin} and the Brunn-Minkowski inequality imply:
\begin{equation*} 
    \mathcal{L}(K(\vec{0}))\geq \mathcal{L}((1/2)K(\vec{t})+(1/2)K(-\vec{t})) \geq (\mathcal{L}(K(\vec{t})))^{1/2}(\mathcal{L}(K(-\vec{t})))^{1/2} = \mathcal{L}(K(\vec{0}))
\end{equation*}
which yields $\mathcal{L}((1/2)K(\vec{t})+(1/2)K(-\vec{t})) = (\mathcal{L}(K(\vec{t})))^{1/2}(\mathcal{L}(K(-\vec{t})))^{1/2}$. Hence by the characterization of the equality case of the Brunn-Minkowski inequality, there exists $\vec{v}\in\mathbb{R}^{N}$ such that $K(\vec{t}) + \vec{v} = K(\vec{0})$. By definition, this means, if $\vec{x}\in\mathbb{R}^{N}$ and $i\in \{1, \text{ ... }, M\}$:
\begin{equation} \label{trans}
    L_{i}(\vec{x})\in J_{i} + t_{i} \iff L_{i}(\vec{x}+\vec{v}) \in J_{i}.
\end{equation}
Let $\vec{z} = \vec{x}+\vec{v}$, then by \eqref{trans}, $|L_{i}(\vec{z})|\leq l_{i}$ implies,
\begin{equation} \label{nearlin1}
    |L_{i}(\vec{z}) - t_{i} - L_{i}(\vec{v})| \leq l_{i}.
\end{equation}
Fix $m\in\{1,\text{ ... },M\}$. Admissibility assumption implies that there exist $\vec{x}_{m,\pm}\in\mathbb{R}^{N}$ such that $|L_{i}(\vec{x}_{m,\pm})| \leq l_{i}$ and $L_{m}(\vec{x}_{m,\pm}) = \pm l_{m}$. In particular, \eqref{nearlin1} implies 
\begin{equation} \label{nearlin2}
    |L_{m}(\vec{x}_{m,\pm}) - t_{m} - L_{m}(\vec{v})| \leq l_{m}. 
\end{equation}
Suppose $t_{m} + L_{m}(\vec{v}) < 0$. Then $L_{m}(\vec{x}_{m,+}) - t_{m} - L_{m}(\vec{v}) = l_{m} - (t_{m} + L_{m}(\vec{v})) > l_{m}$, which poses a contradiction to \eqref{nearlin2}. Similarly, $t_{m} + L_{m}(\vec{v}) > 0$ will also imply the contradiction to \eqref{nearlin2} since $L_{m}(\vec{x}_{m,-}) - (t_{m} + L_{m}(\vec{v})) = -l_{m} - (t_{m} + L_{m}(\vec{v})) < -l_{m}$. Hence $t_{m} = L_{m}(-\vec{v})$ for every $m\in\{1,\text{ ... }, M\}$. We've completed the proof of the lemma as the "if" direction is apparent by the discussion at the beginning of this section.    
\end{proof}

\noindent
To finish the proof of Lemma \ref{lem:H}, let $M =2^{k}-1$, $N = k$. We note that for each $\vec{0}\not=\alpha\in\{0,1\}^{k}$, $L_{\alpha}:\mathbb{R}^{k}\to\mathbb{R}$ defined by $L_{\alpha}(\vec{h}) = \alpha\cdot\vec{h}$ is a surjective linear mapping. Moreover, for each $\vec{0}\not=\alpha\in\{0,1\}^{k}$, let $|\alpha| = \sum_{i=1}^{k} \alpha_{i}$ and $\vec{h}_{\alpha}\in\mathbb{R}^{k}$ be such that $(\vec{h}_{\alpha})_{i} = 1/|\alpha|$ if $\alpha_{i} = 1$ and $(\vec{h}_{\alpha})_{i} = 0$ otherwise. Then it's easily checked that $L_{\alpha}(\pm\vec{h}_{\alpha}) =\pm 1$ and $|L_{\beta}(\vec{h}_{\alpha})|\leq 1$ if $\alpha\not=\beta\in\{0,1\}^{k}\setminus\{\vec{0}\}$. Hence the tuple $((L_{\alpha},1)_{\alpha}:\alpha\in\{0,1\}^{k}\setminus\{\vec{0}\})$ is admissible. Note that with $\vec{x} = (x)_{i}\in\mathbb{R}^{\{0,1\}^{k}\setminus\{\vec{0}\}}$:
\begin{equation*} H(x) = \int_{\mathbb{R}^{k}} \prod_{\alpha \in \{0,1\}^{k};\alpha \not= \vec{0}} 1_{I}(x+ \alpha \cdot \vec{h}) \, d\vec{h} = \int_{\mathbb{R}^{k}} \prod_{\alpha \in \{0,1\}^{k};\alpha \not= \vec{0}} 1_{I-x}(\alpha \cdot \vec{h}) \, d\vec{h} = \Psi(-\vec{x}).\end{equation*}
Then from Lemma \ref{lem:personal}, $H(0) = \Psi(\vec{0}) \geq \Psi(-\vec{x}) = H(x)$. If $H(0) = H(x)$ for some $x > 0$ then there must exist $\vec{v}\in\mathbb{R}^{k}$ such that $L_{\alpha}(\vec{v}) = \alpha\cdot\vec{v} = -x$, for all $\vec{0}\not=\alpha\in\{0,1\}^{k}$, which is clearly impossible. Hence $H(0) > H(x)$ for all $x > 0$, as we wish to conclude. $\qed$

\subsection{Alignment of super-level sets} \label{sec:alignment}

Now Subsection \ref{sec:mono} concludes that, if $I$ and $J$ are two intervals "compatible" in size, and if 
\begin{equation*} 
    \mathcal{T}_{k} (J, I, \text{ ... }, I)= \int_{\mathbb{R}^{k+1}} 1_{i}(x) \prod_{\alpha \in \{ 0,1\}^{k};\alpha \not= \vec{0}} 1_{I}(x+\alpha \cdot \vec{h}) \, d\vec{h}dx
\end{equation*} 
is nearly maximized over all tuples of intervals of the same sizes, then the centers of $I$ and $J$ must be "close" to each other. The compatibility condition is given by allowing $\mathcal{L}(J)=\mathcal{L}(I)+\eta$ with $\eta\in [0, \frac{2}{k-1}\mathcal{L}(I)]$.\\

\noindent
Let $0<\epsilon\ll \frac{2}{k-1}$. Then the discussion in the beginning of Subsection \ref{sec:prep} and Lemma \ref{lem:rearlow} conclude that there exist $\delta > 0$ and $\eta = \eta(\delta) > 0$ satisfying the following properties:\\
\begin{enumerate}
    \item If $f$ is a nonnegative $(1-\delta)$ near extremizer and $\|\mathcal{F}-f\|_{p_{k}}\leq\delta$ for a centered Gaussian $\mathcal{F}$, then for every $s\in \Omega(\eta) = [\eta, \| \mathcal{F}\|_{\infty} - \eta]$ there exists an interval $I_{s}$ such that $\mathcal{L}(I_{s} \Delta F_{s}) < c\epsilon\mathcal{L}(F_{s})$ and $|\mathcal{L}(\mathcal{F}_{s}) - \mathcal{L}(F_{s})| < c\epsilon$. 
    \item Furthermore, $\mathcal{T}_{k}(1_{I_{\nu_{\alpha}}}:\alpha\in\{0,1\}^{k}) \geq (1-\epsilon)\mathcal{T}_{k}(1_{I^{*}_{\nu_{\alpha}}}:\alpha\in\{0,1\}^{k})$, if $\nu_{\alpha}\in\Omega(\eta)$, $\alpha\in\{0,1\}^{k}$.
\end{enumerate}
Let $\mathcal{F}_{s}$ denote the super-level set of the Gaussian $\mathcal{F}$ associated with the value $s$. There exist $s_1 = \eta < s_2 < \text{ ... } < s_{N} = \| \mathcal{F}\|_{\infty} - \eta \in [\eta, \|\mathcal{F}\|_{\infty}-\eta]$ with $N=N(\eta)$ such that $\mathcal{L}(\mathcal{F}_{s_{i}} \Delta \mathcal{F}_{s_{i+1}}) \asymp \frac{1}{k-1} \mathcal{L}(\mathcal{F}_{s_{i}})$ for all $i\in\{1,\text{ ... }, N-1\}$. Since super-level sets are nested,  $\mathcal{L}(\mathcal{F}_{s} \Delta \mathcal{F}_{s_{i}}) \leq\frac{c}{k-1}\mathcal{L}(\mathcal{F}_{s_{i}})$, $c$ being some sufficiently small constant, if $s_{i}\leq s\leq s_{i+1}$, and consequently, $|\mathcal{L}(I_{s}) - \mathcal{L}(I_{s_{i}})| \leq \frac{c}{k-1}\mathcal{L}(I_{s_{i}})$, $i\in\{1,\text{ ... }, N-1\}$. That is to say, the size compatibility condition is also satisfied by the sub-intervals with $s\in [s_{i},s_{i+1}]\subset\Omega(\eta)$. Recall that $c_{s}$ denotes the center of an interval $I_{s}$. If $\epsilon$ is sufficiently small, then by the second property, the previous paragraph and the satisfaction of the compatibility condition, $|c_{s}-c_{s_{i}}| =o(\epsilon)\mathcal{L}(I_{s_{i}})$ if $s_{i}\leq s \leq s_{i+1}$, $i\in\{1,\text{ ... }, N-1\}$. If $s\in\Omega(\eta)$, then $s$ must lie in one such interval $[s_{i}, s_{i+1}]$; hence $|c_{s}-c_{\eta}|=o(\epsilon)\mathcal{L}(I_{\eta}) = o(\epsilon)(\log(1/\eta))^{C}$ for $s\in\Omega(\eta)$ since $\mathcal{L}(\mathcal{F}_{\eta}) = C(\log (1/\eta))^{C}$. Note that if $\epsilon\to 0$ then $\delta\to 0$, and if $\eta' > \eta$ then $\Omega(\eta)\subset\Omega(\eta')$. We now require $\eta=\eta(\delta)\to 0$ sufficiently slow so that $|c_{s}-c_{\eta}|=o(\epsilon)\mathcal{L}(I_{\eta}) = o(\epsilon)(\log (1/\eta))^{C} = o(\delta)$ as $\delta\to 0$. The intervals $I_{s}$ might change as the parameters change, but the size estimates still hold. We now obtain the following lemma:\\

\begin{lemma}
Let $\delta, \eta=\eta(\delta)$ and $c_{s}$ be as above. Then $|c_{s} - c_{\eta}| = o(\delta)$ if $s\in\Omega(\eta)$.    
\end{lemma}

\subsection{A compactness result}

Suppose we have a nonnegative extremizing sequence $\{f_{i}\}_{i}$ such that $\|f_{i}\|_{p_{k}}=1$ and
\begin{equation*} 
    \| f_{i}\|_{U^{k}} \geq (1-\delta_{i})A(k)\|f_{i}\|_{p_{k}} = (1-\delta_{i})A(k)
\end{equation*}
and that $\|\mathcal{F}- f^{*}_{i}\|_{p_{k}} \leq \delta_{i}$ for some centered Gaussian $\mathcal{F}$ with $\delta_{i}\to 0$ as $i\to\infty$. Recall that in Subsection \ref{sec:prep} we defined for such a sequence $\{f_{i}\}_{i}$ an associate superposition sequence $h_{i}(x) =\int_{\eta_{i}}^{\| \mathcal{F}\|_{\infty} - \eta_{i}} 1_{I_{i,s}}(x) \, ds$. We select $\eta_{i} = \eta(\delta_{i})$ satisfying $\eta_{i}\to 0$ as $i\to\infty$ such that $|c_{s} - c_{\eta_{i}}| = o(\delta_{i})$ if $s\in\Omega(\eta_{i})$, as in Subsection \ref{sec:alignment}. The conclusions of Lemma \ref{lem:2funcomp} and Lemma \ref{lem:precompact} stay unchanged with this selection of $\{ \eta_{i}\}_{i}$. 

\begin{remark} 
It's possible to select $I_{i,s_{i}}$ so that $(x,s_{i})\mapsto 1_{I_{i,s_{i}}}(x)$ is a measurable function. Fix $i$ and denote $\eta_{i}$ as $\eta$, $f_{i}$ as $f$, $F_{i,s_{i}}$ as $F_{s}$ and $I_{i,s_{i}}$ as $I_{s}$. The set $E = \cup_{s\in\Omega(\eta)} \{s\}\times F_{s}$ is a measurable subset of $\Omega(\eta)\times\mathbb{R}$. Let $a_{s} < b_{s}$ denote the endpoints of $I_{s}$, $s\in\Omega(\eta)$, then it comes down to the ability to select the endpoints $a_{s}, b_{s}$ of the intervals $I_{s}$ in a measurable manner. We suppose for a moment that $f$ is a continuous function. Decompose the range of values $\Omega(\eta)$ into a finite number $C_{f}(\eta)$ of smaller sub-ranges of values $\Omega'(\eta)$ such that $\mathcal{L}(F_{s}\Delta F_{s'}) = o(\eta)$ for each pair $s, s'\in\Omega'(\eta)$. Fix one such sub-range $\Omega'(\eta)$ and let $F_{s_{t}}, F_{s_{b}}$ denote the super-level sets of $f$ associated with the largest and smallest values of the range, respectively. By replacing $F_{s_{t}}$ with $I_{s_{t}}$ and $F_{s_{b}}$ with $I_{s_{b}}$, we assume $F_{s_{t}}$ and $F_{s_{b}}$ are both intervals. Let $a_{s_{t}}$ denote the left endpoint of $F_{s_{t}}$ and $a_{s_{b}}$ the left endpoint of $F_{s_{b}}$; note that $|a_{s_{t}}-a_{s_{b}}|=C_{f}o(\eta)$. Consider the part of the graph of $f$ inside the rectangle box $\Omega'(\eta)\times [a_{s_{t}},a_{s_{b}}]$; call this set $R_{\eta}$. By the Measurable Choice Theorem \cite{bogachev2007measure}, there exists a measurable function $s\mapsto a(s)$, so that $(s,a(s))\in R_{\eta}$, for every $s\in\Omega'(\eta)$. Take $a(s)$ to be the left endpoint of our interval $I_{s}$. Proceed similarly to obtain a measurable function $s\mapsto b(s)$ for the right endpoint of $I_{s}$. Note that by construction $\mathcal{L}(I_{s}\Delta F_{s}) = o(\eta)\mathcal{L}(F_{s})$, if $s\in\Omega'(\eta)$. Continue this procedure for each of these sub-ranges $\Omega'(\eta)$ and concatenate the obtained left endpoint and right endpoint functions to obtain a measurabe function $(x,s)\mapsto 1_{I_{s}}(x)$, $s\in\Omega(\eta)$. For the general case we approximate $f$ with a positive continuous function $g$ so that $\|f-g\|_{p_{k}} = \rho(\eta)$ with $\rho\ll\eta$ in order for us to have $\mathcal{L}(G_{s}\Delta F_{s})=o(\eta)$, for a.e $s\in\Omega(\eta)$. Then apply the described procedure with $G_{s}$ in place of $F_{s}$. \\
Another and easier way is to construct a piecewise constant function $(x,s)\mapsto I_{s}(x)$; ie, $I_{s}(x) = I_{s'}(x)$ if $s,s'\in\Omega'(\eta)$. This ensures measurability but might come at the expense of increasing the (still finite) number of sub-ranges so that $\mathcal{L}(I_{s}\Delta F_{s}) = o(\eta)\mathcal{L}(F_{s})$ is still guaranteed.   
\end{remark}

\begin{proposition} \label{prop:shiftcomp}
Let $\{h_{i}\}_{i}$ be as above. There exists $\{ a_{i}\}_{i}$ such that $\{ h_{i}(\cdot - a_{i})\}_{i}$ is precompact in $L^{p_{k}}$.    
\end{proposition}
    
\begin{proof}
Let $a_{i} = c_{\eta_{i}}$ with $c_{\eta_{i}}$ being the center of the interval $I_{\eta_{i}}$. Define:
\begin{equation*} 
    g_{i}(x) = h_{i}(x - a_{i}) = \int_{\eta_{i}}^{\|\mathcal{F}\|_{\infty} - \eta_{i}} 1_{I_{i, s}}(x-a_{i})\, ds.
\end{equation*} 
Let $\epsilon > 0$ be a small number and select $\eta>0$ such that:
\begin{equation} \label{etatail}
    \| \int_0^{\eta} 1_{\mathcal{F}_{s}} \, ds \|_{p_{k}} + \| \int_{\|\mathcal{F}\|_{\infty}-\eta}^{\|\mathcal{F}\|_{\infty}} 1_{\mathcal{F}_{s}} \, ds \|_{p_{k}} \leq \epsilon.
\end{equation}
For simplicity, we assume that $\eta_{i} \leq \eta$ for all $i$. Recall that $|\mathcal{L}(\mathcal{F}_{s}) - \mathcal{L}(I_{i,s})| = o(\delta_{i})\mathcal{L}(\mathcal{F}_{s})$ if $s\in \Omega(\eta_{i})$. Then \eqref{etatail} entails:
\begin{equation} \label{etatailI}
    \| \int_{\eta_{i}}^{\eta} 1_{I_{i,s}}(\cdot - a_{i}) \, ds \|_{p_{k}} + \| \int_{\|\mathcal{F}\|_{\infty} - \eta_{i}}^{\|\mathcal{F}\|_{\infty} - \eta} 1_{I_{i,s}}(\cdot - a_{i}) \, ds \|_{p_{k}} \leq C\epsilon.
\end{equation}
The fact that $|\mathcal{L}(\mathcal{F}_{s}) - \mathcal{L}(I_{i,s})|=o(\delta_{i})\mathcal{L}(\mathcal{F}_{s})$ also implies that the intervals $I_{i,s}(x - a_{i})$ are contained within some compact interval $[-B,B]$ if $s \in \Omega(\eta) = [\eta, \| \mathcal{F}\|_{\infty} - \eta] \subset \Omega(\eta_{i}) = [\eta_{i}, \| \mathcal{F}\|_{\infty} - \eta_{i}]$. By Lemma \ref{lem:precompact}, there exist a subsequence of $\{ \int_{\eta}^{\|\mathcal{F}\|_{\infty} - \eta} 1_{I_{i,s}}(\cdot - a_{i}) \, ds \}_{i}$, for which we still use the same subscript, and $\mathcal{G}\in L^{p_{k}}(\mathbb{R})$, such that,
\begin{equation} \label{limi}
    \lim_{i\to\infty} \| \int_{\eta}^{\|\mathcal{F}\|_{\infty} - \eta} 1_{I_{i,s}}(\cdot - a_{i}) \, ds - \mathcal{G}\|_{p_{k}} = 0.
\end{equation}
Combining \eqref{etatailI} with \eqref{limi}, we have:
\begin{equation*} 
    \limsup_{i\to\infty} \| h_{i}(\cdot - a_{i}) - \mathcal{G}\|_{p_{k}} \leq C\epsilon.
\end{equation*} 
Since this holds for every $\epsilon$, we have the conclusion.     
\end{proof}

\noindent
From the conclusions of Lemma \ref{lem:2funcomp} and Proposition \ref{prop:shiftcomp}:
\begin{equation}  \label{Gcomp1}
    \| f_{i}(\cdot - a_{i}) -\mathcal{G} \|_{p_{k}} \leq \| (f_{i} - h_{i})(\cdot - a_{i})\|_{p_{k}} + \| h_{i}(\cdot - a_{i}) - \mathcal{G}\|_{p_{k}}\to 0
\end{equation}
as $i\to\infty$, for some $\mathcal{G}\in L^{p_{k}}(\mathbb{R})$. Since $\{f_{i}\}_{i}$ is an extremizing sequence:
\begin{equation} \label{Gcomp2}
    \| f_{i}(\cdot - a_{i})\|_{U^{k}} \geq (1-\delta_{i})A(k)\|f_{i}(\cdot - a_{i})\|_{p_{k}}
\end{equation}
\eqref{Gcomp1}, \eqref{Gcomp2} and Gowers-Host-Kra norm inequality give:
\begin{equation*} 
    \| \mathcal{G} \|_{U^{k}} = A(k)\|\mathcal{G}\|_{p_{k}}
\end{equation*}
which means $\mathcal{G}$ must be a Gaussian, by the characterization of nonnegative extremizers of the Gowers-Host-Kra norm inequality. We've now finished the proof of Theorem \ref{thm:main} for nonnegative near extremizers in one dimension.

\subsection{A remark on admissibility}

Admissibility plays a central role, as demonstrated, in obtaining our result in Subsection \ref{sec:lemHpf}. Admissibility condition highlights the interrelations between the intervals involved, in terms of their centers and lengths. For instance, it's equivalent to the condition that the functional $\Psi$ strictly decreases once a subset of $M$ centered intervals $J_{i}$ is translated. Admissibility can be considered a boundary case of strict admissibility, whose definition speaks of the following: The lengths $l_{i}$ must be selected so as to the condition $|L_{m}(\vec{x})|\leq l_{m}$ is not redundant for any particular $m$. In the case $N=2$, $M=3$ and $L_1((x_1,x_2))=x_1, L_2((x_1,x_2))=x_2, L_3((x_1,x_2)) = x_1+x_2$, we recover the classic Riesz-Sobolev inequality \cite{lieb2001analysis}. In this case, strict admissibility is equivalent to the strict admissibility polygon condition given by Burchard \cite{burchard1996cases}, which is: $l_{i} < l_{j} + l_{k}$ for any permutation $(i, j, k)$ of $(1,2,3)$. Indeed, if the set $\{(x_1,x_2): |x_1| \leq l_1, |x_2| \leq l_2, |x_1+x_2|\leq l_3\}$ is a nonempty proper subset of $\{(x_1,x_2): |x_1|\leq l_1, |x_2| \leq l_2\}$, i.e. the condition $|L_3((x_1,x_2))|\leq l_3$ is not redundant, then $l_3 > l_2+l_1$, and conversely.

\section{Nonnegative extremizers of Gowers product inequality}

\noindent
Our induction step to higher dimensions will need a complete characterization of an arbitrary nonnegative extremizing tuple $\vec{f} = (f_{\alpha}:\alpha\in\{0,1\}^{k})$ of the Gowers product inequality in one dimension. In fact, we obtain a characterization for all dimensions.\\

\noindent
Let $m, n\geq 1$ be integers. For $1\leq i\leq m$, let $H, H_{i}$ be vector spaces, $B_{i}: H\to H_{i}$ be a surjective linear mapping, with $ker(B_{i})\cap ker(B_{j}) = \{0\}$, $1\leq i\not= j\leq m$, and $1\leq p_{i}\leq\infty$. Then $(\vec{B},\vec{p}) = ((B_1,\text{ ... }, B_{m}), (p_1,\text{ ... }, p_{m}))$ is called a Brascamp-Lieb datum \cite{bennett2008brascamp}. A Brascamp-Lieb inequality is an inequality of the form:
\begin{equation*} 
    \int_{H} \prod_{1\leq i\leq m} f_{i}\circ B_{i}(x)\, dx \leq BL(\vec{B},\vec{p})\|f_{i}\|_{L^{p_{i}}(H_{i})}.
\end{equation*}
$BL(\vec{B},\vec{p})$ is the smallest constant such that the above inequality is satisfied for all input tuples $\vec{f} = (f_1,\text{ ... }, f_{m})$ with measurable $f_{i}: H_{i}\to\mathbb{R}_{\geq 0}$. An extremizing tuple $\vec{f}$ is an input tuple with which the equal sign happens. We quote the following result:

\begin{theorem} \cite{bennett2008brascamp}
Let $(\vec{B},\vec{p})$ be an extremizable Brascamp-Lieb datum with $1\leq p_{i}<\infty$ for all $i$. Suppose also that $B_{i}^{*}H_{i}\cap B_{j}^{*}H_{j} = \{0\}$ whenever $1\leq i < j\leq m$. Then if $\vec{f} = (f_{i})$ is an extremizing input, then all the $f_{i}$ are Gaussians, thus there exist real numbers $C, c_{i}>0$, positive definite transformations $M_{i}: H_{i}\to H_{i}$, and points $x_{i}\in H_{i}$ such that $f_{i}(x) = c_{i}\exp (-C\langle A_{i}(x-x_{i}), (x-x_{i})\rangle_{H_{i}})$. Moreover, $x_{i} = B_{i}w$, for some $w\in H$.    
\end{theorem}

\noindent
Our example of $(\vec{B},\vec{p}) = ((B_{\alpha}:\alpha\in\{0,1\}^{k}), (p_{k}, \text{ ... }, p_{k}))$ with $B_{\alpha}:\mathbb{R}^{(k+1)n}\to\mathbb{R}^{n}$ defined by $B_{\alpha}(x,\vec{h}) = x+\alpha\cdot\vec{h}$ is a Brascamp-Lieb datum. Hence the quoted result allows us to say, if $\vec{f} = (f_{\alpha}:\alpha\in\{0,1\}^{k})$, $f_{\alpha}:\mathbb{R}^{n}\to\mathbb{R}_{\geq 0}$, is a nonnegative extremizing tuple for the inequality, 
\begin{equation*} 
    |\mathcal{T}_{k}(\vec{f})| = |\mathcal{T}_{k}(f_{\alpha}: \alpha\in\{0,1\}^{k})| \leq A(k,n)^{2^{k}} \prod_{\alpha\in\{0,1\}^{k}} \|f_{\alpha}\|_{p_{k}} = A(k,n)^{2^{k}}
\end{equation*}
then each $f_{\alpha}$ is a Gaussian of the described form, provided that the hypotheses are met. Indeed, for each $\alpha\in\{0,1\}^{k}$, define $B_{\alpha}^{*}:\mathbb{R}^{n}\to\mathbb{R}^{(k+1)n}$ by $B_{\alpha}^{*}(x) = (x, \vec{X}^{\alpha})$ with $\vec{X}^{\alpha} = (X^{\alpha}_1, \text{ ... }, X^{\alpha}_{k})$ and $\mathbb{R}^{n}\ni X^{\alpha}_{i} = x$ if $\alpha_{i} = 1$ and $X^{\alpha}_{i} = 0$ otherwise. It's easy to check that $B^{*}_{\alpha}$ is indeed the adjoint of $B_{\alpha}$ and that $B^{*}_{\alpha}\mathbb{R}^{n}\cap B^{*}_{\beta}\mathbb{R}^{n} = \{0\}$, $\alpha\not=\beta\in\{0,1\}^{k}$. As indicated by the Gowers-Host-Kra norm inequality, the datum $((B_{\alpha}:\alpha\in\{0,1\}^{k}), (p_{k}, \text{ ... }, p_{k}))$ is extremizable. Hence $f_{\alpha}(x) = m_{\alpha}\exp (-C\langle M_{\alpha}(x-c_{\alpha}), (x-c_{\alpha})\rangle_{\mathbb{R}^{n}})$, with $C, m_{\alpha}>0$, $M_{\alpha}:\mathbb{R}^{n}\to\mathbb{R}^{n}$ a positive definite transformation, and $c_{\alpha} = c_0 + \alpha\cdot\vec{c}$, for some $(c_0,\vec{c})\in\mathbb{R}^{(k+1)n}$.\\

\noindent 
Furthermore, in our case, we claim that there exists a positive definite transformation $M$ such that $M_{\alpha} = M$ for all $\alpha\in\{0,1\}^{k}$. To see this, note that if $\vec{f}=(f_{\alpha}:\alpha\in\{0,1\}^{k})$ is an extremizing tuple then so is its symmetric rearrangement tuple, $\vec{f}^{*} = (f_{\alpha}^{*}:\alpha\in\{0,1\}^{k})$; this is simply a consequence of the general rearrangement inequality, $\mathcal{T}_{k}(\vec{f}) = \mathcal{T}_{k}(f_{\alpha}: \alpha\in\{0,1\}^{k}) \leq \mathcal{T}_{k}(f_{\alpha}^{*}: \alpha\in\{0,1\}^{k}) = \mathcal{T}_{k}(\vec{f}^{*})$ and the fact that $\|f_{\alpha}\|_{p_{k}} = \|f^{*}_{\alpha}\|_{p_{k}}$. It then suffices to assume $f_{\alpha}(x)= m_{\alpha}\exp(-C\langle M_{\alpha}x, x\rangle_{\mathbb{R}^{n}})$. For $k=2$, the Cauchy-Schwarz inequality, Young's convolution inequality and the fact that $A(2,n) = ((C_{4/3}^2/C_2)^{1/2})^{n} A(1,n)^{1/2} = ((C_{4/3}^2/C_2)^{1/2})^{n}$ - see again Section \ref{sec:basics} - give us:
\begin{equation*} 
    \mathcal{T}_2(f_1, f_2, f_3, f_4) = \int_{\mathbb{R}^{n}} (f_1 \ast f_2)\cdot (f_3 \ast f_4)(x)\, dx \leq \| f_1 \ast f_2 \|_2 \| f_3 \ast f_4 \|_2 \leq A(2)^4 \prod_{i=1}^4 \| f_{i} \|_{4/3}.
\end{equation*}
The equal sign is a simple consequence of change of variables. The conclusion for the case $k=2$ then follows from the characterization of extremizers of the Cauchy-Schwarz inequality and of Young's convolution inequality. Assume the claim is true for index $k$ or lower. Then:
\begin{align*}
    \mathcal{T}_{k+1} (f_{\gamma}: \gamma\in\{0,1\}^{k+1}) &= \int_{\mathbb{R}^{n}} \mathcal{T}_{k} (T^{h} f_{(\alpha,1)} \cdot f_{(\alpha, 0)}: \alpha\in\{0,1\}^{k}) \, dh \\
    &\leq \int_{\mathbb{R}^{n}} \prod_{\alpha\in\{0,1\}^{k}} \| T^{h}f_{(\alpha,1)} \cdot f_{(\alpha,0)} \|_{U^{k}} \, dh\\
    &= A(k,n)^{2^{k}} \int_{\mathbb{R}^{n}} \prod_{\alpha\in\{0,1\}^{k}} \| T^{h} f_{(\alpha,1)} \cdot f_{(\alpha,0)}\|_{p_{k}} \, dh \\
    &\leq A(k,n)^{2^{k}} \prod_{\alpha\in\{0,1\}^{k}} \left ( \int_{\mathbb{R}^{n}} \| T^{h}f_{(\alpha,1)} \cdot f_{(\alpha,0)} \|_{p_{k}}^{2^{k}} \, dh \right )^{1/2^{k}} \\
    &\leq A(k+1,n)^{2^{k+1}} \prod_{\alpha\in\{0,1\}^{k}} \| f_{(\alpha, 0)} \|_{p_{k+1}} \| f_{(\alpha,1)} \|_{p_{k+1}} \\
    &= A(k+1,n)^{2^{k+1}} \prod_{\gamma\in\{0,1\}^{k+1}} \| f_{\gamma} \|_{p_{k+1}}.
\end{align*}

\noindent
The first inequality is the Gowers product inequality. The second equality is due to the fact that Gaussians are extremizers of the Gowers-Host-Kra norm inequality. The third inequality follows from H\"older's inequality and the fourth from the sharp Young's inequality, as discussed in Section \ref{sec:basics}. If $(f_{\gamma}: \gamma\in\{0,1\}^{k+1})$ is an extremizing tuple with $f_{\gamma}(x) = m_{\gamma}\exp(-C\langle M_{\gamma}x, x\rangle_{\mathbb{R}^{n}})$, this forces the first and third inequalities in the display to become equalities. Due to the Gowers-Host-Kra norm inequality and the fact that all the integrands involved are continuous, equal sign in the first inequality happens only when for all $h\in\mathbb{R}^{n}$ and $\alpha\in\{0,1\}^{k}$,
\begin{equation} \label{multlinexpand}
    \mathcal{T}_{k} (T^{h} f_{(\alpha,1)} \cdot f_{(\alpha,0)}) = A(k,n)^{2^{k}}\prod_{\alpha\in\{0,1\}^{k}} \| T^{h}f_{(\alpha,1)} \cdot f_{(\alpha,0)} \|_{p_{k}}.
\end{equation}
$T^{h}f_{(\alpha,1)}\cdot f_{(\alpha,0)}$ is still a Gaussian; hence by the induction hypothesis, \eqref{multlinexpand} implies in particular for $h=0$:
\begin{equation} \label{M}
    M_{(\alpha,0)} + M_{(\alpha,1)} = M_{(\beta,0)}+ M_{(\beta,1)} 
\end{equation}
for all $\alpha,\beta\in\{0,1\}^{k}$. On the other hand, equal sign in the third inequality gives, for each $\alpha\in\{0,1\}^{k}$:
\begin{equation} \label{Tksplit}
    \left ( \int_{\mathbb{R}^{n}} \| T^{h}f_{(\alpha,1)} \cdot f_{(\alpha, 0)} \|_{p_{k}}^{2^{k}} \, dh \right )^{1/2^{k}}= A(k,n)^{2^{k}} \| f_{(\alpha,0)} \|_{p_{k}} \| f_{(\alpha,1)} \|_{p_{k}}.
\end{equation}
The characterization of extremizers of the sharp Young's inequality \cite{eisner2012large} and \eqref{Tksplit} implies there exist $m_0, m_1\in\mathbb{R}_{>0}$, a positive definite transformation $M:\mathbb{R}^{n}\to\mathbb{R}^{n}$ and $c\in\mathbb{R}^{n}$ such that $\exp (-p_{k}\langle M(x-c), (x-c)\rangle_{\mathbb{R}^{n}}) = m_0 f_{(\alpha,0)} (x) = m_1 f_{(\alpha,1)} (x)$, which in turns implies 
\begin{equation} \label{2M}
    M_{(\alpha,0)} = M_{(\alpha,1)} = M
\end{equation} 
for all $\alpha\in\{0,1\}^{k}$. Then \eqref{M} and \eqref{2M} conclude $M_{\gamma} = M_{\gamma'}$ for all $\gamma, \gamma' \in \{0,1\}^{k+1}$ and hence the induction step. In particular, we obtain:

\begin{corollary}
Let $k\geq 2, n\geq 1$ be integers and $\vec{f} = (f_{\alpha}: \alpha\in\{0,1\}^{k})$ with measurable $f_{\alpha}:\mathbb{R}^{n}\to\mathbb{R}_{\geq 0}$, $\|f_{\alpha}\|_{p_{k}} = 1$. If 
\begin{equation*} 
    \mathcal{T}_{k}(\vec{f}) = \int_{\mathbb{R}^{(k+1)n}} \prod_{\alpha\in\{0,1\}^{k}} f_{\alpha}(x+\alpha\cdot\vec{h}) \, dxd\vec{h} = A(k,n)^{2^{k}}\prod_{\alpha\in\{0,1\}^{k}} \|f_{\alpha}\|_{p_{k}}
\end{equation*}
then $f_{\alpha} = m\exp(-\langle M(x-c_{\alpha}), (x-c_{\alpha})\rangle_{\mathbb{R}^{n}})$, with $m>0$, $M:\mathbb{R}^{n}\to\mathbb{R}^{n}$ a positive definite transformation, and $c_{\alpha} = c_0 + \alpha\cdot\vec{c}$, for some $(c_0,\vec{c})\in\mathbb{R}^{(k+1)n}$.
\end{corollary}

\begin{remark}
Another proof for the fact that $c_{\alpha} = c_0 + \alpha\cdot\vec{c}$, $\alpha\in\{0,1\}^{k}$, for some $(c_0,\vec{c})\in\mathbb{R}^{k+1}$ (hence only applicable for $n=1$) is as follows. If $f_{\alpha}(x) = m\exp(-a(x-c_{\alpha})^2)$ and
\begin{equation*} 
    \mathcal{T}_{k}(\vec{f}) =  \int \mathcal{T}_{k}(1_{F_{\alpha, v_{\alpha}}}: \alpha\in\{0,1\}^{k}) \, d\vec{v} = A(k)^{2^{k}} = \mathcal{T}_{k}(\vec{f}^{*}) = \int \mathcal{T}_{k}(1_{F^{*}_{\alpha, v_{\alpha}}}: \alpha\in\{0,1\}^{k}) \, d\vec{v}
\end{equation*}
then this entails $\mathcal{T}_{k}(1_{F_{\alpha, v}}: \alpha\in\{0,1\}^{k}) = \mathcal{T}_{k}(1_{F^{*}_{\alpha,v}}:\alpha\in\{0,1\}^{k})$ for all $v\in (0,m)$. Note that $(F_{\alpha,v})^{*} = F^{*}_{\alpha,v}$. Fix $v$ and let $\mathcal{L}(F_{\alpha,v}) = \mathcal{L}(F^{*}_{\alpha,v}) = l_{\alpha}$. Then apply Lemma \ref{lem:personal} to the admissble tuple $(B_{\alpha},l_{\alpha})_{\alpha\in\{0,1\}^{k}}$ to obtain the desired conclusion.   
\end{remark}

\subsection{Gaussian near extremizers in one dimension} \label{sec:nearG}

We now characterize Gaussian near extremizers of the Gowers product inequality in one dimension. Given a Gaussian tuple $\vec{f} = (f_{\alpha} = ca_{\alpha}^{1/2p_{k}}\exp(-a_{\alpha}(x-c_{\alpha})^2): \alpha\in\{0,1\}^{k})$, so that $\|f_{\alpha}\|_{p_{k}} = 1$. Suppose that for some $\delta >0$ small,
\begin{equation*} 
    \mathcal{T}_{k}(\vec{f}) = \mathcal{T}_{k}(f_{\alpha}: \alpha\in\{0,1\}^{k}) \geq (1-\delta) A(k)^{2^{k}} \prod_{\alpha\in\{0,1\}^{k}} \|f_{\alpha}\|_{p_{k}} = (1-\delta) A(k)^{2^{k}}.
\end{equation*}
We claim that there exist $a, \Gamma = \Gamma(k) >0$ and a nonnegative function $\eta = \eta_{k}(\delta)$ that is increasing for small values of $\delta$ such that $|a_{\alpha}/a - \Gamma|\leq\eta$, for all $\alpha\in\{0,1\}^{k}$. As before, we can first assume that $c_{\alpha} = 0$ for all $\alpha\in\{0,1\}^{k}$. We start with the induction step for the case $k+1$ and assume the claim is true for index $k$ or lower. If:
\begin{equation*} 
    \mathcal{T}_{k+1}(f_{\beta}:\beta\in\{0,1\}^{k+1}) \geq (1-\delta) A(k+1)^{2^{k+1}}.
\end{equation*}
It then follows from the calculations in Section \ref{sec:shortpf} that there exists $c(k)>0$ such that for each $\alpha\in\{0,1\}^{k}$, 
\begin{equation*} 
    \| f^{p_{k}}_{(\alpha,0)}\ast f^{p_{k}}_{(\alpha,1)}\|_{k+1}^{k+1}\geq (1-c(k)\delta)Y(k+1)^{k+1}\|f^{p_{k}}_{(\alpha,0)}\|_{q}^{k+1}\|f_{(\alpha,1)}^{p_{k}}\|_{q}^{k+1}.
\end{equation*} 
Here, $q = p_{k+1}/p_{k}$ and $Y(k+1)$ is the optimal constant of Young's convolution inequality for the involved exponents. Since all the functions involved are centered Gaussians in one dimension and $p_{k} = 2^{k}/(k+1)$, it's a simple calculation to show that there exist $\delta_0 >0$ sufficiently small, a function $\eta = \eta_{k}(\delta)$ that is increasing for $0<\delta\leq\delta_0$ and $\Gamma = \Gamma(k)$, such that $|a_{(\alpha,0)}/a_{(\alpha,1)} - \Gamma|\leq\eta$, for all $\alpha\in\{0,1\}^{k}$. By symmetry, this means that if $\beta,\beta'\in\{0,1\}^{k+1}$ are such that $\beta_{i}=\beta'_{i}$ for all but one single index $i\in\{1,\text{ ... }, k+1\}$, then $|a_{\beta}/a_{\beta'} - \Gamma|\leq\eta$; in other words, there exists $a>0$ such that $|a_{\beta}/a-\Gamma|\leq\eta$ for all $\beta\in\{0,1\}^{k+1}$. The case $k=2$ is proved using similar arguments. \\

\noindent
We now claim that if $\delta$ is small enough then there exist $c_0\in\mathbb{R}$ and $\vec{c}\in\mathbb{R}^{k}$ such that $|c_{\alpha} - (c_0 + \alpha\cdot\vec{c})| = o(\delta)$. To see this, we now assume $a_{\alpha} = a>0$ for all $\alpha\in\{0,1\}^{k}$ as permitted by the above reasoning. For each $\alpha\in\{0,1\}^{k}$, let $\vec{\xi}^{\alpha}\in\mathbb{R}^{k+1}$ be such that $B_{\alpha}(\vec{\xi}^{\alpha}) = c_{\alpha}$. Define $\mathbb{R}^{k+1}\ni \vec{T}$ so that $\sum_{\alpha\in\{0,1\}^{k}} B^{*}_{\alpha} B_{\alpha}\vec{T} = \sum_{\alpha\in\{0,1\}^{k}} B_{\alpha}^{*}B_{\alpha}\vec{\xi}^{\alpha}$ (it's an easy calculation that $|det(\sum_{\alpha\in\{0,1\}^{k}} B_{\alpha}^{*}B_{\alpha})|>0$). Here, $B_{\alpha}$ is defined as above with $n=1$. Then by a change of variables:
\begin{align}
    \nonumber \int_{\mathbb{R}^{k+1}} \prod_{\alpha\in\{0,1\}^{k}} f_{\alpha}\circ B_{\alpha}(\vec{x})\, d\vec{x} &= C\int_{\mathbb{R}^{k+1}} \exp \big\{ -a\sum_{\alpha\in\{0,1\}^{k}} (B_{\alpha}(\vec{x}-\vec{\xi}^{\alpha}))^2 \big \} \, d\vec{x} \\
    \nonumber &= C\int_{\mathbb{R}^{k+1}} \exp \big\{ -a\sum_{\alpha\in\{0,1\}^{k}} (B_{\alpha}(\vec{x}-(\vec{\xi}^{\alpha}-\vec{T})))^2 \big \} \, d\vec{x} \\
    \nonumber &= C\exp \big\{ -a\sum_{\alpha\in\{0,1\}^{k}} (B_{\alpha}(\vec{\xi}^{\alpha}-\vec{T}))^2 \big \}\int_{\mathbb{R}^{k+1}} \exp \big\{ -a\sum_{\alpha\in\{0,1\}^{k}} (B_{\alpha}(\vec{x}))^2 \big \} \, d\vec{x}\\
    \label{fB} &\leq C\int_{\mathbb{R}^{k+1}} \exp \big\{ -a\sum_{\alpha\in\{0,1\}^{k}} (B_{\alpha}(\vec{x}))^2 \big \} \, d\vec{x}. 
\end{align}
Note that by definition of $\vec{T}$, $\int_{\mathbb{R}^{k+1}} \exp \big\{ 2a\sum_{\alpha\in\{0,1\}^{k}} B_{\alpha}(\vec{x})\cdot B_{\alpha}(\vec{\xi}^{\alpha}-\vec{T})\big\} \, d\vec{x} = 1$, hence the last equality in \eqref{fB} follows. It's also now clear from \eqref{fB} and the fact that its last expression is $\|c\exp(-ax^2)\|_{U^{k}} = A(k)^{2^{k}}\|c\exp(-ax^2)\|_{p_{k}}^{2^{k}} = A(k)^{2^{k}}$, that if $\mathcal{T}_{k}(\vec{f}) = \int_{\mathbb{R}^{k+1}} \prod_{\alpha\in\{0,1\}^{k}} f_{\alpha}\circ B_{\alpha}(\vec{x})\, d\vec{x}$ is nearing its optimal value then $|c_{\alpha} - B_{\alpha}(\vec{T})| = o(\delta)$ for all $\alpha\in\{0,1\}^{k}$. Hence the claim follows. 

\begin{remark}
If $\vec{f} = (f_{\alpha}:\alpha\in\{0,1\}^{k})$ is such that 
\begin{equation*} 
    \mathcal{T}_{k}(\vec{f})\geq (1-\delta) A(k)^{2^{k}}\prod_{\alpha\in\{0,1\}^{k}}\|f_{\alpha}\|_{p_{k}},
\end{equation*} 
then it follows from the Gowers-Cauchy-Schwarz inequality that for each $\alpha\in\{0,1\}^{k}$, $\|f_{\alpha}\|_{U^{k}}\geq (1-o(\delta))A(k)\|f_{\alpha}\|_{p_{k}}$. If furthermore, $f_{\alpha}\geq 0$ then by the result for dimension one, there exists a Gaussian $g_{\alpha}$ such that $\|g_{\alpha} - f_{\alpha}\|_{p_{k}} = o(\delta)\|f_{\alpha}\|_{p_{k}}$. Then from the Gowers product inequality, $\vec{g} = (g_{\alpha}:\alpha\in\{0,1\}^{k})$ is also a near extremizing tuple: $\mathcal{T}_{k}(\vec{g})\geq (1-o(\delta)) A(k)^{2^{k}}\prod_{\alpha\in\{0,1\}^{k}}\|g_{\alpha}\|_{p_{k}}$, and the analytic descriptions of the $g_{\alpha}$ are given above. We note that all of these arguments can be generalized to higher dimensions. For now, a characterization of nonnegative near extremizing tuples for the Gowers product inequality is sufficient for an induction process in Section \ref{sec:hidim} below. We also note that the arguments given in this section are stronger in the sense that they establish analytic properties of near extremizers, not just extremizers.
\end{remark}

\section{Extension to higher dimensions} \label{sec:hidim}

\subsection{An additive relation}

We first quote a result in \cite{christ2019near}:

\begin{proposition} \label{prop:christm} \cite{christ2019near} 
Let $n\geq 1$. There exists a positive constant $K = K(n)> 0$ with the following property. Let $B$ be a ball of positive, finite radius. Let $\alpha, \beta, \gamma:\mathbb{R}^{n}\to\mathbb{C}$ be measurable functions. Let $\tau\in (0,\infty)$ and $\delta\in (0,1]$. Suppose that,
\begin{equation*} 
    \mathcal{L}( \{(x,y) \in B^2: |\alpha(x) + \beta(y) - \gamma(x+y)| > \tau\} ) < \delta (\mathcal{L}(B))^2. 
\end{equation*}
Then there exists an affine function $L: \mathbb{R}^{n}\to\mathbb{C}$ such that
\begin{equation*} 
    \mathcal{L}(\{ x\in B: |\alpha(x) - L(x)| > K\tau\}) < K\delta \mathcal{L}(B).
\end{equation*}    
\end{proposition}

\noindent
Using this proposition, we prove the following:

\begin{proposition} \label{prop:nearlin}
Let $n\geq 1$ and $l\geq 2$ be an integer. Let $C_{i} > 0$ be such that $C_{i}\asymp C_{j}$ for $i,j\in\{1,\text{ ... }, l\}$. There exists $K > 0$ with the following property. Let $B_{i}$ be a ball in $\mathbb{R}^{n}$ such that $\mathcal{L}(B_{i})=C_{i}$ and let $f_{i}: \mathbb{R}^{n}\to\mathbb{C}$ be a measurable function, $i\in\{1,\text{ ... }, l\}$. Let $\tau\in (0,\infty)$ and $\delta\in (0,1]$. Suppose that,
\begin{equation} \label{nearlylin}
    |f_{l+1}(x_1+\text{ ... }+x_{l}) - \sum_{i=1}^{l} f_{i}(x_{i})| \leq \tau 
\end{equation}
for all $(x_1, \text{ ... }, x_{l}) \in B_1 \times \text{ ... } \times B_{j}$ outside a subset of measure $\delta \prod_{i=1}^{l} \mathcal{L}(B_{i})$. Then there exist an affine function $a:\mathbb{R}^{n}\to\mathbb{C}$ and a positive function $\eta$ satisfying $\lim_{\delta\to 0} \eta(\delta) = 0$, such that
\begin{equation*} 
    |f_1(x) - a(x)| \leq K\tau 
\end{equation*} 
for all $x\in B_1$ outside a subset of measure $K\eta(\delta)\mathcal{L}(B_1)$.    
\end{proposition}

\begin{proof}
Since the conclusion doesn't change after applying a finite number of translations and dilations, we assume that $B_1 = \text{ ... } = B_{l} = B$ - here $B$ is a ball of positive radius. Then \eqref{nearlylin} gives that, for all $(x_1, \text{ ... }, x_{l}) \in \times_{i=1}^{l} B$ outside a subset of measure $\delta(\mathcal{L}(B))^{l}$, the following holds:
\begin{equation} \label{nearlylin1}
    f_{l+1}(x_1+x_2 + \sum_{i=3}^{l} x_{i}) - \sum_{i=3}^{l}f_{i}(x_{i}) = f_1(x_1) + f_2(x_2) + O(\tau). 
\end{equation}
That means there exists a positive function $\eta$ satisfying $\lim_{\delta\to 0} \eta(\delta) = 0$ such that the following holds. For each $(x_3, \text{ ... }, x_{l}) \in\times_{i=3}^{l} B$ outside a subset of measure at most $\eta(\delta)(\mathcal{L}(B))^{l-2}$, \eqref{nearlylin1} is true for all $(x_1, x_2)\in B\times B$ except for a subset of measure at most $\eta(\delta)(\mathcal{L}(B))^2$. Take such a point $(x_3, \text{ ... }, x_{l}) \in\times_{i=3}^{l} B$. For this point, define $\tilde{f}: \mathbb{R}^{n}\to\mathbb{C}$ by $\tilde{f}_{l+1}(u) = f_{l+1}(u + \sum_{i=3}^{l} x_{i}) - \sum_{i=3}^{l}f_{i}(x_{i})$. Then \eqref{nearlylin1} becomes:
\begin{equation*} \tilde{f}_{l+1}(x_1+x_2) = f_1(x_1) + f_2(x_2) + O(\tau)\end{equation*}
for all $(x_1, x_2) \in B\times B$ outside a subset of measure at least $(1-\eta(\delta))(\mathcal{L}(B))^2$. Now apply Proposition \ref{prop:christm} to $\tilde{f}_{l+1}, f_1, f_2$ - and an appropriate translation and dilation if necessary - to obtain an affine function $a:\mathbb{R}^{n}\to\mathbb{C}$ such that $f_1(x) = a(x) + O(\tau)$ for all $x\in B_1$ outside a subset of measure $K\eta(\delta)\mathcal{L}(B_1)$.    
\end{proof}

\subsection{Extension to higher dimensions}

Let $n\geq 1$. Assume Theorem \ref{thm:main} is true for nonnegative functions and for dimensions $n$ and lower. Let $f(x,s): \mathbb{R}^{n}\times\mathbb{R} \to \mathbb{R}_{\geq 0}$ be a $(1-\delta)$ near extremizer of the Gowers-Host-Kra norm inequality. Assume that $\| f \|_{p_{k}} = 1$. Define $F: \mathbb{R}^{n}\to\mathbb{R}_{\geq 0}$ by $F(x) = \| f(x, \cdot) \|_{p_{k}}$. Define $f_{x}:\mathbb{R}\to\mathbb{R}_{\geq 0}$ by $f_{x}(s) = \frac{f(x,s)}{F(x)}$ if $F(x) \not\in\{0, \infty\}$, and $f_{x}(s) \equiv 0$ if $F(x) \in\{ 0, \infty\}$ - which happens only for a null subset of $spt(F) \subset \mathbb{R}^{n}$, outside of which, $\| f_{x} \|_{p_{k}} = \| f \|_{p_{k}} = 1$. Thus, $\| F\|_{p_{k}} = \| f\|_{p_{k}} = 1$. From the definition of the Gowers-Host-Kra norms and Fubini's theorem,
\begin{equation*} 
    \| f \|_{U^{k}}^{2^{k}} = \int_{\mathbb{R}^{(k+1)n}} \prod_{\alpha \in \{0,1\}^{k}} F(x+\alpha \cdot \vec{h}) \mathcal{T}_{k} (f_{x+\alpha\cdot \vec{h}}: \alpha\in\{0,1\}^{k}) \, dx d\vec{h}.
\end{equation*}
The Gowers-Cauchy-Schwarz inequality then gives,
\begin{align} 
    \nonumber (1-\delta)A(k,n+1)^{2^{k}} &\leq\int_{\mathbb{R}^{(k+1)n}} \prod_{\alpha \in \{0,1\}^{k}} F(x+\alpha \cdot \vec{h}) \mathcal{T}_{k} (f_{x+\alpha\cdot \vec{h}}: \alpha\in\{0,1\}^{k}) \, dx d\vec{h}\\
    \label{neardelta} &\leq A(k)^{2^{k}} \int_{\mathbb{R}^{(k+1)n}} \prod_{\alpha \in \{0,1\}^{k}} F(x+\alpha \cdot \vec{h}) \, dx d\vec{h} = A(k)^{2^{k}} \| F\|_{U^{k}}^{2^{k}}.
\end{align}
Since $A(k,m) = A(k)^{m}$, \eqref{neardelta} implies $\| F\|_{U^{k}} \geq (1-\delta)A(k,n)^{2^{k}}$. Then by the inductive assumption, there exists a Gaussian $\mathcal{F}:\mathbb{R}^{n}\to\mathbb{R}_{>0}$ such that $\|\mathcal{F}\|_{p_{k}} = 1$ and $\| \mathcal{F} - F\|_{p_{k}} = o(\delta)$.\\

\noindent
Let $B_{R}\subset\mathbb{R}^{n}$ denote a centered ball of radius $R$. Take $\eta > 0$ small. There exist $\delta = \delta(\eta) > 0$ small and $B_{R}$ with $R = R(\eta) > 0$ large such that $\delta\to 0$ and $R\to\infty$ as $\eta\to 0$ with the following extra property. Let $\mathcal{F}^{\dagger}$ denote the standard centered Gaussian on $\mathbb{R}^{n}$. Suppose $F^{\dagger}$ is such that $\|F^{\dagger}\|_{p_{k}} = 1$ and $\| \mathcal{F}^{\dagger} - F^{\dagger}\|_{p_{k}} \leq \delta$. Then:
\begin{equation} \label{Fdagger}
    \int_{\mathbb{R}^{(k+1)n}\setminus \tilde{B}_{R}^{k+1}} \prod_{\alpha \in \{0,1\}^{k}} F^{\dagger}(x+\alpha \cdot \vec{h}) \, dx d\vec{h} = \int_{\mathbb{R}^{(k+1)n}} \prod_{\alpha\in\{0,1\}^{k}} F^{\dagger}\chi_{\mathbb{R}^{n}\setminus B_{R}} (x+\alpha\cdot\vec{h}) \, dx d\vec{h}< \eta.
\end{equation}
Indeed, since $\mathcal{F}^{\dagger}$ is the standard centered Gaussian, for every $\eta > 0$, there exists $B_{R}$ with $R$ sufficiently large so that:
\begin{equation} \label{Fdaggereta}
    \int_{\mathbb{R}^{(k+1)n}\setminus \tilde{B}_{R}^{k+1}} \prod_{\alpha \in \{0,1\}^{k}} \mathcal{F}^{\dagger} (x+\alpha \cdot \vec{h}) \, dx d\vec{h} = \int_{\mathbb{R}^{(k+1)n}} \prod_{\alpha\in\{0,1\}^{k}} \mathcal{F}^{\dagger}\chi_{\mathbb{R}^{n}\setminus B_{R}} (x+\alpha\cdot\vec{h}) \, dx d\vec{h} < \eta/2. 
\end{equation}
Applying the Gowers product inequality to have:
\begin{equation} \label{Fdaggerdelta}
    \int_{\mathbb{R}^{(k+1)n}} \prod_{\alpha\in\{0,1\}^{k}} |\mathcal{F}^{\dagger}- F^{\dagger}|\chi_{\mathbb{R}^{n}\setminus B_{R}} (x+\alpha\cdot\vec{h}) \, dx d\vec{h} \leq A(k,n)^{2^{k}}\prod_{\alpha\in\{0,1\}^{k}}\|\mathcal{F}^{\dagger}-F^{\dagger}\|_{p_{k}} < A(k,n)^{2^{k}} \delta. 
\end{equation}
Simply take $\delta$ small enough so that $A(k,n)^{2^{k}} \delta < \eta/2$. Then \eqref{Fdagger} follows from \eqref{Fdaggereta} and \eqref{Fdaggerdelta}.\\

\noindent
Suppose now that our Gaussian $\mathcal{F}$ is indeed the standard centered one; this can be done by applying an affine transformation and scaling by a suitable constant. Then there exists a sufficiently large $R$ so that
\begin{equation} \label{Finteta}
    \int_{\mathbb{R}^{(k+1)n}\setminus \tilde{B}_{R}^{k+1}} \prod_{\alpha \in \{0,1\}^{k}} F(x+\alpha \cdot \vec{h}) \, dx d\vec{h} < \eta.
\end{equation}
Now by the Gowers product inequality:
\begin{equation*} 
    \mathcal{T}_{k}(f_{x+\alpha \cdot \vec{h}}: \alpha\in\{0,1\}^{k}) \leq A(k)^{2^{k}} \prod_{\alpha\in\{0,1\}^{k}} \| f_{x+\alpha\cdot\vec{h}}\|_{p_{k}} = A(k)^{2^{k}}
\end{equation*}
which then with \eqref{Finteta}, implies
\begin{equation} \label{Fintmulteta}
    \int_{\mathbb{R}^{(k+1)n}\setminus \tilde{B}_{R}^{k+1}} \prod_{\alpha \in \{0,1\}^{k}} F(x+\alpha \cdot \vec{h}) \mathcal{T}_{k}(f_{x+\alpha\cdot h}:\alpha\in\{0,1\}^{k}) \, dx d\vec{h} < \eta A(k)^{2^{k}}.
\end{equation}
Note that if $R\to\infty$ and $\delta\to 0$ then $\eta = \eta(R,\delta)\to 0$ in \eqref{Finteta} and \eqref{Fintmulteta}. We've set up the case for restricting our analysis in $\tilde{B}_{R}^{k+1}$ with $R$ big enough.

\begin{lemma} \label{lem:extraset}
Let $f, F, f_{x}, R$ be as above. Then there exist $\Omega \subset \tilde{B}_{R}^{k+1}$ and $\omega \subset B_{R}$ such that $\mathcal{L}(\Omega) + \mathcal{L}(\omega) = o(\delta)$ satisfying:\\
For $(x,\vec{h})\in\tilde{B}_{R}^{k+1}\setminus\Omega$ and $\alpha\in\{0,1\}^{k}$,
\begin{equation} \label{lemextraset1}
    \| \phi_{x+\alpha\cdot\vec{h}} -f_{x+\alpha\cdot\vec{h}} \|_{p_{k}} = o(\delta)
\end{equation}
with $\phi_{x+\alpha\cdot\vec{h}}(s) = ca(x+\alpha\cdot\vec{h})^{1/2p} \exp (-a(x+\alpha\cdot\vec{h})(s-d(x+\alpha\cdot\vec{h}))^2)$. The functions $a: \mathbb{R}^{n}\to\mathbb{R}_{>0}$ and $d: \mathbb{R}^{n}\to\mathbb{R}$ are measurable and satisfy the following properties:\\
There exists a scalar $a > 0$ such that for $x \in B_{R}\setminus\omega$,
\begin{equation} \label{lemextraset2}
    | a(x) - a| = o(\delta).
\end{equation} 
For $(x,\vec{h}) \in \tilde{B}_{R}^{k+1}\setminus\Omega$,
\begin{equation} \label{lemextraset3}
    | d(x+\alpha\cdot\vec{h}) - (1,\alpha)\cdot (d(x), d(h_1), \text{ ... }, d(h_{k})) | = o(\delta).
\end{equation}
Lastly,
\begin{equation} \label{lemextraset4}
    \int_{\Omega} \prod_{\alpha\in\{0,1\}^{k}} F(x+\alpha \cdot \vec{h}) \, dx d\vec{h} = o(\delta). 
\end{equation}    
\end{lemma}

\begin{proof}
From the near extremizing hypothesis and the Gowers-Host-Kra norm inequality,
\begin{align} 
    \nonumber \| f\|_{U^{k}}^{2^{k}} &= \int_{\mathbb{R}^{(k+1)n}} \prod_{\alpha\in\{0,1\}^{k}} F(x+\alpha \cdot \vec{h}) \mathcal{T}_{k}(f_{x+\alpha \cdot \vec{h}}: \alpha\in\{0,1\}^{k}) \, dx d\vec{h} \\ 
    \nonumber &\geq (1-\delta)A(k,n+1)^{2^{k}} = (1-\delta)A(k)^{2^{k}}A(k,n)^{2^{k}} \| F\|_{p_{k}}^{2^{k}} \\ 
    \label{fUk} &\geq (1-\delta)A(k)^{2^{k}} \| F\|_{U^{k}}^{2^{k}} = (1-\delta)A(k)^{2^{k}} \int_{\mathbb{R}^{(k+1)n}} \prod_{\alpha\in\{0,1\}^{k}} F(x+\alpha \cdot \vec{h}) \, dx d\vec{h}.
\end{align}
Take $R$ appropriately so that, as in \eqref{Finteta} and \eqref{Fintmulteta}, 
\begin{align}
    \label{odelta1} \int_{\mathbb{R}^{(k+1)n}\setminus\tilde{B}_{R}^{k+1}} \prod_{\alpha\in\{0,1\}^{k}} F(x+\alpha \cdot \vec{h}) \,dx d\vec{h} &= o(\delta)\\
    \label{odelta2} \int_{\mathbb{R}^{(k+1)n}\setminus \tilde{B}_{R}^{k+1}} \prod_{\alpha \in \{0,1\}^{k}} F(x+\alpha \cdot \vec{h}) \mathcal{T}_{k}(f_{x+\alpha\cdot h}:\alpha\in\{0,1\}^{k}) \, dx d\vec{h} &= o(\delta).
\end{align}
There might still exist $\Omega_0\subset\tilde{B}_{R}^{k+1}$ such that, 
\begin{equation} \label{odelta3}
    \int_{\Omega_0} \prod_{\alpha\in\{0,1\}^{k}} F(x+ \alpha \cdot \vec{h}) \, dx d\vec{h} = o(\delta) 
\end{equation} 
and hence, due to the Gowers product inequality,
\begin{equation} \label{deltaAk}
    \int_{\Omega_0} \prod_{\alpha \in \{0,1\}^{k}} F(x+\alpha \cdot \vec{h}) \mathcal{T}_{k}(f_{x+\alpha\cdot h}:\alpha\in\{0,1\}^{k}) \, dx d\vec{h} =o(\delta) A(k)^{2^{k}}.
\end{equation} 
Since $F$ is $o(\delta)$ close in $L^{p_{k}}$ norm to the standard centered Gaussian on $\mathbb{R}^{n}$, $\mathcal{F}$, \eqref{odelta3} implies a similar inequality for $\mathcal{F}$: 
\begin{equation*} 
    \int_{\Omega_0} \prod_{\alpha\in\{0,1\}^{k}} \mathcal{F}(x+ \alpha \cdot \vec{h}) \, dx d\vec{h}  = o(\delta)
\end{equation*}
which then implies, 
\begin{equation*} 
    \mathcal{L}(\Omega_0) = o(\delta).
\end{equation*}
Now \eqref{fUk}, \eqref{odelta1}, \eqref{odelta2} and \eqref{deltaAk} imply that $\Omega_0$ also has the following properties:
\begin{multline} \label{odeltaAk}
    \int_{\tilde{B}_{R}^{k+1}\setminus\Omega_0} \prod_{\alpha\in\{0,1\}^{k}} F(x+\alpha \cdot \vec{h}) \mathcal{T}_{k}(f_{x+\alpha \cdot \vec{h}}: \alpha\in\{0,1\}^{k}) \, dx d\vec{h} \\ \geq (1-o(\delta))A(k)^{2^{k}} \int_{\tilde{B}_{R}^{k+1}\setminus\Omega_0} \prod_{\alpha\in\{0,1\}^{k}} F(x+\alpha \cdot \vec{h}) \, dx d\vec{h}.
\end{multline}
Since $\mathcal{T}_{k}(f_{x+\alpha \cdot \vec{h}}: \alpha\in\{0,1\}^{k}) \leq A(k)^{2^{k}}$ for all $(x,\vec{h})\in\mathbb{R}^{(k+1)n}$, by the Gowers inner product inequality, \eqref{deltaAk} and \eqref{odeltaAk} imply that for a.e $(x,\vec{h}) \in \tilde{B}_{R}^{k+1} \setminus \Omega_0$,
\begin{equation*} 
    \mathcal{T}_{k}(f_{x+\alpha \cdot \vec{h}}: \alpha\in\{0,1\}^{k}) \geq (1-o(\delta))A(k)^{2^{k}}
\end{equation*}
which, by the Gowers-Cauchy-Schwarz inequality, entails
\begin{equation*} 
    \prod_{\alpha\in\{0,1\}^{k}} \| f_{x+\alpha \cdot \vec{h}}\|_{U^{k}} \geq \mathcal{T}_{k}(f_{x+\alpha \cdot \vec{h}}: \alpha\in\{0,1\}^{k}) \geq (1-o(\delta))A(k)^{2^{k}} \prod_{\alpha\in\{0,1\}^{k}} \| f_{x+\alpha \cdot \vec{h}} \|_{p_{k}} 
\end{equation*}
which then gives that, for each $\alpha\in\{0,1\}^{k}$ and a.e $(x,\vec{h}) \in \tilde{B}_{R}^{k+1} \setminus \Omega_0$:
\begin{equation} \label{xalphah}
    \| f_{x+\alpha \cdot \vec{h}} \|_{U^{k}} \geq (1-o(\delta))A(k) \| f_{x+\alpha \cdot \vec{h}}\|_{p_{k}}.
\end{equation}
Excluding a null subset if necessary, then the inductive hypothesis for dimension $n=1$ and \eqref{xalphah} imply that, if $(x,\vec{h}) \in \tilde{B}_{R}^{k+1} \setminus \Omega_0$ then $f_{x+\alpha \cdot \vec{h}}$ is $o(\delta)$ close in $L^{p_{k}}$ norm to a Gaussian $\phi_{x+\alpha \cdot \vec{h}}(s) = ca(x+\alpha\cdot\vec{h})^{1/2p} \exp (-a(x+\alpha\cdot\vec{h})(s-d(x+\alpha\cdot\vec{h}))^2)$. \\

\noindent
For $x \in B_{R}$, let $E^{x} = \{ \vec{h}\in B_{R}^{k}: (x,\vec{h}) \in \tilde{B}_{R}^{k+1}\setminus\Omega_0\}$ and $\omega_0 = \{ x \in B_{R}: \mathcal{L}(E^{x}) \leq (1-\delta')(\mathcal{L}(B_{R}))^{k} \}$, with $\delta'= \delta'(\delta)$ satisfying $\delta'\to 0$ sufficiently slow compared to $\delta\to 0$ so that $\mathcal{L}(\omega_0) = o(\delta)\mathcal{L}(B_{R})$ and for every $x\in B_{R}\setminus \omega_0$, $\mathcal{L}(E^{x}) \geq (1-o(\delta))(\mathcal{L}(B_{R}))^{k}$. From the definitions of $\omega_0$ and $\Omega_0$, if $x \in B_{R} \setminus \omega_0$ and $\vec{h} \in E^{x}$, we have a decomposition, 
\begin{equation} \label{fphirho}
    f_{x+\alpha \cdot \vec{h}} = \phi_{x+\alpha \cdot \vec{h}} + \rho_{x+\alpha \cdot \vec{h}} 
\end{equation} 
for some Gaussian $\phi_{x+\alpha \cdot \vec{h}}(s) = ca(x+\alpha\cdot\vec{h})^{1/2p} \exp (-a(x+\alpha\cdot\vec{h})(s-d(x+\alpha\cdot\vec{h}))^2)$, and $\| \rho_{x+\alpha\cdot\vec{h}}\|_{p_{k}} = o(\delta)$. This decomposition can be done so that $(x,\vec{h})\mapsto a(x+\alpha\cdot\vec{h})$, $(x,\vec{h})\mapsto d(x+\alpha\cdot\vec{h})$, $((x,\vec{h}),s)\mapsto \rho_{x+\alpha\cdot\vec{h}}(s)$ are all measurable. Take $x\in B_{R}\setminus\omega_0$ and $\vec{h}\in E^{x}$. By \eqref{xalphah}, \eqref{fphirho} and the Gowers product inequality:
\begin{equation} \label{Tkphi}
    \mathcal{T}_{k}(\phi_{x+\alpha \cdot \vec{h}}: \alpha\in\{0,1\}^{k}) \geq (1-o(\delta))A(k)^{2^{k}}.
\end{equation}
By the result in Subsection \ref{sec:nearG}, it follows from \eqref{Tkphi} that for every $\alpha\in\{0,1\}^{k}$, $x\in B_{R}\setminus\omega_0$ and $\vec{h}\in E^{x}$,
\begin{align}
    \label{alin1} a(x) = a(x+\alpha\cdot\vec{h}) &+ o(\delta)\\
    \label{alin2} | d(x+\alpha\cdot\vec{h}) - (1,\alpha)\cdot (d(x), d(h_1), \text{ ... }, d(h_{k})) | &= o(\delta).
\end{align}
Since $\mathcal{L}(B_{R}\setminus\omega_0) = (1-o(\delta))\mathcal{L}(B_{R})$, there exists $x_0\in B_{R}\setminus\omega_0$ such that $|x_0| = o(\delta)$. Moreover, $\mathcal{L}(E^{x_0}) \geq (1-o(\delta))(\mathcal{L}(B_{R}))^{k+1}$ implies that the set $V$ of values $x_0 +\alpha\cdot\vec{h}$ for some $\alpha\in\{0,1\}^{k}$ and $\vec{h}\in E^{x_0}$ must take up a measure of $(1-o(\delta))\mathcal{L}(B_{R})$ in $B_{R}$. We define $\omega$ to be the complement of the intersection of $B_{R}\setminus\omega_0$ and $V$. Then by \eqref{alin1}, for every $x\in B_{R}\setminus\omega$,
\begin{equation} \label{alin3}
    a(x) = a(x_0) + o(\delta) = a + o(\delta)
\end{equation}
for some $a > 0$. Define $\Omega$ by letting $\tilde{B}_{R}^{k+1}\setminus\Omega$ to be the set $(x,\vec{h})\in\tilde{B}_{R}^{k+1}$ with $x\in\omega$ and $\vec{h}\in E^{x}$, then $\mathcal{L}(\Omega) + \mathcal{L}(\omega) = o(\delta)$. Hence \eqref{odelta3} is still retained with $\Omega$ in place of $\Omega_0$:
\begin{equation*} 
    \int_{\Omega} \prod_{\alpha\in\{0,1\}^{k}} F(x+ \alpha \cdot \vec{h}) \, dx d\vec{h} = o(\delta) 
\end{equation*} 
which is \eqref{lemextraset4}. Moreover, if $\alpha\in\{0,1\}^{k}$ and $(x,\vec{h})\in\tilde{B}_{R}^{k+1}\setminus\Omega$,
\begin{equation*} 
    \| \rho_{x+\alpha\cdot\vec{h}}\|_{p_{k}} = o(\delta)
\end{equation*} 
which is \eqref{lemextraset1}. Finally, \eqref{alin2} and \eqref{alin3} are \eqref{lemextraset3} and \eqref{lemextraset2}, respectively. Hence this completes the proof of Lemma \ref{lem:extraset}.    
\end{proof}

\begin{remark} \label{rem:decomp}
It's possible to select a decomposition in \eqref{fphirho} so that $(x,\vec{h})\mapsto a(x+\alpha\cdot\vec{h})$, $(x,\vec{h})\mapsto d(x+\alpha\cdot\vec{h})$, $((x,\vec{h}),s)\mapsto \rho_{x+\alpha\cdot\vec{h}}(s)$ are all measurable. Let $\delta >0$ be as above, then since $\mathcal{L}(\omega)=o(\delta)$, at each $x\in B_{R}\setminus\omega$, we can define $\phi_{x} = ca(x)^{1/2p}\exp(-a(x)(s-d(x))^2)$ in a way that both $a(x)$ and $d(x)$ are locally piecewise constant in a sufficiently small neighborhood $N_{x}$ of $x$, so that $\| \phi_{y}-f_{y}\|_{p_{k}} = o(\delta)$ for all $y\in N_{x}$. Then it's clear that the decomposition $f_{x} = \phi_{x} + \rho_{x}$ satisfies the conditions $\| \rho_{x}\|_{p_{k}} = o(\delta)$ and $x\mapsto\rho$ is measurable, if $x\in B_{R}\setminus\omega$.    
\end{remark}

\begin{theorem}
Let $f(x,s): \mathbb{R}^{n}\times\mathbb{R} \to \mathbb{R}_{\geq 0}$ be a $(1-\delta)$ near extremizer of the Gowers-Host-Kra norm inequality. Then there exists a Gaussian $\mathcal{G}:\mathbb{R}^{n}\times\mathbb{R}\to\mathbb{R}_{>0}$ such that $\| \mathcal{G}\|_{p_{k}} = 1$ and $\| \mathcal{G} - f\|_{p_{k}} = o(\delta)$.    
\end{theorem}

\begin{proof}
Let $\psi_{x}(s) = ca^{1/2p}\exp(-a(s-d(x))^2)$, with $a>0$ and the function $d$ be as above; $x\in B_{R}$. By \eqref{Tkphi}, \eqref{alin3} and the continuity of exponential functions, for $x\in B_{R}\setminus\omega$,
\begin{equation} \label{psif}
    \| \psi_{x} - f_{x}\|_{p_{k}} = o(\delta).
\end{equation}
An application of Proposition \ref{prop:nearlin} to \eqref{lemextraset3} gives for every $x\in B_{R}$ except for a subset of measure $o(\delta)\mathcal{L}(B_{R})$:
\begin{equation} \label{dL}
    d(x) = L(x) + o(\delta).
\end{equation}
Here $L$ is an affine function on $\mathbb{R}^{n}$, which we can take to be real-valued, as is the center function $d$. Let $\omega'$ be the union of $\omega$ and this new subset of measure $o(\delta)\mathcal{L}(B_{R})$; it's clear $\mathcal{L}(\omega') = o(\delta)$. Let $\varsigma_{x}(s) = ca^{1/2p}\exp(-a(s-L(x))^2)$. Then \eqref{dL} implies for $x\in B_{R}\setminus\omega'$:
\begin{equation} \label{psisigma}
    \| \psi_{x} - \varsigma_{x}\|_{p_{k}} = o(\delta).
\end{equation}
\eqref{psif} and \eqref{psisigma} then give $\| f_{x} - \varsigma_{x}\|_{p_{k}} = o(\delta)$, if $x\in B_{R}\setminus\omega'$. Let $\mathcal{G}(x,s) = \mathcal{F}(x)\varsigma_{x}(s)$.  It's clear that $\| \mathcal{G}\|_{p_{k}} = 1$. Recall that $f(x,s) = F(x)f_{x}(s)$. Then:
\begin{equation} \label{Gdiff}
    \| \mathcal{G} - f\|_{p_{k}}^{p_{k}} \leq C\int_{\mathbb{R}^{n}} F^{p_{k}}(x)\| f_{x} - \psi_{x}\|_{p_{k}}^{p_{k}} \, dx + C\int_{\mathbb{R}^{n}} F^{p_{k}}(x)\|\psi_{x} - \varsigma_{x}\|_{p_{k}}^{p_{k}} \, dx + C\| \mathcal{F} - F\|_{p_{k}}^{p_{k}}.
\end{equation}
The presence of the last term in \eqref{Gdiff} is due to the fact that $\varsigma_{x}(s)$ is a Schwartz function in terms of $s$ and of $L^{\infty}$ in terms of $x$. Moreover, since $\|\mathcal{F} - F\|_{p_{k}} = o(\delta)$, the contribution of this last term is that of size $o(\delta)$ in absolute value. For the first term and the second term, we split them as follows, respectively:
\begin{multline} \label{Fdiff1}
    \int_{\mathbb{R}^{n}} F^{p_{k}}(x)\| f_{x} - \psi_{x}\|_{p_{k}}^{p_{k}} \, dx \leq \int_{B_{R}\setminus\omega'} |\mathcal{F}^{p_{k}}(x)-F^{p_{k}}(x)|\| f_{x}-\psi_{x}\|_{p_{k}}^{p_{k}} \, dx \\ + \int_{B_{R}\setminus\omega'} \mathcal{F}^{p_{k}}(x)\| f_{x}-\psi_{x}\|_{p_{k}}^{p_{k}} \, dx + C\int_{(\mathbb{R}^{n}\setminus B_{R})\cup\omega'} F^{p_{k}}(x) \, dx
\end{multline}
and,
\begin{multline} \label{Fdiff2}
    \int_{\mathbb{R}^{n}} F^{p_{k}}(x)\| \psi_{x}-\varsigma_{x}\|_{p_{k}}^{p_{k}} \, dx \leq \int_{B_{R}\setminus\omega'} |\mathcal{F}^{p_{k}}(x)-F^{p_{k}}(x)|\| \psi_{x}-\varsigma_{x}\|_{p_{k}}^{p_{k}} \, dx \\ 
    + \int_{B_{R}\setminus\omega'} \mathcal{F}^{p_{k}}(x)\| \psi_{x} - \varsigma_{x}\|_{p_{k}}^{p_{k}} \, dx + C\int_{(\mathbb{R}^{n}\setminus B_{R})\cup\omega'} F^{p_{k}}(x) \, dx.
\end{multline}
The first and second terms in \eqref{Fdiff1} and \eqref{Fdiff2} has the size of $o(\delta)$ in absolute value due to \eqref{psif}, \eqref{psisigma} and the fact that $\|\mathcal{F}-F\|_{p_{k}} = o(\delta)$. The third terms can be further dominated by the sum $C\int_{(\mathbb{R}^{n}\setminus B_{R})\cup\omega'} \mathcal{F}^{p_{k}}(x) \, dx + \| \mathcal{F} - F\|_{p_{k}}$, the first term of which is of size $o(\delta)$ in absolute value by the choice of $R$ and the fact that $\mathcal{L}(\omega') = o(\delta)$. All of these yield $\|\mathcal{G}-f\|_{p_{k}} = o(\delta)$.    
\end{proof}

\noindent
With this theorem, the extension step to higher dimensions for nonnegative near extremizers is now complete.

\section{Complex-valued case} \label{sec:complexcase}

\subsection{Preparation} \label{sec:compprep}

Let $f:\mathbb{R}^{n}\to\mathbb{C}$ be a $(1-\delta)$ near extremizer of Gowers-Host-Kra norm inequality. We write $f(x) = |f|(x)a(x)$, with $a(x) = e^{i2\pi q(x)}$ and $q:\mathbb{R}^{n}\to\mathbb{R}/\mathbb{Z}$. Assume $\|f\|_{p_{k}} = 1$. We first make a few observations. Recall that if $f$ is a $(1-\delta)$ near extremizer then so is $|f|$:
\begin{equation} \label{modulusbd}
    \| |f| \|_{U^{k}} \geq \| f\|_{U^{k}} \geq (1-\delta)A(k,n)\|f\|_{p_{k}} =(1-\delta)A(k,n)\||f|\|_{p_{k}} = (1-\delta)A(k,n).
\end{equation}
Then by the previous sections, $|f|$ is $o(\delta)$ close in $L^{p_{k}}$ to a Gaussian extremizer, which by an affine change of variables, we can assume to be the standard centered Gaussian $\mathcal{F}^{\dagger}$ on $\mathbb{R}^{n}$: $\| \mathcal{F}^{\dagger} - |f| \|_{p_{k}} = o(\delta)$ and $\|\mathcal{F}^{\dagger} \|_{p_{k}} = 1$. Now \eqref{modulusbd} and Gowers-Host-Kra norm inequality imply,
\begin{equation*} 
    A(k,n)\| |f|\|_{p_{k}} = A(k,n) \geq \| |f| \|_{U^{k}} \geq \| f\|_{U^{k}} \geq (1-\delta)A(k,n)
\end{equation*}
which entails $\| f\|_{U^{k}}\geq (1-o(\delta))\| |f|\|_{U^{k}}$. Then since $\| f\|_{U^{k}} > 0$:
\begin{align} 
    \nonumber \| f\|_{U^{k}}^{2^{k}} = Re(\|f\|_{U^{k}}^{2^{k}}) & \\ 
    \nonumber &= \int_{\mathbb{R}^{(k+1)n}} \prod_{\alpha\in\{0,1\}^{k}} |f|(x+\alpha \cdot \vec{h})Re(\prod_{\alpha\in\{0,1\}^{k}} \mathcal{C}^{\omega_{\alpha}} a(x+\alpha\cdot\vec{h})) \, dxd\vec{h} \\ 
    \label{Ref} &\geq (1-o(\delta)) \| |f| \|_{U^{k}}^{2^{k}} = (1-o(\delta)) \int_{\mathbb{R}^{(k+1)n}} \prod_{\alpha\in\{0,1\}^{k}} |f|(x+\alpha \cdot \vec{h}) \, dxd\vec{h}.
\end{align}
On the other hand, from the Gowers product inequality,
\begin{multline} \label{Fdaggerfmod}
    \bigg | \int_{\mathbb{R}^{(k+1)n}} \prod_{\alpha\in\{0,1\}^{k}} \mathcal{F}^{\dagger}(x+\alpha\cdot\vec{h}) \, dx d\vec{h} - \int_{\mathbb{R}^{(k+1)n}} \prod_{\alpha\in\{0,1\}^{k}} |f|(x+\alpha\cdot\vec{h}) \, dx d\vec{h} \bigg | \\ 
    \leq A(k,n)^{2^{k}}\|\mathcal{F}^{\dagger}-|f|\|_{p_{k}}^{2^{k}} = o(\delta)
\end{multline}
and similarly, since $|a|=1$ on $\mathbb{R}^{n}$,
\begin{multline} \label{FdaggerRe1}
    \bigg | \int_{\mathbb{R}^{(k+1)n}} \prod_{\alpha\in\{0,1\}^{k}} \mathcal{F}^{\dagger}(x+\alpha\cdot\vec{h})Re(\prod_{\alpha\in\{0,1\}^{k}} \mathcal{C}^{\omega_{\alpha}} a(x+\alpha\cdot\vec{h})) \, dx d\vec{h} \\ 
    -  \int_{\mathbb{R}^{(k+1)n}} \prod_{\alpha\in\{0,1\}^{k}} |f|(x+\alpha\cdot\vec{h})Re(\prod_{\alpha\in\{0,1\}^{k}} \mathcal{C}^{\omega_{\alpha}} a(x+\alpha\cdot\vec{h})) \, dx d\vec{h} \bigg | = o(\delta).
\end{multline}
Then \eqref{Fdaggerfmod} and \eqref{FdaggerRe1} allow us to replace $|f|$ with $\mathcal{F}^{\dagger}$ in \eqref{Ref}:
\begin{multline} \label{FdaggerRe2}
    \int_{\mathbb{R}^{(k+1)n}} \prod_{\alpha\in\{0,1\}^{k}} \mathcal{F}^{\dagger}(x+\alpha\cdot\vec{h}) Re(\prod_{\alpha\in\{0,1\}^{k}} \mathcal{C}^{\omega_{\alpha}} a(x+\alpha\cdot\vec{h})) \, dx d\vec{h} \\ 
    \geq (1-o(\delta)) \int_{\mathbb{R}^{(k+1)n}}\prod_{\alpha\in\{0,1\}^{k}} \mathcal{F}^{\dagger}(x+\alpha\cdot\vec{h})\, dx d\vec{h}.
\end{multline}
Let $B_{R}\subset\mathbb{R}^{n}$ denote a centered ball of radius $R$. We can find a sufficiently large positive $R$, so that $\int_{\mathbb{R}^{n}\setminus B_{R}} (\mathcal{F}^{\dagger})^{p_{k}} = o(\delta)$, and $\int_{\mathbb{R}^{n}\setminus B_{R}} |f|^{p_{k}} = o(\delta)$, since $\|\mathcal{F}^{\dagger} -|f|\|_{p_{k}} = o(\delta)$. Moreover, we can select this $R$ so that the following properties are also satisfied:
\begin{align}
    \label{FdaggerBR1} \int_{\mathbb{R}^{(k+1)n}\setminus \tilde{B}_{R}^{k+1}} \prod_{\alpha\in\{0,1\}^{k}} \mathcal{F}^{\dagger}(x+\alpha\cdot\vec{h}) \, dx d\vec{h} &= o(\delta)\\
    \label{FdaggerBR2} \int_{\mathbb{R}^{(k+1)n}\setminus \tilde{B}_{R}^{k+1}} \prod_{\alpha\in\{0,1\}^{k}} \mathcal{F}^{\dagger}(x+\alpha\cdot\vec{h}) Re(\prod_{\alpha\in\{0,1\}^{k}} \mathcal{C}^{\omega_{\alpha}} a(x+\alpha\cdot\vec{h})) \, dx d\vec{h} &= o(\delta).
\end{align}
Then \eqref{FdaggerRe2}, \eqref{FdaggerBR1} and \eqref{FdaggerBR2} allow us to reduce our analysis within a bounded region:
\begin{multline} \label{FBR}
    \int_{\tilde{B}_{R}^{k+1}}\prod_{\alpha\in\{0,1\}^{k}} \mathcal{F}^{\dagger}(x+\alpha\cdot\vec{h}) Re(\prod_{\alpha\in\{0,1\}^{k}} \mathcal{C}^{\omega_{\alpha}} a(x+\alpha\cdot\vec{h})) \, dx d\vec{h} \\ 
    \geq (1-o(\delta)) \int_{\tilde{B}_{R}^{k+1}}\prod_{\alpha\in\{0,1\}^{k}} \mathcal{F}^{\dagger}(x+\alpha\cdot\vec{h})\, dx d\vec{h}.
\end{multline}
Since $|a| = 1$ on $\mathbb{R}^{n}$, it follows from \eqref{FBR} that:
\begin{equation*} 
    \bigg | Re(\prod_{\alpha\in\{0,1\}^{k}} \mathcal{C}^{\omega_{\alpha}} a(x+\alpha\cdot\vec{h})) -1 \bigg | = o(\delta)
\end{equation*}
or, 
\begin{equation*} 
    \prod_{\alpha\in\{0,1\}^{k}} \mathcal{C}^{\omega_{\alpha}} a(x+\alpha\cdot\vec{h}) = 1 + o(\delta)
\end{equation*}
for $(x,\vec{h})$ in a subset of $\tilde{B}_{R}^{k+1}$ of measure $(1-o(\delta))\mathcal{L}(\tilde{B}_{R}^{k+1})$. Since $a(x) = e^{i2\pi q(x)}$, this last display is simply
\begin{equation} \label{angle}
    e^{i2\pi \Delta_{h_{k}}\text{ ... }\Delta_{h_1}q(x)} = 1 + o(\delta).
\end{equation}
Here $\Delta_{h}f(x) = f(x+h)-f(x)$.

\begin{theorem} \label{thm:angle}
Let $n, k\geq 1$. There exists $K > 0$ with the following property.\\
Let $B\subset\mathbb{R}^{n}$ be a centered ball of positive radius and $\psi:B\to\mathbb{R}/\mathbb{Z}$ be a measurable function. Let $\eta,\tau > 0$ be small numbers. Suppose
\begin{equation} \label{thmangle}
    \mathcal{L}(\{(x,\vec{h})\in\tilde{B}^{k+1}: |e^{i2\pi \Delta_{h_1}\text{ ... }\Delta_{h_{k}}\psi(x)} -1| >\tau\}) <\eta\mathcal{L}(\tilde{B}^{k+1}).
\end{equation}
Then there exist a polynomial $P:\mathbb{R}^{n}\to\mathbb{R}$ of degree at most $k-1$ and a positive function $\rho$ satisfying $\lim_{\tau\to 0} \rho(\tau) = 0$ such that
\begin{equation*}
    \mathcal{L}(\{x\in B: |e^{i2\pi\psi(x)}e^{-i2\pi P(x)} - 1| > \rho(\tau)\}) < K\eta\mathcal{L}(B). 
\end{equation*}    
\end{theorem}

\subsection{Proof of Theorem \ref{thm:angle}}

Let $\tilde{B}^2 =\{(x,h_1)\in\mathbb{R}^{2n}: x\in B, h_1\in B\text{ and } x+h_1\in B\}$. When $k=1$, \eqref{thmangle} becomes
\begin{equation} \label{psiangle}
    e^{i2\pi(\psi(x+h_1)-\psi(x))}=1+O(\tau)
\end{equation} 
which happens for a subset of $\tilde{B}^2$ of measure $(1-\eta)\mathcal{L}(\tilde{B}^2)$. Then there exists $c\in B$ such that $|c|=o(\eta)$ and that $c$ satisfies \eqref{psiangle} for $h_1\in B$ outside a subset of measure at most $O(\eta)\mathcal{L}(B)$. That means
\begin{equation*} 
    e^{i2\pi\psi(x)}=(1+O(\tau))e^{i2\pi\psi(c)}
\end{equation*}
for $x\in B$ outside a subset of measure at most $O(\eta)\mathcal{L}(B)$. Hence we can take the constant polynomial $P\equiv\psi(c)$. Assuming the conclusion is true for the case $k-1$, we now prove it for the case of index $k$. \\ 

\noindent
Rewrite \eqref{thmangle} as
\begin{equation} \label{tauangle}
    1 +O(\tau) = e^{i2\pi \Delta_{h_1}\text{ ... }\Delta_{h_{k}}\psi(x)} = e^{i2\pi\Delta_{h_1}\text{ ... }\Delta_{h_{k-1}}(\Delta_{h_{k}}\psi(x))}.
\end{equation}
This holds for $(x,\vec{h})\in\tilde{B}^{k+1}$ outside a subset of measure at most $\eta\mathcal{L}(B)$. Applying the induction hypothesis to $\Delta_{h_{k}}\psi(x)$ in \eqref{tauangle}, we conclude that there exists a polynomial in $x$, $P_{h_{k}}(x) = \sum_{j=0}^{k-2} \sum_{|\gamma|=j} a_{\gamma}(h_{k})x^{\gamma}$ satisfying, for $(x,h_{k})\in\tilde{B}^2$ outside a subset of measure at most $O(\eta)\mathcal{L}(\tilde{B}^2)$,
\begin{equation} \label{polyangle}
    e^{i2\pi\Delta_{h_{k}}\psi(x)} e^{-i2\pi P_{h_{k}}(x)}=1+ o(\tau).
\end{equation}
The coefficient functions $a_{\gamma}$ can be selected to be measurable functions in terms of $h_{k}$ (in fact, these coefficient functions can be selected to be locally piecewise constant, in a manner that is described in Remark \ref{rem:decomp}). To resolve \eqref{polyangle}, we prove a sub-claim:\\

\begin{claim} \label{cl:1}
Let $\psi, B, \eta, \tau$ be as above. Let $P_{t}(x) = \sum_{j=0}^{k-1} \sum_{|\gamma|=j} a_{\gamma}(t)x^{\gamma}$, with $a_{\gamma}:\mathbb{R}^{n}\to\mathbb{R}$ being measurable functions. Suppose for $(x,t)\in\tilde{B}^2$ outside a set of measure at most $O(\eta)\mathcal{L}(\tilde{B}^2)$ the following occurs:
\begin{equation} \label{claim1}
    e^{i2\pi\Delta_{t}\psi(x)} e^{-i2\pi P_{t}(x)} = 1 +O(\tau).
\end{equation}
Then there exists a polynomial $Q$ of degree at most $k$ such that, for $x\in B$ outside a subset of measure at most $O(\eta)\mathcal{L}(B)$,
\begin{equation} \label{claim1ans}
    e^{i2\pi\psi(x)}e^{-i2\pi Q(x)}=1+o(\tau).
\end{equation}
\end{claim}

\begin{proof}[Proof of Claim \ref{cl:1}]
We again use induction. If $k=1$ then \eqref{claim1} becomes
\begin{equation} \label{claim1.1}
    e^{i2\pi\psi(x+t)}e^{-i2\pi\psi(x)}e^{-i2\pi a_0(t)} =1+O(\tau).
\end{equation}
Suppose \eqref{claim1.1} holds for $(x,t)\in\tilde{B}^2$ outside a subset of measure at most $\eta\mathcal{L}(\tilde{B}^2)$, we borrow the following result:

\begin{proposition} \label{prop:christnear} \cite{christ2019near} 
Let $n\geq 1$. There exists a positive constant $K = K(n)> 0$ with the following property. Let $B$ be a centered ball of positive, finite radius. Let $\delta > 0$ and $\eta\in (0,1/2]$. Let $f_1, f_2, f_3: 2B \to\mathbb{C}$ be measurable functions that vanish only on sets of Lebesgue measure zero. Suppose that
\begin{equation*} 
    \mathcal{L}(\{(x,y) \in B^2: |f_1(x)f_2(y)f_3^{-1}(x+y) - 1 | > \eta\} ) < \delta (\mathcal{L}(B))^2.
\end{equation*}
Then for each index $j$ there exists a real-linear function $L_{j}: \mathbb{R}^{n}\to\mathbb{C}$ such that
\begin{equation*} 
    \mathcal{L}(\{ x\in B: |f_{j}(x)e^{-L_{j}(x)}| > K\eta^{1/K}\}) < K\delta\mathcal{L}(B).
\end{equation*}  
\end{proposition}

\noindent 
Applying Proposition \ref{prop:christnear} verbatim, there exist a constant $K = K(n)$, an affine function $L:\mathbb{R}^{n}\to\mathbb{R}$ and a subset of $B$ of measure at most $K\eta\mathcal{L}(B)$ outside of which, $|e^{i2\pi\psi(x)}e^{-i2\pi L(x)}| < K\tau^{1/K}$. Hence the claim holds for the base case $k=1$. Assume the claim is true for the case $k-1$, that is, if \eqref{claim1} holds for a polynomial $P_{t}$ of degree at most $k-1$ in $x$, for $(x,t)\in\tilde{B}^2$ outside a subset of measure at most $O(\eta)\mathcal{L}(\tilde{B}^2)$, then \eqref{claim1ans} follows, with a polynomial $Q$ of degree at most $k$, in terms of $x$. \\

\noindent
See also \cite{christ2021young} for the following argument. Suppose now that $P_{t}$ in \eqref{claim1} is a polynomial of degree at most $k$ in $x$. The leading terms in $x$ in $P_{t}$ are $\sum_{|\gamma|=k} a_{\gamma}(t)x^{\gamma}$. Suppose also for the moment that the leading coefficient functions, for $|\gamma|=k$, are, $a_{\gamma}(t) = \langle A_{\gamma}, t\rangle_{\mathbb{R}^{n}}$, for some elements $A_{\gamma}\in\mathbb{R}^{n}$:
\begin{equation} \label{poly}
    \sum_{|\gamma|=k} a_{\gamma}(t)x^{\gamma} = \sum_{|\gamma|=k} \langle A_{\gamma}, t\rangle_{\mathbb{R}^{n}} x^{\gamma}.
\end{equation}
Define:
\begin{equation} \label{qpoly}
    q(x) = \sum_{|\beta|=k+1} \left (\sum_{i=1}^{n} \beta_{i}^{-1} (A_{\beta})_{i} \right ) x^{\beta}.
\end{equation}
The relations between indices $\gamma$ in \eqref{poly} and $\beta$ in \eqref{qpoly}, and correspondingly, $A_{\gamma}$ and $A_{\beta}$, are as follows. Each $\beta$ has the form $\beta = (\beta_1,\text{ ... },\beta_{i}, \text{ ... }, \beta_{n}) = (\gamma_1,\text{ ... }, \gamma_{i}+1,\text{ ... }, \gamma_{n})$, for some $i\in\{1,\text{ ... }, n\}$, and for each $\beta$ that is arised from a $\gamma$ in this way, $A_{\beta} = A_{\gamma}$. Note that $q$ is a polynomial of degree at most $k+1$, and the leading $x$-terms of $\Delta_{t}q$ are,
\begin{equation} \label{polyleading}
    \sum_{|\beta|=k+1} \left ( \sum_{i=1}^{n} (A_{\beta})_{i} t_{i} x_{i}^{\beta_{i}-1}\prod_{j\not=i}x_{j}^{\beta_{j}}\right ) = \sum_{|\gamma|=k} \langle A_{\gamma},t\rangle_{\mathbb{R}^{n}} x^{\gamma}
\end{equation} 
which is to say, $\Delta_{t}q$ is a polynomial of degree at most $k$ in terms of $x$. Consider $\Psi(x) = \psi(x)-q(x)$. Then \eqref{claim1} yields, for $(x,t)\in\tilde{B}^2$ outside a subset of measure at most $O(\eta)\mathcal{L}(\tilde{B}^2)$,
\begin{equation} \label{St}
    e^{i2\pi\Delta_{t}\Psi(x)} = e^{i2\pi\Delta_{t}(\psi - q)(x)} = (1+O(\tau))e^{i2\pi P_{t}(x)}e^{-i2\pi\Delta_{t}q(x)} = (1+O(\tau)) e^{i2\pi S_{t}(x)}.
\end{equation}
Note that \eqref{poly} and \eqref{polyleading} imply that in $S_{t}$ in \eqref{St} takes the following form: 
\begin{equation*} S_{t}(x) = \sum_{j=0}^{k-1}\sum_{|\mu|=j}c_{\mu}(t)x^{\mu}.\end{equation*}
Hence by the induction hypothesis, we conclude from \eqref{St} that there exists a polynomial $Q$ of degree at most $k$ such that $e^{i2\pi\Psi(x)} = (1+o(\tau))e^{i2\pi Q(x)}$, and in turn, $e^{i2\pi\psi(x)} = (1+o(\tau))e^{i2\pi (Q(x)+q(x))}$, for $x\in B$ outside a subset of measure at most $O(\eta)\mathcal{L}(B)$. It's obvious that $Q+q$ is a polynomial of degree at most $k+1$; we now close the induction loop.    
\end{proof}

\noindent
In the above proof of Claim \ref{cl:1}, we need to make an assumption that if $|\gamma|=k$, $a_{\gamma}(t) = \langle A_{\gamma},t\rangle_{\mathbb{R}^{n}}$ for some $A_{\gamma}\in\mathbb{R}^{n}$, for $t\in B$ outside a subset of measure at most $O(\eta)\mathcal{L}(B)$. Now we show that such condition must indeed occur:

\begin{claim} \label{cl:2}
Let $\psi, \eta, \tau$ be as above. Given $P_{t}(x) = \sum_{j=0}^{k} \sum_{|\gamma|=j} a_{\gamma}(t)x^{\gamma}$, with $a_{\gamma}:\mathbb{R}^{n}\to\mathbb{R}$ being measurable functions. Suppose for $(x,t)\in\tilde{B}^2$ outside a subset of measure at most $O(\eta)\mathcal{L}(\tilde{B}^2)$, the following happens:
\begin{equation} \label{claim2}
    e^{i2\pi\Delta_{t}\psi(x)} = (1+O(\tau))e^{i2\pi P_{t}(x)}.
\end{equation}
Then, for $|\gamma|=k$, there exist $A_{\gamma}\in\mathbb{R}^{n}$ such that $a_{\gamma}(t) = \langle A_{\gamma},t\rangle_{\mathbb{R}^{n}}$, for $t\in B$ outside of a subset of measure at most $O(\eta)\mathcal{L}(B)$.    
\end{claim}

\begin{proof}[Proof of Claim \ref{cl:2}]
\eqref{claim2} gives us,
\begin{equation} \label{polyequiv}
    \psi(x+t)-\psi(x) \equiv\sum_{j=0}^{k} \sum_{|\gamma|=j} a_{\gamma}(t)x^{\gamma} +o(\tau).
\end{equation}
By "$u\equiv v$", we mean $u - v \in\mathbb{Z}$. Let $\mathcal{A}^{x}=\{t\in B: (x,t)\in\tilde{B}^2 \text{ and } \eqref{polyequiv} \text{ apply}.\}$, for each $x\in B$. Let $\mathcal{A}$ denote the set of all $x\in (1/2)B$ such that $\mathcal{L}(\mathcal{A}^{x}) > (1-O(\eta))\mathcal{L}(B)$. Then there exists $c>0$ such that, for all sufficiently small $\eta$, $\mathcal{L}(\mathcal{A}) > c$. Then uniformly for $x\in\mathcal{A}$, as $t$ varies in $\mathcal{A}^{x}$, these values $x+t\in B$ occupy a subset of $B$ of measure at least $(1-O(\eta))\mathcal{L}(B)$. Fix $x\in\mathcal{A}$, this set of values $x+t$ has a nonempty intersection with $\mathcal{A}$. Fix one such value $x+t\in\mathcal{A}$. Then $\mathcal{L}(\mathcal{A}^{x+t}) > (1-O(\eta))\mathcal{L}(B)$. That means, the set of values $t+s$, as $s$ varies in $\mathcal{A}^{x+t}$, occupy a subset of $B$ of measure at least $(1-O(\eta))\mathcal{L}(B)$. This set of values $t+s$ will then have a nonempty intersection with $\mathcal{A}^{x}$. In other words, we can select $x, t, s$ so that the following argument is applicable:\\

\noindent
Consider $\psi(x+t+s) -\psi(x+t)$. From \eqref{polyequiv} we can write the said difference in two ways. Firstly:
\begin{align} 
    \nonumber \psi((x+t)+s)-\psi(x+t) &\equiv\sum_{j=0}^{k} \sum_{|\gamma|=j} a_{\gamma}(s)(x+t)^{\gamma} +o(\tau) \\
    \label{psidiff1} &= \sum_{|\gamma|=k} a_{\gamma}(s)(x+t)^{\gamma} + O(x^{k-1}) +o(\tau) =  \sum_{|\gamma|=k} a_{\gamma}(s)x^{\gamma} + O(x^{k-1}) +o(\tau).
\end{align}
By "$O(x^{k-1})$" we mean a linear combination of monomials in terms of $x$ of degree at most $k-1$. Secondly:
\begin{multline} \label{psidiff2}
    \psi(x+t+s)-\psi(x+t) = [\psi(x+(t+s)) - \psi(x)] - [\psi(x+t)-\psi(x)] \\ \equiv \sum_{|\gamma|=k} [a_{\gamma}(t+s)-a_{\gamma}(t)]x^{\gamma} + O(x^{k-1}) +o(\tau).
\end{multline}
Since $\mathcal{L}(\mathcal{A})>c$, a comparison between \eqref{psidiff1} and \eqref{psidiff2} shows that, for $|\gamma|=k$:
\begin{equation} \label{aequal}
    a_{\gamma}(t+s) = a_{\gamma}(t)+ a_{\gamma}(s).
\end{equation}
By the argument presented above, \eqref{aequal} is satisfied for $t,s\in B$ except for a subset of measure at most $O(\eta)\mathcal{L}(B)$. Let $\{\rho_{i}\}_{i}$ be a sequence tending to zero. Then from \eqref{aequal}, for a.e $(t,s)\in B\times B$ outside a subset of measure at most $O(\eta)\mathcal{L}(B\times B)$ and for every $\rho_{i}$, we have, $ |a_{\gamma}(t+s) - a_{\gamma}(t)- a_{\gamma}(s)| <\rho_{i}$. By Proposition \ref{prop:christm} and its proof given in \cite{christ2019near}, there exist affine functions $L_{\rho_{i}}(t) = \langle A_{\rho_{i}},t\rangle_{\mathbb{R}^{n}} + b_{\rho_{i}}$, $A_{\rho_{i}}\in\mathbb{R}^{n}, b_{\rho_{i}}\in\mathbb{R}$, and a subset $U\subset B$, such that $a_{\gamma}(t) - L_{\rho_{i}}(t) = O(\rho_{i})$ for a.e $t\in U$. The relative complement of $U$ has measure at most $O(\eta)\mathcal{L}(B)$, and the implicit constants in the notations $O(\rho_{i}), O(\eta)\mathcal{L}(B)$ are independent of $i$. This implies in particular that the $L_{\rho_{i}}$ are uniformly bounded and converges pointwise in $U$. Since $L_{\rho_{i}}$ are affine functions on $\mathbb{R}^{n}$, this in turn implies, by passing to a subsequence, $A_{\rho_{i}}\to A\in\mathbb{R}^{n}$ and $b_{\rho_{i}}\to b\in\mathbb{R}$ as $i\to\infty$ \cite{rudin1976principles}. Hence $a_{\gamma}(t) = \langle A, t\rangle_{\mathbb{R}^{n}} + b$ for $t\in B$ outside of a subset of measure at most $O(\eta)\mathcal{L}(B)$. But then \eqref{aequal} gives $b=0$, hence the proof of Claim \ref{cl:2} is now complete.    
\end{proof}

\noindent
We have also finished the proof of Theorem \ref{thm:angle}. $\qed$

\subsection{Production of an extremizer}

We combine the results the previous sections, \eqref{angle} and Theorem \ref{thm:angle} to conclude that there exist a centered ball $B_{R}$, with $R$ sufficiently large, so that
\begin{equation} \label{fout}
    \int_{\mathbb{R}^{n}\setminus B_{R}} |f|^{p_{k}} = o(\delta)
\end{equation}
and a polynomial $P$ of degree at most $k-1$ so that,
\begin{equation} \label{polyphase}
    \|(e^{i2\pi P}-a)\cdot\chi_{B_{R}}\|_{p_{k}} = o(\delta).
\end{equation}
Let $\mathcal{G} = \mathcal{F}^{\dagger}e^{i2\pi P}$, with $\mathcal{F}^{\dagger}$ is as in Section \ref{sec:compprep}. Then $\mathcal{G}$ is an extremizer of \eqref{GHKineq}. Observe:
\begin{equation} \label{Gfdecomp}
    \| \mathcal{G}-f \|_{p_{k}}^{p_{k}} = \int_{\mathbb{R}^{n}\setminus B_{R}} |\mathcal{G}-f|^{p_{k}} + \int_{B_{R}} |\mathcal{G} - f|^{p_{k}}.
\end{equation}
By \eqref{fout}, the facts that $f(x) = |f|(x)a(x)$ and $|e^{i2\pi P} - a| \leq C$ uniformly on $\mathbb{R}^{n}$:
\begin{equation*} \int_{\mathbb{R}^{n}\setminus B_{R}} ||f|e^{i2\pi P} - |f|a|^{p_{k}} = o(\delta),\end{equation*}
which together with the fact that $\| \mathcal{F}^{\dagger}-|f|\|_{p_{k}}^{p_{k}} = o(\delta)$ allows us to conclude the following for the first term in \eqref{Gfdecomp}:
\begin{multline*} 
    \int_{\mathbb{R}^{n}\setminus B_{R}} |\mathcal{F}^{\dagger}e^{i2\pi P} - |f|a|^{p_{k}} \leq C(k)\int_{\mathbb{R}^{n}\setminus B_{R}} |\mathcal{F}^{\dagger}e^{i2\pi P} - |f|e^{i2\pi P}|^{p_{k}} \\ + C(k)\int_{\mathbb{R}^{n}\setminus B_{R}} ||f|e^{i2\pi P} - |f|a|^{p_{k}} \leq C(k)\|\mathcal{F}^{\dagger}-|f|\|_{p_{k}}^{p_{k}} + o(\delta) = o(\delta). 
\end{multline*}
In the same spirit, we split the second term in \eqref{Gfdecomp}:
\begin{equation} \label{GfBR}
    \int_{B_{R}} |\mathcal{G} - f|^{p_{k}} \leq C(k)\int_{B_{R}}  |\mathcal{F}^{\dagger}e^{i2\pi P} - \mathcal{F}^{\dagger}a|^{p_{k}} + C(k)\int_{B_{R}} |\mathcal{F}^{\dagger}a - |f|a|^{p_{k}}. 
\end{equation}
The second term in \eqref{GfBR} is at most $o(\delta)$ in absolute value, again by $\| \mathcal{F}^{\dagger}-|f|\|_{p_{k}}^{p_{k}} = o(\delta)$ and $|a| = 1$ on $\mathbb{R}^{n}$. It follows from \eqref{polyphase} that the contribution of first term in \eqref{GfBR} is also at most $o(\delta)$ in absolute value: 
\begin{equation*} 
    \int_{B_{R}} |\mathcal{F}^{\dagger}e^{i2\pi P} - \mathcal{F}^{\dagger}a|^{p_{k}} \leq C(k)\int_{B_{R}} |e^{i2\pi P}-a|^{p_{k}} = o(\delta).
\end{equation*}
Hence we conclude that $\| \mathcal{G}-f\|_{p_{k}} = o(\delta)$. \\

\noindent
We have now obtained the conclusion for Theorem \ref{thm:main} for a general measurable near extremizer $f:\mathbb{R}^{n}\to\mathbb{C}$, $n\geq 1$, of the $k$th Gowers-Host-Kra norm inequality for $k\geq 2$. $\qed$

\end{document}